\documentclass[12pt]{article}
\usepackage[utf8]{inputenc}
\usepackage[T1]{fontenc}

\usepackage[]{amsfonts}
\usepackage{amssymb}
\usepackage{amsmath}

\def\couleur(#1 #2 #3)
	{
\special{ps::[end]}}

\def\underset#1#2{\mathrel{\mathop{\kern0pt #2}\limits_{#1}}}

\def\overset#1#2{\mathrel{\mathop{\kern0pt #2}\limits^{#1}}}

\def\bx#1{\setbox1=\hbox{\kern3pt{#1}\kern3pt}			
 \dimen1=\ht1 \advance\dimen1 by 3pt \dimen2=\dp1 \advance\dimen2 by 3pt
 \setbox1=\hbox{\vrule height\dimen1 depth\dimen2\box1\vrule}%
 \setbox1=\vbox{\hrule\box1\hrule}%
 \advance\dimen1 by .4pt \ht1=\dimen1
 \advance\dimen2 by .4pt \dp1=\dimen2 \box1\relax}

\def\wbb#1{\kern#1em}
\def\vci{\vrule  width.02em height1.47ex depth-.0ex}		
\def\11{{\rm\wbb{.2}\vci\wbb{-.37}1}}

\newtheorem{Thrm}{Theorem}[section]
\newtheorem{Lmm}[Thrm]{Lemma}
\newtheorem{Dfnt}[Thrm]{Definition}
\newtheorem{Prps}[Thrm]{Proposition}
\newtheorem{Crll}[Thrm]{Corollary}
\newtheorem{Rmrq}[Thrm]{Remark}

\begin{document}

\title{Interpolating sequences and Carleson measures in the Hardy-Sobolev spaces of the ball in  $\displaystyle {\mathbb{C}}^{n}.$ }

\author{E. Amar}
\maketitle
 \ \par 
\renewcommand{\abstractname}{Abstract}

\begin{abstract}
In this work we study Hardy Sobolev spaces in the ball of  $\displaystyle
 {\mathbb{C}}^{n}$  with respect to interpolating sequences and
 Carleson measures.\ \par 
\quad  We compare them with the classical Hardy spaces of the ball
 and we stress analogies and differences.\ \par 

\tableofcontents
\end{abstract}

\section{Introduction.}
 We shall work with the Hardy-Sobolev spaces  $\displaystyle
 H_{s}^{p}.$  For  $\displaystyle 1\leq p<\infty $  and  $s\in
 {\mathbb{R}},\ H_{s}^{p}$  is the space of holomorphic functions
 in the unit ball  $\displaystyle {\mathbb{B}}$  in  $\displaystyle
 {\mathbb{C}}^{n}$  such that the following expression is finite\ \par 
\quad \quad \quad \quad \quad 	 $\displaystyle \ {\left\Vert{f}\right\Vert}_{s,p}^{p}:=\sup
 _{r<1}\int_{\partial {\mathbb{B}}}{\left\vert{(I+R)^{s}f(rz)}\right\vert
 ^{p}d\sigma (z)},$ \ \par 
where  $I$  is the identity,  $\displaystyle d\sigma $  is the
 Lebesgue measure on  $\displaystyle \partial {\mathbb{B}}$ 
 and  $R$  is the radial derivative\ \par 
\quad \quad \quad \quad \quad 	 $\displaystyle Rf(z)=\sum_{j=1}^{n}{z_{j}\frac{\partial f}{\partial
 z_{j}}(z)}.$ \ \par 
For  $s\in {\mathbb{N}},$  this norm is equivalent to\ \par 
\quad \quad \quad \quad 	 $\displaystyle \ {\left\Vert{f}\right\Vert}_{H_{s}^{p}}^{p}=\max
 _{0\leq j\leq s}\int_{\partial {\mathbb{B}}}{\left\vert{R^{j}f(z)}\right\vert
 ^{p}d\sigma (z)}.$ \ \par 
This means that  $\displaystyle R^{j}f\in H^{p}({\mathbb{B}}),\
 j=0,...,\ s.$ \ \par 
\quad  We shall prove estimates only in the case  $s\in {\mathbb{N}}$
  but the spaces  $\displaystyle H_{s}^{p}$  form an interpolating
 scale with respect to the parameter  $\displaystyle s,$  see
 in section~\ref{1HS0}, hence, in some cases, this allows to
 extend the results to the case  $s\in {\mathbb{R}}_{+}.$ \ \par 
\ \par 
\quad  If  $\displaystyle sp>n$  the functions in  $\displaystyle H_{s}^{p}$
  are continuous up to the boundary  $\displaystyle \partial
 {\mathbb{B}}$  hence the results we are interested in are essentially
 trivial, so we shall restrict ourselves to the case  $\displaystyle
 s\leq n/p.$ \ \par 
\ \par 
\quad \quad  	For  $\displaystyle s=0$  the Hardy Sobolev spaces  $\displaystyle
 H_{0}^{p}$  are the classical Hardy spaces  $\displaystyle H^{p}({\mathbb{B}})$
  of the unit ball  $\displaystyle {\mathbb{B}}$  and a natural
 question is to study what remains true from classical Hardy
 spaces  $\displaystyle H^{p}({\mathbb{B}})$  to Hardy Sobolev
  $\displaystyle H_{s}^{p}.$ \ \par 
\quad  An important notion is that of Carleson measure.\ \par 

\begin{Dfnt}
 The {\bf measure}  $\mu $  in  $\displaystyle {\mathbb{B}}$
  is {\bf Carleson}  for  $\displaystyle H_{s}^{p},\ \mu \in
 C_{s,p},$  if we have the embedding\par 
\quad \quad \quad \quad \quad 	 $\displaystyle \forall f\in H_{s}^{p},\ \int_{{\mathbb{B}}}{\left\vert{f}\right\vert
 ^{p}d\mu \leq C{\left\Vert{f}\right\Vert}_{H_{s}^{p}}^{p}.}$  
\end{Dfnt}
\quad  Carleson measures where introduced by Carleson~\cite{CarlInt58}
 in his work on interpolating sequences.\ \par 
\quad  We have the following table concerning the known results about
 Carleson measures :\ \par 
\ \par 
{\begin{tabular}{|c|c|c|c|c|c|c|c|c|c|c|c|c|c|c|c|c|c|c|c}
\hline  $H^{p}({\mathbb{D}})$  &   $\displaystyle H^{p}({\mathbb{B}})=H_{0}^{p}({\mathbb{B}})$
  &   $\displaystyle H_{s}^{p}({\mathbb{B}})$  \cr
\hline {\begin{tabular}{cccccccccccccccccccccc}
Characterized \cr
geometrically \cr
by L. Carleson~\cite{CarlInt58}\cr
\end{tabular}}
 & {\begin{tabular}{cccccccccccccccccccccc}
Characterized \cr
geometrically \cr
by L. H\"ormander~\cite{HormPSH67}\cr
\end{tabular}}
 & {\begin{tabular}{cccccccccccccccccccccc}
Studied by C. Cascante \&  \cr
J. Ortega~\cite{CascOrte95} ; characterized \cr
for  $\displaystyle n-1\leq ps\leq n.$  \cr
For  $\displaystyle p=2,$  any  $\displaystyle s$   \cr
characterized by \cr
A. Volberg \&  B. Wick~\cite{VolbWick12}\cr
\end{tabular}}
 \cr
\hline Same for all  $\displaystyle p$  & Same for all  $\displaystyle
 p$  & Depending on  $\displaystyle p$ \cr
\hline 
\end{tabular}}
\ \par 
\ \par 

\begin{Dfnt}
 The {\bf multipliers algebra}  $\displaystyle {\mathcal{M}}_{s}^{p}$
  of  $\displaystyle H_{s}^{p}$  is the algebra of functions
  $m$  on  $\displaystyle {\mathbb{B}}$  such that\par 
\quad \quad \quad \quad \quad 	 $\displaystyle \forall h\in H_{s}^{p},\ mh\in H_{s}^{p}.$ \par 
The norm of a multiplier is its norm as an operator from  $\displaystyle
 H_{s}^{p}$  into  $\displaystyle H_{s}^{p}.$ 
\end{Dfnt}
\quad \quad  	We have the following table of already known results, where
 C.C. means a certain Carleson condition:\ \par 
\ \par 
\ \par 
{\begin{tabular}{|c|c|c|c|c|c|c|c|c|c|c|c|c|c|c|c|c|c|c|c}
\hline   $H^{p}({\mathbb{D}})$   &   $\displaystyle H^{p}({\mathbb{B}})$
   &   $\displaystyle H_{s}^{p}({\mathbb{B}})$   \cr
\hline   ${\mathcal{M}}_{0}^{p}({\mathbb{D}})=H^{\infty }({\mathbb{D}}),\
 \forall p$   &   $\displaystyle {\mathcal{M}}_{0}^{p}({\mathbb{B}})=H^{\infty
 }({\mathbb{B}}),\ \forall p$   & {\begin{tabular}{cccccccccccccccccccccc}
 $\displaystyle {\mathcal{M}}_{s}^{p}=H^{\infty }({\mathbb{B}})\cap C.C.$   \cr
characterized for \cr
 $\displaystyle n-1\leq ps\leq n$   \cr
and for  $\displaystyle p=2$  by Volberg \&  Wick~\cite{VolbWick12} \cr
Depending on  $\displaystyle p$  \cr
\end{tabular}}
\cr
\hline 
\end{tabular}}
\ \par 
\ \par 
\quad  Now we shall deal with sequences of points in the ball  $\displaystyle
 {\mathbb{B}}$  and to state results we shall need several definitions,
 most of them being well known.\ \par 
\quad \quad  	Let  $\displaystyle p'$  be the conjugate exponent for  $\displaystyle
 p,\ \frac{1}{p}+\frac{1}{p'}=1\ ;$  the Hilbert space  $\displaystyle
 H_{s}^{2}$  is equipped with the reproducing kernels :\ \par 
\quad \quad \quad \quad  	\begin{equation}  \forall a\in {\mathbb{B}},\ k_{a}(z)=\frac{1}{(1-\bar
 a\cdot z)^{n-2s}},\ {\left\Vert{k_{a}}\right\Vert}_{H_{s}^{p}}\simeq
 (1-\left\vert{a}\right\vert ^{2})^{s-n/p'}\label{iS20}\end{equation}\ \par 
i.e.  $\displaystyle \forall a\in {\mathbb{B}},\ \forall f\in
 H_{s}^{p},\ f(a)={\left\langle{f,k_{a}}\right\rangle},$  where
  $\displaystyle \ {\left\langle{\cdot ,\cdot }\right\rangle}$
  is the scalar product of the Hilbert space  $\displaystyle
 H_{s}^{2}.$  In the case  $\displaystyle s=n/2$  there is a
 log in  $\displaystyle k_{a}.$ \ \par 

\begin{Dfnt}
 The sequence  $S$  is Carleson in  $\displaystyle H_{s}^{p}({\mathbb{B}}),$
  if the associated measure\par 
\quad \quad \quad \quad \quad 	 $\displaystyle \nu _{S}:=\sum_{a\in S}{{\left\Vert{k_{s,a}}\right\Vert}_{H_{s}^{p'}}^{-p}\delta
 _{a}}$ \par 
is Carleson for  $\displaystyle H_{s}^{p}({\mathbb{B}}).$ \par 

\end{Dfnt}

\begin{Dfnt}
Let  $\displaystyle p>1,$  the {\bf sequence}  $S$  of points
 in  $\displaystyle {\mathbb{B}}$  is {\bf interpolating} in
  $\displaystyle H_{s}^{p}({\mathbb{B}}),$  if there is a  $\displaystyle
 C=C_{p}>0$  such that\par 
\quad \quad \quad 	 $\displaystyle \forall \lambda \in \ell ^{p}(S),\ \exists f\in
 H_{s}^{p}({\mathbb{B}})::\forall a\in S,\ f(a)=\lambda _{a}{\left\Vert{k_{a}}\right\Vert}_{H_{s}^{p'}}=\lambda
 _{a}(1-\left\vert{a}\right\vert ^{2})^{s-n/p},\ {\left\Vert{f}\right\Vert}_{H_{s}^{p}}\leq
 C{\left\Vert{\lambda }\right\Vert}_{p},$ \par 
where  $\displaystyle p'$  is the conjugate exponent for  $\displaystyle
 p,\ \frac{1}{p}+\frac{1}{p'}=1.$ \par 
\quad  If  $\displaystyle p=1,$  we take the limiting case in the above
 definition :  $S$  is IS for  $\displaystyle H_{s}^{1},$  if
 there is a  $\displaystyle C>0$  such that	      	 $\displaystyle
 \forall \lambda \in \ell ^{1}(S),\ \exists f\in H_{s}^{1}({\mathbb{B}})::\forall
 a\in S,\ f(a)=\lambda _{a}(1-\left\vert{a}\right\vert ^{2})^{s-n},\
 {\left\Vert{f}\right\Vert}_{H_{s}^{1}}\leq C{\left\Vert{\lambda
 }\right\Vert}_{\ell ^{1}}.$ 
\end{Dfnt}
\ \par 

\begin{Dfnt}
 The {\bf sequence}  $S$  of points in  $\displaystyle {\mathbb{B}}$
  is {\bf interpolating} in the multipliers algebra  $\displaystyle
 {\mathcal{M}}_{s}^{p}$  of  $\displaystyle H_{s}^{p}({\mathbb{B}})$
  if  there is a  $\displaystyle C>0$  such that\par 
\quad \quad   $\forall \lambda \in \ell ^{\infty }(S),\ \exists m\in {\mathcal{M}}_{s}^{p}::\forall
 a\in S,\ m(a)=\lambda _{a}$ 	and  $\displaystyle \ {\left\Vert{m}\right\Vert}_{{\mathcal{M}}_{s}^{p}}\leq
 C{\left\Vert{\lambda }\right\Vert}_{\infty }.$ \par 

\end{Dfnt}

\begin{Dfnt}
 Let  $S$  be an interpolating sequence in  $\displaystyle {\mathcal{M}}_{s}^{p}\
 ;$  we say that  $S$  has a {\bf bounded linear extension operator,
 BLEO}, if there is a a bounded linear operator  $\displaystyle
 E\ :\ \ell ^{\infty }(S)\rightarrow {\mathcal{M}}_{s}^{p}$ 
 and a  $\displaystyle C>0$  such that\par 
\quad \quad \quad 	 $\displaystyle \forall \lambda \in \ell ^{\infty }(S),\ E(\lambda
 )\in {\mathcal{M}}_{s}^{p},\ {\left\Vert{E(\lambda )}\right\Vert}_{{\mathcal{M}}_{s}^{p}}\leq
 \ C{\left\Vert{\lambda }\right\Vert}_{\infty }\ :\ \forall a\in
 S,\ E(\lambda )(a)=\lambda _{a}.$ 
\end{Dfnt}
\ \par 
\quad \quad  	We have the table of results on interpolating sequences, where
 A.R.S. means Arcozzi, Rochberg and Sawyer~\cite{ArcoRochSaw08}.\ \par 
\ \par 
{\begin{tabular}{|c|c|c|c|c|c|c|c|c|c|c|c|c|c|c|c|c|c|c|c}
\hline   $H^{\infty }({\mathbb{D}})$   &   $\displaystyle H^{\infty
 }({\mathbb{B}})$   &   $\displaystyle {\mathcal{M}}_{s}^{p}({\mathbb{B}})$
   \cr
\hline {\begin{tabular}{cccccccccccccccccccccc}
IS characterized \cr
by L. Carleson\cr
\end{tabular}}
 & No characterisation & {\begin{tabular}{cccccccccccccccccccccc}
Characterized for  $\displaystyle p=2$   \cr
and  $\displaystyle n-1<2s\leq n$   \cr
by A.R.S. and the \cr
Pick property\cr
\end{tabular}}
 \cr
\hline {\begin{tabular}{cccccccccccccccccccccc}
ISM $\displaystyle \Rightarrow $ BLEO \cr
by P. Beurling~\cite{PBeurling62}\cr
\end{tabular}}
 & {\begin{tabular}{cccccccccccccccccccccc}
ISM $\displaystyle \Rightarrow $ BLEO \cr
by A. Bernard~\cite{Bernard71}\cr
\end{tabular}}
 & {\begin{tabular}{cccccccccccccccccccccc}
ISM $\displaystyle \Rightarrow $ BLEO \cr
by E. A. here \cr
\end{tabular}}
\cr
\hline 
\end{tabular}}
\ \par 
\ \par 
\quad  In the case of the classical Hardy spaces  $\displaystyle H^{p}$
 , whose multiplier algebra is  $\displaystyle H^{\infty }({\mathbb{B}}),$
  we know(~\cite{amBerg78} theorem 5, p. 712 and the lines following
 it) that if  $S$  is interpolating for  $\displaystyle H^{\infty
 }({\mathbb{B}})$  then  $S$  is interpolating for  $\displaystyle
 H^{p}({\mathbb{B}})$  ; this is still true in the case of Hardy
 Sobolev spaces.\ \par 

\begin{Thrm}
 Let  $S$  be an interpolating sequence for the multipliers algebra
  ${\mathcal{M}}_{s}^{p}$  of  $\displaystyle H_{s}^{p}({\mathbb{B}})$
  then  $S$  is also an interpolating sequence for  $\displaystyle
 H_{s}^{p},$  with a bounded linear extension operator.
\end{Thrm}
\quad  In the classical case  $\displaystyle s=0$  i.e.  $\displaystyle
 {\mathcal{M}}_{0}^{p}=H^{\infty }({\mathbb{B}}),\ H_{0}^{p}=H^{p}({\mathbb{B}}),$
  N. Varopoulos~\cite{Varo72} proved that  $S$  interpolating
 in  $\displaystyle H^{\infty }({\mathbb{B}})$  implies that
  $S$  is Carleson in  $\displaystyle {\mathbb{B}}$  and P. Thomas~\cite{Thomas87}
 (see also~\cite{AmarWirtBoule07}) proved that  $S$  interpolating
 in  $\displaystyle H^{p}({\mathbb{B}})$  implies that  $S$ 
 is Carleson in  $\displaystyle {\mathbb{B}}.$  The next results
 generalise this fact to  $\displaystyle {\mathcal{M}}_{s}^{p}$
  for  $\displaystyle p\leq 2$  and any real values of  $\displaystyle
 s\in \lbrack 0,n/p\rbrack .$ \ \par 

\begin{Thrm}
 Let  $S$  be an interpolating sequence for  $\displaystyle {\mathcal{M}}_{s}^{p}$
  with  $\displaystyle p\leq 2,$  then  $S$  is Carleson  $\displaystyle
 H_{s}^{p}({\mathbb{B}}).$ 
\end{Thrm}
\quad  And\ \par 

\begin{Crll}
 Let  $S$  be an interpolating sequence for  $\displaystyle H_{s}^{2}$
  with  $\displaystyle n-2s\leq 1,$  i.e.  $\displaystyle {\mathcal{M}}_{s}^{2}$
  is a Pick algebra, then  $S$  is Carleson for  $\displaystyle
 H_{r}^{2},\ \forall r\leq s.$ 
\end{Crll}
\ \par 
\quad  Because  $\displaystyle {\mathcal{M}}_{s}^{2}$  is an operators
 algebra in a Hilbert space, then we know~\cite{AmarThesis77}
 that the union  $S$  of two interpolating sequences in   $\displaystyle
 {\mathcal{M}}_{s}^{2}$  is still interpolating in  $\displaystyle
 {\mathcal{M}}_{s}^{2}$  if  $S$  is separated. This generalises
 a theorem of Varopoulos~\cite{Varo71} done for uniform algebras.\ \par 
\quad  We prove the analogous result but we shall have to use a more
 precise notion of separation, see section~\ref{6IM3}.\ \par 

\begin{Thrm}
 Let  $s\in {\mathbb{N}}\cap \lbrack 0,n/p\rbrack $  and  $\displaystyle
 S_{1}$  and  $\displaystyle S_{2}$  be two completely separated
 interpolating sequences in  $\displaystyle {\mathcal{M}}_{s}^{p}$
  then  $\displaystyle S:=S_{1}\cup S_{2}$  is still an interpolating
 sequence in  $\displaystyle {\mathcal{M}}_{s}^{p}.$ 
\end{Thrm}
\quad  In  $\displaystyle {\mathbb{C}}^{n},\ n\geq 2$  we know~\cite{DenAm78}
 that the union  $S$  of two interpolating sequences in  $\displaystyle
 H^{2}({\mathbb{B}})$  is not in general an interpolating sequence
 even if  $S$  is separated, so the next result is in complete
 opposition to this fact.\ \par 

\begin{Crll}
 If  $\displaystyle {\mathcal{M}}_{s}^{2}$  is a Pick algebra,
 i.e. if  $\displaystyle s\geq \frac{n-1}{2},$  and   $\displaystyle
 S,\ S'$  are two interpolating sequences for  $\displaystyle
 H_{s}^{2}$  such that  $\displaystyle S\cup S'$  is separated
 then  $\displaystyle S\cup S'$  is still interpolating for 
 $\displaystyle H_{s}^{2}.$ 
\end{Crll}
\quad  We shall need the following notion.\ \par 

\begin{Dfnt}
The {\bf sequence}  $S$  of points in  $\displaystyle {\mathbb{B}}$
  is {\bf dual bounded} (or minimal, or weakly interpolating)
 in  $\displaystyle H_{s}^{p}({\mathbb{B}})$  if there is a bounded
 sequence  $\displaystyle \lbrace \rho _{a}\rbrace _{a\in S}\subset
 H_{s}^{p}$  such that\par 
\quad \quad \quad \quad \quad 	 $\displaystyle \forall a,b\in S,\ \rho _{a}(b)=\delta _{ab}{\left\Vert{k_{a}}\right\Vert}_{H_{s}^{p'}}.$
 
\end{Dfnt}
\quad  	Clearly if  $S$  is interpolating for  $\displaystyle H_{s}^{p}$
  then it is dual bounded in  $\displaystyle H_{s}^{p}.$ \ \par 
\quad  This notion characterizes interpolating sequences for the classical
 Hardy spaces in the unit disc  ${\mathbb{D}}\ ;$  the question
 is open for the Hardy spaces  $\displaystyle H^{p}({\mathbb{B}})$
  in the ball in  $\displaystyle {\mathbb{C}}^{n}\geq 2.$  Nevertheless
 we know~\cite{AmarCarlType07} that if  $\displaystyle S\subset
 {\mathbb{B}}$  is dual bounded for  $\displaystyle H^{p}({\mathbb{B}})$
  then it is interpolating for  $\displaystyle H^{q}({\mathbb{B}}),\
 \forall q<p.$ \ \par 
\quad  The next results generalise only partially this result and we
 get an analogous result to theorem 6.1 in~\cite{AmarExtInt06}.\ \par 

\begin{Dfnt}
 We shall say that  $S$  is a  $\displaystyle H_{s}^{p}$  weighted
 interpolating sequence for the weight  $\displaystyle w=\lbrace
 w_{a}\rbrace _{a\in S}$  if \par 
\quad \quad \quad  $\displaystyle \forall \lambda \in \ell ^{p}(S),\ \exists f\in
 H_{s}^{p}::\forall a\in S,\ f(a)=\lambda _{a}w_{a}{\left\Vert{k_{a}}\right\Vert}_{H_{s}^{p'}}.$
 
\end{Dfnt}

\begin{Thrm}
 Let  $S$  be a sequence of points in  $\displaystyle {\mathbb{B}}$
  such that, with  $\displaystyle \ \frac{1}{r}=\frac{1}{p}+\frac{1}{q},$
  and  $\displaystyle p\leq 2,$  \par 
\quad \quad \quad 	 $\displaystyle \bullet $   $S$  is dual bounded in  $\displaystyle
 H_{s}^{p}.$ \par 
\quad \quad \quad  $\displaystyle \bullet $    $S$  is Carleson in  $\displaystyle
 H_{s}^{q}({\mathbb{B}}).$ \par 
Then  $S$  is a  $\displaystyle H_{s}^{r}$  weighted interpolating
 sequence for the weight  $\displaystyle \lbrace (1-\left\vert{a}\right\vert
 ^{2})^{s}\rbrace _{a\in S}$  with the bounded linear extension property.
\end{Thrm}
\quad  This work was exposed in Oberwolfach workshop "Hilbert Modules
 and Complex Geometry " in April 2014, and also in the conference
 in honor of A. Bonami in June 2014 in Orleans, France. This
 is an improved version of these talks.\ \par 
This work is presented the following way.\ \par 
\ \par 
\quad  In the next section we study the basics of Hardy Sobolev spaces
  $\displaystyle H_{s}^{p}\ :$  they make an interpolating scale
 with respect to  $\displaystyle s,p\ ;$  they have the same
 type and cotype than  $\displaystyle L^{p}$  spaces.\ \par 
\ \par 
\quad  In section~\ref{iSH12} we start the study of the multipliers
 algebra  $\displaystyle {\mathcal{M}}_{s}^{p}$  of  $\displaystyle
 H_{s}^{p}.$  We prove that  $\displaystyle {\mathcal{M}}_{s}^{p}$
  is invariant by the automorphisms of the ball.\ \par 
\ \par 
\quad  In the following section we study Carleson measures and Carleson
 sequences.\ \par 
\ \par 
\quad  In the following section we study links between  $p$  interpolating
 sequences of vectors in a general Banach space  $B$  and Carleson
 measures and basic sequence in  $\displaystyle \ell ^{p}.$ 
 We study also algebras of operators on  $B$  which diagonalize
 along a sequence of vectors in  $B.$  Application to  $\displaystyle
 H_{s}^{p}$  are done.\ \par 
\ \par 
\quad  Then, in the next harmonic analysis section, we develop a very
 useful method due to S. Drury~\cite{Drury70}, for the union
 of two Sidon sets, to fit Hardy Sobolev spaces.\ \par 
\ \par 
\quad  The following section contains the results on interpolating
 sequences of points for the multipliers algebra  $\displaystyle
 {\mathcal{M}}_{s}^{p}.$ \ \par 
\quad  In section~\ref{iSH13} we study the notion of dual boundedness
 in the framework of Hardy Sobolev spaces.\ \par 
\ \par 
\quad  Finally in the appendix we put technical lemmas to ease the
 reading of section~\ref{6IM3}.\ \par 
\ \par 
\quad  In the sequel we shall deal only with finite sequences of points
  $\displaystyle S\subset {\mathbb{B}}$  but with estimates not
 depending of the number of points in  $\displaystyle S.$  The
 results for infinite sequences is then got by a normal family argument.\ \par 

\section{Hardy Sobolev spaces.~\label{1HS0}}
\quad  By a result of J. Ortega and J. Fabrega~\cite{OrtFab97}, corollary
 3.4, (see also E. Ligocka~\cite{Ligocka87}), we have that the
 Hardy Sobolev spaces  $\displaystyle H_{s}^{p}$  form a interpolating
 scale with respect to  $s$  and  $p.$  This means that for 
 $\displaystyle 1<p_{0},p_{1}<\infty ,\ 0\leq s_{0},s_{1},\ 0<\theta
 <1$  and  $\displaystyle \ \frac{1}{p}=\frac{(1-\theta )}{p_{0}}+\frac{\theta
 }{p_{1}},\ s=(1-\theta )s_{0}+\theta s_{1},$  we have\ \par 
\quad \quad \quad  \begin{equation}  (H_{s_{0}}^{p_{0}},H_{s_{1}}^{p_{1}})_{\lbrack
 \theta \rbrack }=H_{s}^{p}.\label{iSH18}\end{equation}\ \par 
\quad  We shall use this result in relation with the Banach spaces
 interpolation method. In particular we shall prove results essentially
 when  $s$  is an integer and, by use of it, we shall get the
 same results for  $s$  real.\ \par 

\subsection{Similarity between  $\displaystyle H_{s}^{p}$  and
  $\displaystyle H^{p}.$ }

\begin{Dfnt}
 Let  $S$  be a sequence, we set  $\epsilon :=\lbrace \epsilon
 _{a},\ a\in S\rbrace \in {\mathcal{R}}(S),$  a Rademacher sequence,where
 the random variables  $\epsilon _{a}$  are Bernouilli independent
 and such that  $\displaystyle P(\epsilon _{a}=1)=P(\epsilon _{a}=-1)=1/2.$ 
\end{Dfnt}
\quad  Let  $\displaystyle f(\epsilon ,z)\in H_{s}^{p}$  for any value
 of the random variable  $\epsilon \in {\mathcal{R}}(S),$  then :\ \par 

\begin{Lmm}
 ~\label{S1}We have\par 
\quad \quad \quad  $\displaystyle \ {\left\Vert{{\mathbb{E}}(f)}\right\Vert}_{H_{s}^{p}}^{p}\lesssim
 {\mathbb{E}}({\left\Vert{f}\right\Vert}_{H_{s}^{p}}^{p}).$ 
\end{Lmm}
\quad  Proof.\ \par 
We can take as an equivalent norm in  $\displaystyle H_{s}^{p}$
  the sum of the  $\displaystyle H^{p}$  norms of the  $\displaystyle
 R^{k}$  derivatives, i.e. with  $\displaystyle p<\infty ,$ \ \par 
\quad \quad \quad  $\displaystyle \ {\left\Vert{f}\right\Vert}_{H_{s}^{p}}^{p}\simeq
 \sum_{k=0}^{s}{{\left\Vert{R^{k}(f)}\right\Vert}_{H^{p}}^{p}}.$ \ \par 
Hence, because  ${\mathbb{E}}$  is linear, we have\ \par 
\quad \quad \quad  $\displaystyle R^{k}{\mathbb{E}}(f)={\mathbb{E}}(R^{k}f).$ \ \par 
On the other hand\ \par 
\quad \quad \quad  $\ \left\vert{{\mathbb{E}}(g)}\right\vert ^{p}\leq ({\mathbb{E}}(\left\vert{g}\right\vert
 ))^{p}\leq {\mathbb{E}}(\left\vert{g}\right\vert ^{p})$ \ \par 
hence\ \par 
\quad \quad \quad  $\displaystyle \ {\left\Vert{{\mathbb{E}}(g)}\right\Vert}_{H^{p}}^{p}=\int_{\partial
 {\mathbb{B}}}{\left\vert{{\mathbb{E}}(g)}\right\vert ^{p}d\sigma
 }\leq \int_{\partial {\mathbb{B}}}{{\mathbb{E}}(\left\vert{g}\right\vert
 ^{p})d\sigma }={\mathbb{E}}(\int_{\partial {\mathbb{B}}}{\left\vert{g}\right\vert
 ^{p}d\sigma })={\mathbb{E}}({\left\Vert{g}\right\Vert}_{H^{p}}^{p}).$ \ \par 
So applying this with  $\displaystyle g=R^{k}f$  we get\ \par 
\quad \quad \quad  $\ {\left\Vert{R^{k}({\mathbb{E}}(f))}\right\Vert}_{H^{p}}^{p}={\left\Vert{{\mathbb{E}}(R^{k}(f))}\right\Vert}_{H^{p}}^{p}\leq
 {\mathbb{E}}({\left\Vert{R^{k}(f)}\right\Vert}_{H^{p}}^{p}),$ \ \par 
and\ \par 
\quad \quad \quad  $\displaystyle \ {\left\Vert{{\mathbb{E}}(f)}\right\Vert}_{H_{s}^{p}}^{p}\lesssim
 {\mathbb{E}}({\left\Vert{f}\right\Vert}_{H_{s}^{p}}^{p}).$ 
  $\displaystyle \hfill\blacksquare $ \ \par 
\ \par 

\begin{Prps}
 ~\label{1HS1}The spaces  $\displaystyle H_{s}^{p}$  have, for
 any  $s\in {\mathbb{N}},$  the same type as  $\displaystyle H^{p}.$ 
\end{Prps}
\quad  Proof.\ \par 
We can prove it by use of the fact that the Sobolev spaces 
 $\displaystyle W_{k}^{p}$  have this property by~\cite{Cobos86},
 but because the problem is on the boundary of the ball which
 is not isotropic with respect to the derivatives, we shall prove
 it directly.\ \par 
\quad  So let  $\displaystyle p\leq 2$  we want to prove that  $\displaystyle
 H_{s}^{p}$  is of type  $\displaystyle p,$  which means, with
  $\epsilon \in {\mathcal{R}}(1,...,N)$  a Rademacher sequence
 and  ${\mathbb{E}}$  the expectation,\ \par 
\quad \quad \quad  $\displaystyle ({\mathbb{E}}({\left\Vert{\sum_{j=1}^{N}{\epsilon
 _{j}f_{j}}}\right\Vert}_{H_{s}^{p}}^{2}))^{1/2}\leq T_{p}(\sum_{j=1}^{N}{{\left\Vert{f_{j}}\right\Vert}_{H_{s}^{p}}^{p}})^{1/p}.$
 \ \par 
We can take as a norm in  $\displaystyle H_{s}^{p}$  the sum
 of the  $\displaystyle H^{p}$  norms of the  $\displaystyle
 R^{k}$  derivatives, hence, because  ${\mathbb{E}}$  is linear,
 it suffices to have\ \par 
\quad \quad \quad  $\displaystyle ({\mathbb{E}}({\left\Vert{\sum_{j=1}^{N}{\epsilon
 _{j}R^{k}(f_{j})}}\right\Vert}_{H^{p}}^{2}))^{1/2}\leq T_{p}(\sum_{j=1}^{N}{{\left\Vert{f_{j}}\right\Vert}_{H_{s}^{p}}^{p}})^{1/p}.$
 \ \par 
But  $\displaystyle H^{p}$  being a subspace of  $\displaystyle
 L^{p}(\partial {\mathbb{B}}),$  it is already of type  $p$  hence\ \par 
\quad \quad \quad  $\displaystyle ({\mathbb{E}}({\left\Vert{\sum_{j=1}^{N}{\epsilon
 _{j}R^{k}(f_{j})}}\right\Vert}_{H^{p}}^{2})^{1/2}\leq T_{p}(\sum_{j=1}^{N}{{\left\Vert{R^{k}(f_{j})}\right\Vert}_{H^{p}}^{p}})^{1/p},$
 \ \par 
So, because  $\displaystyle f\in H_{s}^{p}$  implies\ \par 
\quad \quad \quad \quad  $\displaystyle \forall k\leq s,\ R^{k}(f)\in H^{p},\ {\left\Vert{R^{k}(f)}\right\Vert}_{H^{p}}\leq
 {\left\Vert{f}\right\Vert}_{H_{s}^{p}},$ \ \par 
we get\ \par 
\quad \quad \quad  $\displaystyle ({\mathbb{E}}({\left\Vert{\sum_{j=1}^{N}{\epsilon
 _{j}R^{k}(f_{j})}}\right\Vert}_{H^{p}}^{2}))^{1/2}\leq T_{p}(\sum_{j=1}^{N}{{\left\Vert{f_{j}}\right\Vert}_{H_{s}^{p}}^{p}})^{1/p},$
 \ \par 
and, adding a finite number of terms, we get\ \par 
\quad \quad \quad  $\displaystyle ({\mathbb{E}}({\left\Vert{\sum_{j=1}^{N}{\epsilon
 _{j}f_{j}}}\right\Vert}_{H_{s}^{p}}^{2}))^{1/2}\leq (s+1)T_{p}(\sum_{j=1}^{N}{{\left\Vert{f_{j}}\right\Vert}_{H_{s}^{p}}^{p}})^{1/p}.$
 \ \par 
If  $\displaystyle p>2,$  then the dual space of  $\displaystyle
 H_{s}^{p}$  is  $\displaystyle H_{s}^{p'}$  with  $\displaystyle
 p'<2$  hence  $\displaystyle H_{s}^{p'}$  is of type  $\displaystyle
 p'\ ;$  this implies that the dual of  $\displaystyle H_{s}^{p'},$
  namely  $\displaystyle H_{s}^{p}$  is of cotype  $\displaystyle
 p.$   $\displaystyle \hfill\blacksquare $ \ \par 
\quad  Using it we get the following theorem.\ \par 

\begin{Thrm}
 ~\label{S3}The spaces  $\displaystyle H_{s}^{p}$  have, for
 any  $s\in {\mathbb{R}}_{+},$  the same type as  $\displaystyle H^{p}.$ 
\end{Thrm}
\quad  Proof.\ \par 
Fix  $N\in {\mathbb{N}}$  and consider the space  $\displaystyle
 (H_{s}^{p})^{N}$  with the following  $\displaystyle \ell ^{p}$  norm :\ \par 
\quad \quad \quad  $\displaystyle \forall f=(f_{1},...,f_{N})\in (H_{s}^{p})^{N},\
 {\left\Vert{f}\right\Vert}_{p}:=(\sum_{j=1}^{N}{{\left\Vert{f_{j}}\right\Vert}_{H_{s}^{p}}^{p}})^{1/p}.$
 \ \par 
Consider the linear operator  $T\ :\ {\mathcal{R}}(1,...,N){\times}(H_{s}^{p})^{N}\rightarrow
 H_{s}^{p}$  defined by\ \par 
\quad \quad \quad  $\displaystyle \forall f=(f_{1},...,f_{N})\in (H_{s}^{p})^{N},\
 T_{N}(\epsilon ,f):=\sum_{j=1}^{N}{\epsilon _{j}f_{j}}\in H_{s}^{p}.$ \ \par 
To say that  $\displaystyle H_{s}^{p}$  is of type  $p$  means
 that, for any  $\displaystyle N\geq 1,$ \ \par 
\quad \quad \quad \quad  $\displaystyle ({\mathbb{E}}({\left\Vert{T_{N}(\epsilon ,f)}\right\Vert}_{H_{s}^{p}}^{2}))^{1/2}\leq
 C{\left\Vert{f}\right\Vert}_{p},$ \ \par 
i.e. the linear operator  $\displaystyle T_{N}$  is bounded
 from  $\displaystyle F_{s}:=(H_{s}^{p})^{N}$  equipped with
 the norm  $\displaystyle \ {\left\Vert{\cdot }\right\Vert}_{p}$
  to  $\displaystyle L^{2}(\Omega ,H_{s}^{p}),$  the space  $L^{2}(\Omega
 ,{\mathcal{A}},P)$  with value in  $\displaystyle H_{s}^{p}.$
  Because the  $\displaystyle H_{s}^{p}$  form an interpolating
 scale with respect to the parameter  $s\in {\mathbb{R}}_{+},$
  we have the same for the spaces  $\displaystyle F_{s}$  and
  $\displaystyle L^{2}(\Omega ,H_{s}^{p}).$ \ \par 
\quad  Fix  $\displaystyle p\leq 2$  and  $s\in {\mathbb{N}}\ ;$  by
 the proposition~\ref{1HS1} we have that  $T_{N}$  is bounded
 from  $\displaystyle F_{s}$  to  $\displaystyle L^{2}(\Omega
 ,H_{s}^{p}),$  and from  $\displaystyle F_{0}$  to  $\displaystyle
 L^{2}(\Omega ,H_{0}^{p}),$  the constant being independent of
  $N\in {\mathbb{N}},$  hence by interpolation  $\displaystyle
 T_{N}$  is bounded from  $\displaystyle F_{r}$  to  $\displaystyle
 L^{2}(\Omega ,H_{r}^{p}),$  for any  $\displaystyle 0\leq r\leq
 s,$  with a constant independent of  $\displaystyle N\in {\mathbb{N}}.$
  This proves that  $\displaystyle H_{r}^{p}$  is of type  $p$
  for any real  $\displaystyle r\in \lbrack 0,s\rbrack .$  By
 duality as in proposition~\ref{1HS1} we have that for  $\displaystyle
 p>2,\ H_{s}^{p}$  is of cotype  $\displaystyle p.$   $\displaystyle
 \hfill\blacksquare $ \ \par 
\ \par 
\quad  Up to a constant, we have the Young inequalities.\ \par 

\begin{Prps}
 ~\label{S2}We have, with  $\displaystyle \ \frac{1}{r}=\frac{1}{p}+\frac{1}{q},$
 \par 
\quad \quad \quad  $\displaystyle \forall f\in H_{s}^{p},\ \forall g\in H_{s}^{q},\
 fg\in H_{s}^{r}$  and  $\displaystyle \ {\left\Vert{fg}\right\Vert}_{H_{s}^{r}}\leq
 C_{s}{\left\Vert{f}\right\Vert}_{H_{s}^{p}}{\left\Vert{g}\right\Vert}_{H_{s}^{q}}.$
 
\end{Prps}
\quad  Proof.\ \par 
We have to compute the  $\displaystyle H^{p}$  norm, for  $\displaystyle
 j=0,...,s,$  of, by Leibnitz formula,\ \par 
\quad \quad \quad  \begin{equation}  R^{j}(fg)=\sum_{k=0}^{j}{C_{j}^{k}R^{k}(f)R^{(j-k)}(g)}.\label{S0}\end{equation}\
 \par 
By Minkowski inequality it is enough to control the norm of\ \par 
\quad \quad \quad  $\displaystyle R^{k}(f)R^{(j-k)}(g).$ \ \par 
But by Young inequality\ \par 
\quad \quad \quad  $\displaystyle \ {\left\Vert{R^{k}(f)R^{(j-k)}(g)}\right\Vert}_{H^{r}}\leq
 {\left\Vert{R^{k}(f)}\right\Vert}_{H^{p}}{\left\Vert{R^{(j-k)}(g)}\right\Vert}_{H^{q}}.$
 \ \par 
Now  $\displaystyle f\in H_{s}^{p}$  implies\ \par 
\quad \quad \quad \quad  $\displaystyle \forall k\leq s,\ R^{k}(f)\in H^{p},\ {\left\Vert{R^{k}(f)}\right\Vert}_{H^{p}}\leq
 {\left\Vert{f}\right\Vert}_{H_{s}^{p}}.$ \ \par 
The same  $\displaystyle g\in H_{s}^{q}$  implies\ \par 
\quad \quad \quad \quad  $\displaystyle \forall k\leq s,\ R^{k}(g)\in H^{p},\ {\left\Vert{R^{k}(g)}\right\Vert}_{H^{q}}\leq
 {\left\Vert{g}\right\Vert}_{H_{s}^{q}}.$ \ \par 
So\ \par 
\quad \quad \quad  $\displaystyle \forall j=0,...,s,\ \forall k\leq j,\ {\left\Vert{R^{k}(f)R^{(j-k)}(g)}\right\Vert}_{H^{r}}\leq
 {\left\Vert{f}\right\Vert}_{H_{s}^{p}}{\left\Vert{g}\right\Vert}_{H_{s}^{q}}.$
 \ \par 
Because we have a finite number of terms in~(\ref{S0}) we get
 the existence of a constant  $\displaystyle C_{s}>0$  such that\ \par 
\quad \quad \quad  $\displaystyle \ {\left\Vert{fg}\right\Vert}_{H_{s}^{r}}\leq
 C_{s}{\left\Vert{f}\right\Vert}_{H_{s}^{p}}{\left\Vert{g}\right\Vert}_{H_{s}^{q}},$
 \ \par 
which proves the proposition.  $\displaystyle \hfill\blacksquare $ \ \par 

\section{The multipliers algebra of  $\displaystyle H_{s}^{p}.$ ~\label{iSH12}}
\quad  Recall that the multipliers algebra  $\displaystyle {\mathcal{M}}_{s}^{p}$
  of  $\displaystyle H_{s}^{p}$  is the algebra of functions
  $m$  on  $\displaystyle {\mathbb{B}}$  such that 	 $\displaystyle
 \forall h\in H_{s}^{p},\ mh\in H_{s}^{p},$  and its norm is
 its norm as an operator from  $\displaystyle H_{s}^{p}$  into
  $\displaystyle H_{s}^{p}.$ \ \par 
\ \par 
\quad  As an easy corollary of the interpolating result~(\ref{iSH18}),
 we get the following theorem.\ \par 

\begin{Thrm}
 ~\label{iSH19}We have the embedding :\par 
\quad \quad \quad  ${\mathcal{M}}_{s}^{p}\subset {\mathcal{M}}_{r}^{p},$  for 
 $\displaystyle 1<p<\infty $  and  $\displaystyle 0\leq r\leq
 s,$  with  $\displaystyle \forall m\in {\mathcal{M}}_{s}^{p},\
 {\left\Vert{m}\right\Vert}_{{\mathcal{M}}_{r}^{p}}\leq {\left\Vert{m}\right\Vert}_{{\mathcal{M}}_{s}^{p}}.$
 
\end{Thrm}
\quad  Proof.\ \par 
Let  $\displaystyle m\in {\mathcal{M}}_{s}^{p}$  then  $m$ 
 is also in  $\displaystyle {\mathcal{M}}_{0}^{p}=H^{\infty }({\mathbb{B}})$
  which means that  $m$  is a bounded operator on  $\displaystyle
 H_{s}^{p}$  and on  $\displaystyle H_{0}^{p}.$  Hence  $m$ 
 is bounded on  $\displaystyle H_{r}^{p}$  for any  $\displaystyle
 r\in \lbrack 0,s\rbrack ,$  by Banach spaces interpolation~(\ref{iSH18}),
 which means that  $\displaystyle m\in {\mathcal{M}}_{r}^{p}.$
  Moreover we have  $\displaystyle \ {\left\Vert{m}\right\Vert}_{H^{\infty
 }({\mathbb{B}})}\leq {\left\Vert{m}\right\Vert}_{{\mathcal{M}}_{s}^{p}}$
  hence  $\displaystyle \ {\left\Vert{m}\right\Vert}_{{\mathcal{M}}_{r}^{p}}\leq
 {\left\Vert{m}\right\Vert}_{{\mathcal{M}}_{s}^{p}}.$   $\displaystyle
 \hfill\blacksquare $ \ \par 
\ \par 

\subsection{Invariance by automorphisms.}
\quad \quad  	Let  $\displaystyle e_{a}(z):=\frac{(1-\left\vert{a}\right\vert
 ^{2})^{\rho /2}}{(1-\bar a\cdot z)^{\rho }},\ \rho :=n-2s,$
  the normalized reproducing kernel for the point  $\displaystyle
 a\in {\mathbb{B}}$  in  $\displaystyle H_{s}^{2}.$ \ \par 
\quad \quad  	We shall show the following theorem which is true for any 
 $s\in {\mathbb{R}}_{+}.$ \ \par 

\begin{Thrm}
Let  $\varphi $  be an automorphism of the ball  $\displaystyle
 {\mathbb{B}}\ ;$  for any  $\displaystyle a\in {\mathbb{B}},$
  there is a number  $\displaystyle \eta (\varphi ,a)$  of modulus
 one such that, setting  $\displaystyle U(\varphi )e_{a}:=\eta
 (\varphi ,a)e_{\varphi (a)},$   $\displaystyle U(\varphi )$
  extends as an unitary representation of   $\displaystyle \mathrm{A}\mathrm{u}\mathrm{t}({\mathbb{B}})$
  in  $\displaystyle {\mathcal{L}}(H_{s}^{2}).$ 
\end{Thrm}
\quad \quad  	Proof.\ \par 
We shall adapt the proof of theorem 2 p. 35 in~\cite{AmarThesis77}.
 We know that  $\displaystyle \mathrm{A}\mathrm{u}\mathrm{t}({\mathbb{B}})$
  is isomorphic to  $\displaystyle U(n,1)$  the group of isometries
 for the sesquilinear form of  $\displaystyle {\mathbb{C}}^{n+1}$  :\ \par 
\quad \quad \quad \quad \quad 	 $\displaystyle (z,w):=\sum_{j=1}^{n}{z_{j}\bar w_{j}}-z_{0}\bar
 w_{0}.$ \ \par 
\quad \quad  	Let  $T\in U(n,1)\ ;$  in the canonical basis of  $\displaystyle
 {\mathbb{C}}^{n+1}$  its matrix  $\displaystyle \lbrack T\rbrack
 $  can be written by blocs :\ \par 
\quad \quad \quad \quad \quad 	 $\displaystyle \lbrack T\rbrack ={\left[{
\begin{matrix}
{A}&{B}\cr 
{C}&{D}\cr 
\end{matrix}
}\right]}\ ;$ \ \par 
where  $\displaystyle A$  is a  $\displaystyle n{\times}n$ 
 matrix,  $B$  is  $\displaystyle n{\times}1,$   $C$  is  $\displaystyle
 1{\times}n$  and  $D$  is  $\displaystyle 1{\times}1.$  The
 automorphism associated to  $T$  is then\ \par 
\quad \quad \quad \quad \quad 	 $\displaystyle \forall z\in {\mathbb{B}},\ \varphi (z):=\frac{AZ+B}{CZ+D},$
  where  $Z=(z_{1},...,\ z_{n}).$ \ \par 
If  $\alpha ,\beta $  are two vectors in  $\displaystyle {\mathbb{C}}^{n},$
  we denote by  $\alpha \cdot \bar \beta $  their scalar product
 ; the scalar product in  $\displaystyle H_{s}^{2}$  is still
 denoted by  $\displaystyle \ {\left\langle{\cdot ,\cdot }\right\rangle}.$
 \ \par 
We have\ \par 
\quad \quad \quad \quad \quad 	 $\displaystyle \ {\left\langle{e_{\varphi (a)},\ e_{\varphi
 (b)}}\right\rangle}=\frac{(1-\left\vert{\varphi (a)}\right\vert
 ^{2})^{\rho /2}(1-\left\vert{\varphi (b)}\right\vert ^{2})^{\rho
 /2}}{(1-\varphi (a)\cdot \bar \varphi (b))^{\rho }}.$ \ \par 
But\ \par 
\quad \quad \quad \quad \quad 	 $\displaystyle 1-\varphi (a)\cdot {\overline{\varphi (b)}}=1-\frac{Aa+B}{Ca+D}\cdot
 {\overline{{\left({\frac{Ab+B}{Cb+D},}\right)}}}=$ \ \par 
\quad \quad \quad \quad \quad \quad \quad \quad \quad \quad \quad \quad \quad 	 $\displaystyle =\frac{1}{(Ca+D){\overline{(Cb+D)}}}\lbrack
 (Ca+D){\overline{(Cb+D)}}-(Aa+B)({\overline{Ab+B}})\rbrack .$ \ \par 
Let  $\displaystyle (X,t)$  and  $\displaystyle (Y,v)$  two
 elements in  $\displaystyle {\mathbb{C}}^{n+1}$  and set  $\displaystyle
 \alpha =T(X,t),\ \beta =T(Y,v)$  we get\ \par 
\quad \quad \quad \quad \quad 	 $\displaystyle (\alpha ,\beta )=(AX+Bt)({\overline{AY+Bv}})-(CX+Dt)({\overline{CY+Dv}})=X\cdot
 \bar Y-t\bar v,$ \ \par 
because  $T$  let  $\displaystyle (\cdot ,\cdot )$  invariant.\ \par 
\quad \quad  	Back to the inhomogeneous coordinates  $\displaystyle a=X/t,\
 b=Y/v$  we get\ \par 
\quad \quad \quad \quad \quad 	 $\displaystyle (Ca+D){\overline{(Cb+D)}}-(Aa+B)({\overline{Ab+B}})=1-a\cdot
 \bar b,$ \ \par 
hence, putting it in  $\displaystyle \ {\left\langle{e_{\varphi
 (a)},\ e_{\varphi (b)}}\right\rangle}$  we get\ \par 
\quad \quad \quad \quad \quad 	 $\displaystyle \ {\left\langle{e_{\varphi (a)},\ e_{\varphi
 (b)}}\right\rangle}=\frac{(Ca+D)^{\rho }}{\left\vert{Ca+D}\right\vert
 ^{\rho }}{\times}\frac{({\overline{Cb+D}})^{\rho }}{\left\vert{Cb+D}\right\vert
 ^{\rho }}{\left\langle{e_{a},e_{b}}\right\rangle}.$ \ \par 
The linear combinations of  $\displaystyle \lbrace e_{c},\ c\in
 {\mathbb{B}}\rbrace $  being dense in  $\displaystyle H_{s}^{2},$
  we define on them the operator  $\displaystyle U(\varphi )$  by\ \par 
\quad \quad \quad \quad \quad 	 $\displaystyle U(\varphi )e_{a}:=\eta (\varphi ,a)e_{\varphi (a)},$ \ \par 
where  $\displaystyle \eta (\varphi ,a):=\frac{({\overline{Ca+D}})^{\rho
 }}{\left\vert{Ca+D}\right\vert ^{\rho }}$  is of modulus  $\displaystyle
 1.$ \ \par 
\quad  The previous computation gives\ \par 
\quad \quad \quad \quad \quad 	 $\displaystyle \ {\left\langle{U(\varphi )e_{a},U(\varphi )e_{b}}\right\rangle}={\left\langle{e_{a},e_{b}}\right\rangle}$
 \ \par 
hence  $\displaystyle U(\varphi )$  is unitary. Moreover  $\displaystyle
 U(\varphi )$  is a representation of  $\displaystyle \mathrm{A}\mathrm{u}\mathrm{t}({\mathbb{B}}).$
  To see this we have to show that :\ \par 
\quad \quad \quad \quad \quad 	 $\displaystyle \forall \psi ,\varphi \in \mathrm{A}\mathrm{u}\mathrm{t}({\mathbb{B}}),\
 \forall a\in {\mathbb{B}},\ \eta (\psi \circ \varphi ,\ a)=\eta
 (\psi ,\varphi (a)){\times}\eta (\varphi ,a).$ \ \par 
Setting\ \par 
\quad \quad \quad \quad \quad 	 $\displaystyle \varphi (a):=\frac{Aa+B}{Ca+D},\ \psi (b):=\frac{Ab+B}{Cb+D},$
 \ \par 
the computation is easy.  $\displaystyle \hfill\blacksquare $ \ \par 

\begin{Rmrq}
We can use equivalently the following identities (Theorem 2.2.2
 p. 26 in~\cite{RudinBall81})\par 
\quad \quad \quad \quad \quad 	 $\displaystyle 1-\varphi (a)\cdot {\overline{\varphi (b)}}=\frac{(1-\left\vert{\mu
 }\right\vert ^{2})}{(1-\bar \mu \cdot a)}{\times}\frac{(1-\bar
 b\cdot a)}{(1-\mu \cdot \bar b)},$ \par 
\quad \quad \quad \quad \quad 	 $\displaystyle 1-\left\vert{\varphi (a)}\right\vert ^{2}=\frac{(1-\left\vert{\mu
 }\right\vert ^{2})(1-\left\vert{a}\right\vert ^{2})}{\left\vert{(1-\bar
 a\cdot z)}\right\vert ^{2}},$ \par 
where  $\displaystyle \varphi (z)=\varphi _{\mu }(z)$  is the
 automorphism exchanging  $\mu $  and  $\displaystyle 0.$  In
 any case we get\par 
\quad \quad \quad \quad \quad 	 $\displaystyle \eta (\varphi _{\mu },a)=\frac{(1-\bar \mu \cdot
 a)^{\rho }}{\left\vert{1-\bar \mu \cdot a}\right\vert ^{\rho }}.$ 
\end{Rmrq}
\ \par 

\begin{Crll}
 The space of multipliers  $\displaystyle {\mathcal{M}}_{s}^{2}$
  of  $\displaystyle H_{s}^{2}$  is invariant by  $\displaystyle
 \mathrm{A}\mathrm{u}\mathrm{t}({\mathbb{B}}).$ 
\end{Crll}
\quad \quad  	Proof.\ \par 
Let  $\displaystyle m\in {\mathcal{M}}_{s}^{2}$  then we have\ \par 
\quad \quad \quad  $\displaystyle \forall a\in {\mathbb{B}},\ m^{*}k_{a}=\bar m(a)k_{a}$ \ \par 
because\ \par 
\quad \quad \quad \quad  $\displaystyle \forall h\in H_{s}^{2},\ {\left\langle{h,m^{*}k_{a}}\right\rangle}={\left\langle{mh,k_{a}}\right\rangle}=m(a)h(a)=m(a){\left\langle{h,k_{a}}\right\rangle}.$
 \ \par 
Hence\ \par 
\quad \quad \quad \quad \quad 	 $\displaystyle m^{*}U(\varphi )e_{a}=m^{*}(\eta (\varphi ,a)e_{\varphi
 (a)})=\eta (\varphi ,a)m^{*}e_{\varphi (a)}=\eta (\varphi ,a){\overline{m(\varphi
 (a))}}e_{\varphi (a)}.$ \ \par 
So\ \par 
\quad \quad \quad \quad  	\begin{equation}  U^{-1}(\varphi )m^{*}U(\varphi )e_{a}=U^{-1}(\varphi
 )(\eta (\varphi ,a){\overline{m(\varphi (a))}}e_{\varphi (a)})=\eta
 (\varphi ,a){\overline{m(\varphi (a))}}U^{-1}(\varphi )e_{\varphi
 (a)}.\label{MI0}\end{equation}\ \par 
But from  $\displaystyle U(\varphi )e_{a}:=\eta (\varphi ,a)e_{\varphi
 (a)},$  we get\ \par 
\quad \quad \quad \quad \quad 	 $\displaystyle e_{a}=U^{-1}Ue_{a}=\eta U^{-1}e_{\varphi (a)}\Rightarrow
 U^{-1}e_{\varphi (a)}=\eta ^{-1}e_{a}$ \ \par 
and putting this in~(\ref{MI0}) we get\ \par 
\quad \quad \quad \quad \quad 	 $\displaystyle U^{-1}(\varphi )m^{*}U(\varphi )e_{a}=\eta {\overline{m(\varphi
 (a))}}\eta ^{-1}e_{a}={\overline{m(\varphi (a))}}e_{a}=(m\circ
 \varphi )^{*}e_{a}.$ \ \par 
So by the density of the linear combinations of the  $\displaystyle
 \lbrace e_{a},\ a\in {\mathbb{B}}\rbrace $  we get\ \par 
\quad \quad \quad  $\displaystyle (m\circ \varphi )^{*}=U^{-1}(\varphi )m^{*}U(\varphi ).$ \ \par 
Now because  $U(\varphi )$  is unitary on  $\displaystyle H_{s}^{2}$
  we have\ \par 
\quad \quad \quad  $\ {\left\Vert{(m\circ \varphi )^{*}}\right\Vert}_{{\mathcal{L}}(H_{s}^{2})}={\left\Vert{m^{*}}\right\Vert}_{{\mathcal{L}}(H_{s}^{2})}\Rightarrow
 {\left\Vert{m\circ \varphi }\right\Vert}_{{\mathcal{M}}_{s}^{2}}={\left\Vert{m}\right\Vert}_{{\mathcal{M}}_{s}^{2}}.$
 	  $\displaystyle \hfill\blacksquare $ \ \par 

\section{Carleson measures and Carleson sequences.}
\quad  Let  $\displaystyle Q(\zeta ,h):=\lbrace z\in \bar {\mathbb{B}},\
 \left\vert{1-\bar \zeta z}\right\vert <h\rbrace $  be the "pseudo
 ball" centered at  $\displaystyle \zeta \in \partial {\mathbb{B}}$
  and of radius  $\displaystyle h>0.$ \ \par 
\quad  We shall use the following well known lemma.\ \par 

\begin{Lmm}
 ~\label{3CS0}If  $\displaystyle \mu $  is a Carleson measure
 for  $\displaystyle H_{s}^{p},$  then  $\displaystyle \mu (Q(\zeta
 ,h))\lesssim h^{n-sp}=\left\vert{Q(\zeta ,h)\cap \partial {\mathbb{B}}}\right\vert
 ^{1-p\frac{s}{n}}.$ 
\end{Lmm}
\quad  Proof.\ \par 
Because  $\displaystyle \mu $  is a Carleson measure for  $\displaystyle
 H_{s}^{p},$  we have  $\displaystyle \ \int_{{\mathbb{B}}}{\left\vert{k_{a}(z)}\right\vert
 ^{p}d\mu }\lesssim {\left\Vert{k_{a}}\right\Vert}_{s,p}^{p}\ ;$ \ \par 
recall that  $\displaystyle k_{a}(z)=\frac{1}{(1-\bar az)^{\rho
 }}$  with  $\rho =n-2s,$  then we get, with\ \par 
\quad \quad \quad \quad \quad   $\displaystyle Q_{a}:=Q(\frac{a}{\left\vert{a}\right\vert },1-\left\vert{a}\right\vert
 )\iff \lbrace z\in {\mathbb{B}}::\left\vert{1-\bar a\cdot z}\right\vert
 <h\rbrace ,\ h:=(1-\left\vert{a}\right\vert ),$ \ \par 
that\ \par 
\quad \quad \quad  $\displaystyle \ \int_{Q_{a}}{\left\vert{\frac{1}{(1-\bar az)^{\rho
 }}}\right\vert ^{p}d\mu }\leq \int_{{\mathbb{B}}}{\left\vert{k_{a}(z)}\right\vert
 ^{p}d\mu }\lesssim {\left\Vert{k_{a}}\right\Vert}_{s,p}^{p}\ ;$ \ \par 
hence\ \par 
\quad \quad \quad \quad \quad 	 $\displaystyle \ \frac{1}{h^{\rho p}}\mu (Q_{a})\lesssim {\left\Vert{k_{a}}\right\Vert}_{s,p}^{p}\simeq
 (1-\left\vert{a}\right\vert ^{2})^{-\rho p-sp+n}\Rightarrow
 \mu (Q_{a})\lesssim h^{n-sp}.$   $\displaystyle \hfill\blacksquare $ \ \par 
\ \par 
\quad  Let us recall the definitions of Carleson sequences.\ \par 

\begin{Dfnt}
 The sequence  $S$  is Carleson in  $\displaystyle H_{s}^{p}({\mathbb{B}}),$
  if the associated measure\par 
\quad \quad \quad \quad \quad 	 $\displaystyle \nu _{S}:=\sum_{a\in S}{(1-\left\vert{a}\right\vert
 ^{2})^{n-sp}\delta _{a}}$ \par 
is Carleson for  $\displaystyle H_{s}^{p}({\mathbb{B}}).$ 
\end{Dfnt}
\quad  At this point we notice that the coefficients of the measure
  $\nu _{S}$  depend on the parameter  $\displaystyle s.$ \ \par 

\begin{Lmm}
 ~\label{SH4}Let  $S$  be  sequence in  $\displaystyle {\mathbb{B}}$
  which is Carleson for  $\displaystyle H_{s}^{p}$  and for 
 $\displaystyle H_{0}^{p}=H^{p}$  then  $S$  is Carleson for
  $\displaystyle H_{r}^{p},\ 0\leq r\leq s.$ 
\end{Lmm}
\quad \quad  	Proof.\ \par 
Consider the linear operator\ \par 
\quad \quad \quad \quad \quad 	 $\displaystyle T\ :\ H_{r}^{p}\rightarrow \ell ^{p}(w_{r}),\
 Tf:=\lbrace f(a)\rbrace _{a\in S}$ \ \par 
with the weight  $\displaystyle w_{r}(a):=(1-\left\vert{a}\right\vert
 ^{2})^{n-pr}.$  Because  $S$  is Carleson  $\displaystyle H_{s}^{p}$
 we have that  $T$  is bounded from  $\displaystyle H_{s}^{p}$
  to  $\displaystyle \ell ^{p}(w_{s}),$  i.e.\ \par 
\quad \quad \quad \quad \quad 	 $\displaystyle \ \sum_{a\in S}{(1-\left\vert{a}\right\vert
 ^{2})^{n-ps}\left\vert{f(a)}\right\vert ^{p}}\lesssim {\left\Vert{f}\right\Vert}_{H_{s}^{p}}.$
 \ \par 
The same for  $\displaystyle s=0,$  i.e.\ \par 
\quad \quad \quad \quad \quad 	 $\displaystyle \ \sum_{a\in S}{(1-\left\vert{a}\right\vert
 ^{2})^{n}\left\vert{f(a)}\right\vert ^{p}}\lesssim {\left\Vert{f}\right\Vert}_{H^{p}}$
 \ \par 
hence, because  $\displaystyle w_{\zeta }(a)=(1-\left\vert{a}\right\vert
 ^{2})^{n-p\zeta }$  is holomorphic in the strip  $\displaystyle
 0\leq \Re \zeta \leq s$  and the scale of  $\displaystyle \lbrace
 H_{s}^{p}\rbrace _{s>0}$  forms an interpolating scale by the
 interpolating result~(\ref{iSH18}), we have that  $\displaystyle
 T$  is bounded from  $\displaystyle H_{r}^{p}$  to  $\displaystyle
 \ell ^{p}(w_{r})$  which means exactly that  $S$  is Carleson
 for  $\displaystyle H_{r}^{p},\ 0\leq r\leq s.$   $\displaystyle
 \hfill\blacksquare $ \ \par 
\ \par 
\quad  If  $\mu $  is a Carleson measure for  $\displaystyle H_{s}^{p},$
  then it is a Carleson measure for  $\displaystyle H_{r}^{p},\
 \forall r\geq s,$  simply because  $\displaystyle \ {\left\Vert{f}\right\Vert}_{H_{s}^{p}}\leq
 {\left\Vert{f}\right\Vert}_{H_{r}^{p}}.$  For the Carleson sequences,
 this goes the opposite way.\ \par 

\begin{Thrm}
 ~\label{3CS1}If the sequence  $S$  is Carleson in  $\displaystyle
 H_{s}^{p}({\mathbb{B}}),$  then  $S$  is Carleson in  $\displaystyle
 H_{r}^{p}$  for all  $r,\ 0\leq r\leq s.$ 
\end{Thrm}
\quad  Proof.\ \par 
We first show that the measure  $\displaystyle \mu :=\sum_{a\in
 S}{(1-\left\vert{a}\right\vert ^{2})^{n}\delta _{a}}$  is Carleson
  $\displaystyle V^{1},$  i.e. that\ \par 
\quad \quad \quad  $\displaystyle \ \sum_{a\in S\cap Q(\zeta ,h)}{(1-\left\vert{a}\right\vert
 ^{2})^{n}}\lesssim h^{n}.$ \ \par 
For this we have that  $\nu _{S}$  Carleson in  $\displaystyle
 H_{s}^{p}({\mathbb{B}})$  implies that  $\nu _{S}$  is finite,
 just using lemma~\ref{3CS0} with  $\displaystyle Q_{a}\supset
 {\mathbb{B}}.$  So we have  $\displaystyle \ \sum_{a\in S}{(1-\left\vert{a}\right\vert
 ^{2})^{n-sp}}\leq C.$  Now still with lemma~\ref{3CS0} we get\ \par 
\quad \quad \quad  $\displaystyle \ \sum_{a\in S\cap Q(\zeta ,h)}{(1-\left\vert{a}\right\vert
 ^{2})^{n-sp}}=\mu (Q(\zeta ,h))\lesssim h^{n-sp}.$ \ \par 
But  $\displaystyle a\in Q(\zeta ,h)\Rightarrow (1-\left\vert{a}\right\vert
 ^{2})<h$  hence, with\ \par 
\quad \quad \quad \quad  $\displaystyle (1-\left\vert{a}\right\vert ^{2})^{n}=(1-\left\vert{a}\right\vert
 ^{2})^{sp}(1-\left\vert{a}\right\vert ^{2})^{n-sp}\leq h^{sp}(1-\left\vert{a}\right\vert
 ^{2})^{n-sp}$ \ \par 
we get\ \par 
\quad \quad \quad  $\displaystyle \ \sum_{a\in S\cap Q(\zeta ,h)}{(1-\left\vert{a}\right\vert
 ^{2})^{n}}\leq h^{sp}\sum_{a\in S\cap Q(\zeta ,h)}{(1-\left\vert{a}\right\vert
 ^{2})^{n-sp}}=h^{sp}\mu (Q(\zeta ,h))\lesssim h^{n}.$ \ \par 
This is valid for all  $\displaystyle Q(\zeta ,h)$  so we get
 that the measure  $\displaystyle \mu :=\sum_{a\in S}{(1-\left\vert{a}\right\vert
 ^{2})^{n}\delta _{a}}$  is Carleson  $\displaystyle V^{1},$
  or, equivalently Carleson  $\displaystyle H^{p}:=H_{0}^{p}.$ \ \par 
\quad  Now we apply lemma~\ref{SH4} to end the proof of the theorem.
  $\displaystyle \hfill\blacksquare $ \ \par 

\section{General results}
\quad  We shall establish a link between Carleson sequences and sequences
 like canonical basis of  $\displaystyle \ell ^{p}.$ \ \par 
Let  $B$  be a Banach space,  $\displaystyle B'$  its dual.\ \par 
\ \par 

\begin{Dfnt}
 We say that the sequence of bounded vectors  $\displaystyle
 \lbrace e_{a}\rbrace _{a\in S}$  in  $B$  is equivalent to a
 canonical basis of  $\displaystyle \ell ^{p}$  if\par 
\quad \quad \quad  $\displaystyle \exists B_{p}>0,\ \forall \lambda \in \ell ^{p}(S),\
 \frac{1}{B_{p}}{\left\Vert{\lambda }\right\Vert}_{\ell ^{p}}\leq
 {\left\Vert{\sum_{a\in S}{\lambda _{a}e_{a}}}\right\Vert}_{B}\leq
 B_{p}C{\left\Vert{\lambda }\right\Vert}_{\ell ^{p}}.$ \par 

\end{Dfnt}

\begin{Dfnt}
 We say that the sequence of bounded vectors  $\displaystyle
 \lbrace e_{a}\rbrace _{a\in S}$  in  $B$  is  $p$  interpolating if\par 
\quad \quad \quad  $\displaystyle \exists I_{p}>0,\ \forall \mu \in \ell ^{p'}(S),\
 \exists h\in B',\ {\left\Vert{h}\right\Vert}_{B'}\leq I_{p}{\left\Vert{\mu
 }\right\Vert}_{\ell ^{p'}}::\forall a\in S,\ {\left\langle{h,e_{a}}\right\rangle}=\mu
 _{a}.$ \par 

\end{Dfnt}

\begin{Dfnt}
 We say that the sequence of bounded vectors  $\displaystyle
 \lbrace e_{a}\rbrace _{a\in S}$  in  $B$  is dual bounded if\par 
\quad \quad \quad  $\displaystyle \exists C>0,\ \exists \lbrace \rho _{a}\rbrace
 _{a\in S}\subset B'::\forall a\in S,\ {\left\Vert{f_{a}}\right\Vert}_{B'}\leq
 C,\ {\left\langle{\rho _{a},e_{b}}\right\rangle}=\delta _{ab}.$ 
\end{Dfnt}
\quad  Clearly if  $\displaystyle \lbrace e_{a}\rbrace _{a\in S}$ 
 is  $p$  interpolating then it is dual bounded : just interpolate
 the basic sequence of  $\displaystyle \ell ^{p}(S).$ \ \par 

\begin{Dfnt}
 We say that the sequence of bounded vectors  $\displaystyle
 \lbrace e_{a}\rbrace _{a\in S}$  in  $B$  is  $p$  Carleson if\par 
\quad \quad \quad  $\displaystyle \exists C_{p}>0,\ \forall h\in B',\ \sum_{a\in
 S}{\left\vert{{\left\langle{h,e_{a}}\right\rangle}}\right\vert
 ^{p'}}\leq C_{p}^{p'}{\left\Vert{h}\right\Vert}_{B'}^{p'}.$ 
\end{Dfnt}
\quad  We have :\ \par 

\begin{Lmm}
 ~\label{cB1}Let  $\displaystyle \lbrace e_{a}\rbrace _{a\in
 S}$  be a sequence in  $B$  of bounded vectors, then the following
 assertions are equivalent :\par 
\quad  (i)  $\lbrace e_{a}\rbrace _{a\in S}$  is  $p$  Carleson in
  $\displaystyle B.$ \par 
\quad  (ii)  $\displaystyle \lbrace e_{a}\rbrace _{a\in S}$  verifies
  $\displaystyle \exists C>0,\ \forall \lambda \in \ell ^{p}(S),\
 {\left\Vert{\sum_{a\in S}{\lambda _{a}e_{a}}}\right\Vert}_{B}\leq
 C{\left\Vert{\lambda }\right\Vert}_{\ell ^{p}(S)}.$ 
\end{Lmm}
\quad  Proof.\ \par 
Suppose that  $\lbrace e_{a}\rbrace _{a\in S}$  verifies the
 {\sl (i)} of the lemma, then using the duality  $\displaystyle
 B-B'$  we have\ \par 
\quad \quad \quad  $\displaystyle \forall \lambda \in \ell ^{p}(S),\ \forall h\in
 B',\ \left\vert{\sum_{a\in S}{\lambda _{a}{\left\langle{e_{a},h}\right\rangle}}}\right\vert
 \leq C{\left\Vert{\lambda }\right\Vert}_{\ell ^{p}}{\left\Vert{h}\right\Vert}_{B'}.$
 \ \par 
By the duality  $\displaystyle \ell ^{p}-\ell ^{p'}$  we get then\ \par 
\quad \quad \quad  $\displaystyle \ \sum_{a\in S}{\left\vert{{\left\langle{e_{a},h}\right\rangle}}\right\vert
 ^{p'}}\leq C^{p'}{\left\Vert{h}\right\Vert}_{B'}^{p'},$ \ \par 
which is the definition of  $\lbrace e_{a}\rbrace _{a\in S}\
 p$  Carleson in  $\displaystyle B.$ \ \par 
Suppose now that  $\lbrace e_{a}\rbrace _{a\in S}$  verifies
 the {\sl (ii)} of the lemma, this means\ \par 
\quad \quad \quad  $\displaystyle \ \sum_{a\in S}{\left\vert{{\left\langle{e_{a},h}\right\rangle}}\right\vert
 ^{p'}}\leq C^{p'}{\left\Vert{h}\right\Vert}_{H_{s}^{p'}}^{p'}$ \ \par 
which leads by the duality  $\displaystyle \ell ^{p}-\ell ^{p'}$   to\ \par 
\quad \quad \quad  $\displaystyle \forall \lambda \in \ell ^{p}(S),\ \ \left\vert{\sum_{a\in
 S}{\lambda _{a}{\left\langle{e_{a},h}\right\rangle}}}\right\vert
 \leq C{\left\Vert{\lambda }\right\Vert}_{\ell ^{p}}{\left\Vert{h}\right\Vert}_{H_{s}^{p'}}$
 \ \par 
and with the duality  $\displaystyle B-B'$  to the {\sl (i)}
 of the lemma.  $\displaystyle \hfill\blacksquare $ \ \par 

\begin{Thrm}
 Let  $\lbrace e_{a}\rbrace _{a\in S}$  be a  $p$  interpolating
 sequence for  $\displaystyle B$  and suppose moreover that 
 $\lbrace e_{a}\rbrace _{a\in S}$  is  $p$  Carleson for  $\displaystyle
 B$  then  $\displaystyle \lbrace e_{a}\rbrace _{a\in S}$  makes
 a system equivalent to a canonical basis in  $\displaystyle \ell ^{p}.$ 
\end{Thrm}
\quad  Proof.\ \par 
We have to show that\ \par 
\quad \quad \quad \quad \quad  $\displaystyle \ \forall \lambda \in \ell ^{p},\ {\left\Vert{\sum_{a\in
 S}{\lambda _{a}e_{a}}}\right\Vert}_{B}\simeq {\left\Vert{\lambda
 }\right\Vert}_{\ell ^{p}}.$ \ \par 
We have\ \par 
\quad \quad \quad  $\displaystyle \ {\left\Vert{\sum_{a\in S}{\lambda _{a}e_{a}}}\right\Vert}_{B}=\sup
 \ _{h\in B',\ {\left\Vert{h}\right\Vert}\leq 1}\left\vert{\sum_{a\in
 S}{\lambda _{a}{\left\langle{e_{a},h}\right\rangle}}}\right\vert $ \ \par 
but by H\"older\ \par 
\quad \quad \quad  $\displaystyle \ \left\vert{\sum_{a\in S}{\lambda _{a}{\left\langle{e_{a},h}\right\rangle}}}\right\vert
 \leq {\left\Vert{\lambda }\right\Vert}_{\ell ^{p}}(\sum_{a\in
 S}{\left\vert{{\left\langle{e_{a},h}\right\rangle}}\right\vert
 ^{p'}})^{1/p'}$ \ \par 
and because  $S$  is  $p$  Carleson we have\ \par 
\quad \quad \quad  $\displaystyle (\sum_{a\in S}{\left\vert{{\left\langle{e_{a},h}\right\rangle}}\right\vert
 ^{p'}})^{1/p'}\leq C_{p}{\left\Vert{h}\right\Vert}_{B'},$ \ \par 
hence\ \par 
\quad \quad \quad  $\displaystyle \ {\left\Vert{\sum_{a\in S}{\lambda _{a}e_{a}}}\right\Vert}_{B}\leq
 C_{p}{\left\Vert{\lambda }\right\Vert}_{\ell ^{p}}.$ \ \par 
For the other direction we still have\ \par 
\quad \quad \quad  $\displaystyle \ {\left\Vert{\sum_{a\in S}{\lambda _{a}e_{a}}}\right\Vert}_{B}=\sup
 \ _{h\in B',\ {\left\Vert{h}\right\Vert}\leq 1}\left\vert{\sum_{a\in
 S}{\lambda _{a}{\left\langle{e_{a},h}\right\rangle}}}\right\vert $ \ \par 
but, because  $\lbrace e_{a}\rbrace _{a\in S}$  is  $\displaystyle
 p$  interpolating, we can find a  $\displaystyle h\in B'$  such that\ \par 
\quad \quad \quad  $\displaystyle \forall a\in S,\ {\left\langle{h,e_{a}}\right\rangle}=\mu
 _{a},\ {\left\Vert{h}\right\Vert}_{B'}\leq I_{p}{\left\Vert{\mu
 }\right\Vert}_{\ell ^{p'}}.$ \ \par 
So we get\ \par 
\quad \quad \quad  $\displaystyle \ {\left\Vert{\sum_{a\in S}{\lambda _{a}e_{a}}}\right\Vert}_{B}\geq
 \frac{1}{I_{p}}\left\vert{\sum_{a\in S}{\lambda _{a}{\left\langle{e_{a},h}\right\rangle}}}\right\vert
 =\frac{1}{I_{p}}\left\vert{\sum_{a\in S}{\lambda _{a}\bar \mu
 _{a}}}\right\vert ,$ \ \par 
and we choose  $\mu $  such that  $\displaystyle \ {\left\Vert{\mu
 }\right\Vert}_{\ell ^{p'}}=1$  and  $\displaystyle \ \sum_{a\in
 S}{\lambda _{a}\bar \mu _{a}}={\left\Vert{\lambda }\right\Vert}_{\ell
 ^{p}}.$ \ \par 
So we get\ \par 
\quad \quad \quad  $\displaystyle \ {\left\Vert{\sum_{a\in S}{\lambda _{a}e_{a}}}\right\Vert}_{B}\geq
 \frac{1}{I_{p}}{\left\Vert{\lambda }\right\Vert}_{\ell ^{p}}.$
   $\displaystyle \hfill\blacksquare $ \ \par 

\begin{Thrm}
 Let  $\lbrace e_{a}\rbrace _{a\in S}$  makes a system equivalent
 to the canonical basis in  $\displaystyle \ell ^{p}$  and suppose
 moreover that :\par 
\quad \quad \quad  $\displaystyle P_{S}\ :\ \varphi \in B\rightarrow P_{S}\varphi
 :=\sum_{a\in S}{{\left\langle{\varphi ,\rho _{a}}\right\rangle}e_{a}}$ \par 
is bounded, then  $\displaystyle \lbrace e_{a}\rbrace _{a\in
 S}$  is  $p$  Carleson and  $p$  interpolating with a bounded
 linear extension operator.
\end{Thrm}
\quad  Proof.\ \par 
Because  $\lbrace e_{a}\rbrace _{a\in S}$  makes a system equivalent
 to a canonical basis in  $\displaystyle \ell ^{p}$  means\ \par 
\quad \quad \quad  \begin{equation}  \forall \lambda \in \ell ^{p},\ {\left\Vert{\sum_{a\in
 S}{\lambda _{a}e_{a}}}\right\Vert}_{B}\simeq {\left\Vert{\lambda
 }\right\Vert}_{\ell ^{p}}\label{cB3}\end{equation}\ \par 
we have in particular that\ \par 
\quad \quad \quad  \begin{equation}  \ {\left\Vert{\sum_{a\in S}{\lambda _{a}e_{a}}}\right\Vert}_{B}\leq
 C_{p}{\left\Vert{\lambda }\right\Vert}_{\ell ^{p}}\label{cB2}\end{equation}\
 \par 
which, by lemma~\ref{cB1} gives that  $\displaystyle \lbrace
 e_{a}\rbrace _{a\in S}$  is  $p$  Carleson in  $\displaystyle B.$ \ \par 
\quad  Suppose first that  $S$  is finite, then there is a dual system
  $\displaystyle \lbrace \rho _{a}\rbrace _{a\in S}$  in  $\displaystyle
 B'.$  Set\ \par 
\quad \quad \quad  $\displaystyle \forall \mu \in \ell ^{p'}(S),\ h:=\sum_{a\in
 S}{\mu _{a}\rho _{a}}\ ;$ \ \par 
we have  $\displaystyle \ {\left\langle{h,e_{b}}\right\rangle}=\sum_{a\in
 S}{\mu _{a}{\left\langle{\rho _{a},e_{b}}\right\rangle}}=\mu _{b}$ \ \par 
hence  $h$  interpolates  $\mu .$  It remains to control its
 norm. We have\ \par 
\quad \quad \quad  $\displaystyle P_{S}\varphi =\sum_{a\in S}{{\left\langle{\varphi
 ,\rho _{a}}\right\rangle}e_{a}},\ {\left\Vert{P_{S}\varphi }\right\Vert}_{B}\leq
 C{\left\Vert{\varphi }\right\Vert}_{B},$ \ \par 
and by use of~(\ref{cB3}) we get\ \par 
\quad \quad \quad  $\displaystyle \ {\left\Vert{P_{S}\varphi }\right\Vert}_{B}\geq
 \frac{1}{B_{p}}(\sum_{a\in S}{\left\vert{{\left\langle{\varphi
 ,\rho _{a}}\right\rangle}}\right\vert ^{p}})^{1/p}$ \ \par 
hence\ \par 
\quad \quad \quad  \begin{equation}  (\sum_{a\in S}{\left\vert{{\left\langle{\varphi
 ,\rho _{a}}\right\rangle}}\right\vert ^{p}})^{1/p}\leq B_{p}C{\left\Vert{\varphi
 }\right\Vert}_{B}\label{cB4}\end{equation}\ \par 
which means that  $\displaystyle \lbrace \rho _{a}\rbrace _{a\in
 S}$  is  $\displaystyle p'$  Carleson.\ \par 
\quad  Now let us estimate the norm of  $h$ \ \par 
\quad \quad \quad  $\displaystyle \ {\left\Vert{h}\right\Vert}_{B'}={\left\Vert{\sum_{a\in
 S}{\mu _{a}\rho _{a}}}\right\Vert}_{B'}=\sup \ _{\varphi \in
 B,\ {\left\Vert{\varphi }\right\Vert}\leq 1}\left\vert{\sum_{a\in
 S}{\mu _{a}{\left\langle{\rho _{a},\varphi }\right\rangle}}}\right\vert
 $ \ \par 
but\ \par 
\quad \quad \quad  $\displaystyle \ \left\vert{\sum_{a\in S}{\mu _{a}{\left\langle{\rho
 _{a},\varphi }\right\rangle}}}\right\vert \leq {\left\Vert{\mu
 }\right\Vert}_{\ell ^{p'}}(\sum_{a\in S}{\left\vert{{\left\langle{\rho
 _{a},\varphi }\right\rangle}}\right\vert ^{p}})^{1/p}$ \ \par 
and by~(\ref{cB4}) we get\ \par 
\quad \quad \quad  $\displaystyle \ \left\vert{\sum_{a\in S}{\mu _{a}{\left\langle{\rho
 _{a},\varphi }\right\rangle}}}\right\vert \leq {\left\Vert{\mu
 }\right\Vert}_{\ell ^{p'}}B_{p}C{\left\Vert{\varphi }\right\Vert}_{B}$ \ \par 
so we have\ \par 
\quad \quad \quad  $\displaystyle \ {\left\Vert{h}\right\Vert}_{B'}\leq B_{p}C{\left\Vert{\mu
 }\right\Vert}_{\ell ^{p'}}.$ \ \par 
The bounded linear extension operator is then\ \par 
\quad \quad \quad  $\displaystyle \mu \in \ell ^{p'}\rightarrow E(\mu ):=\sum_{a\in
 S}{\mu _{a}\rho _{a}},\ {\left\Vert{E(\cdot )}\right\Vert}_{\ell
 ^{p'}\rightarrow B'}\leq B_{p}C.$ \ \par 
Hence we prove the theorem.  $\displaystyle \hfill\blacksquare $ \ \par 

\begin{Rmrq}
 The fact that  $\displaystyle P_{S}$  is bounded implies that
  $\displaystyle E_{S}:=\mathrm{S}\mathrm{p}\mathrm{a}\mathrm{n}(e_{a},\
 a\in S)$  is complemented in  $\displaystyle B.$  Just set :\par 
\quad \quad \quad  $\displaystyle \forall \varphi \in B,\ \varphi _{1}:=P_{S}\varphi
 \in E_{S},\ \varphi _{2}:=\varphi -\varphi _{1}.$ 
\end{Rmrq}

\begin{Lmm}
 ~\label{cB6}Let  $\displaystyle \lbrace e_{a}\rbrace _{a\in
 S}$  be dual bounded and such that  $\displaystyle \lbrace \rho
 _{a}\rbrace _{a\in S}$  is  $\displaystyle p'$  Carleson, then
  $\displaystyle \lbrace e_{a}\rbrace _{a\in S}$  is  $p$  interpolating
 with a bounded linear extension operator.
\end{Lmm}
\quad  Proof.\ \par 
Because  $\displaystyle \lbrace \rho _{a}\rbrace _{a\in S}$
  is  $\displaystyle p'$  Carleson we have\ \par 
\quad \quad \quad  \begin{equation}  \forall \varphi \in B\subset B'',\ (\sum_{a\in
 S}{\left\vert{{\left\langle{\rho _{a},\varphi }\right\rangle}}\right\vert
 ^{p}})^{1/p}\leq C_{p'}{\left\Vert{\varphi }\right\Vert}_{B}.\label{cB5}\end{equation}\
 \par 
Now take  $\mu \in \ell ^{p'}$  and set  $\displaystyle h:=\sum_{a\in
 S}{\mu _{a}\rho _{a}}$  we have\ \par 
\quad \quad \quad  $\displaystyle \ {\left\Vert{h}\right\Vert}_{B'}={\left\Vert{\sum_{a\in
 S}{\mu _{a}\rho _{a}}}\right\Vert}_{B'}=\sup \ _{\varphi \in
 B,\ {\left\Vert{\varphi }\right\Vert}\leq 1}\left\vert{\sum_{a\in
 S}{\mu _{a}{\left\langle{\rho _{a},\varphi }\right\rangle}}}\right\vert
 $ \ \par 
but\ \par 
\quad \quad \quad  $\displaystyle \ \left\vert{\sum_{a\in S}{\mu _{a}{\left\langle{\rho
 _{a},\varphi }\right\rangle}}}\right\vert \leq {\left\Vert{\mu
 }\right\Vert}_{\ell ^{p'}}(\sum_{a\in S}{\left\vert{{\left\langle{\rho
 _{a},\varphi }\right\rangle}}\right\vert ^{p}})^{1/p}$ \ \par 
and by~(\ref{cB5}) we get\ \par 
\quad \quad \quad  $\displaystyle \ \left\vert{\sum_{a\in S}{\mu _{a}{\left\langle{\rho
 _{a},\varphi }\right\rangle}}}\right\vert \leq {\left\Vert{\mu
 }\right\Vert}_{\ell ^{p'}}B_{p}C{\left\Vert{\varphi }\right\Vert}_{B}$ \ \par 
so we have\ \par 
\quad \quad \quad  $\displaystyle \ {\left\Vert{h}\right\Vert}_{B'}\leq B_{p}C{\left\Vert{\mu
 }\right\Vert}_{\ell ^{p'}}.$ \ \par 
The bounded linear extension operator is then\ \par 
\quad \quad \quad  $\displaystyle \mu \in \ell ^{p'}\rightarrow E(\mu ):=\sum_{a\in
 S}{\mu _{a}\rho _{a}},\ {\left\Vert{E(\cdot )}\right\Vert}_{\ell
 ^{p'}\rightarrow B'}\leq B_{p}C.$ \ \par 
Hence we prove the lemma.  $\displaystyle \hfill\blacksquare $ \ \par 

\subsection{Diagonalizing operators algebras.}
\quad  Let  $B$  be a Banach space and  $\displaystyle \lbrace e_{a}\rbrace
 _{a\in S}$  be a sequence of bounded vectors in  $B$  ; we shall
 work with operators  $M$  such that  $\displaystyle M\ :\ B\rightarrow
 B$  is bounded and\ \par 
\quad \quad \quad  $\displaystyle \forall a\in S,\ Me_{a}=m_{a}e_{a}.$ \ \par 
Let  ${\mathcal{A}}$  be a commutative algebra of operators
 on  $B$  diagonalizing on  $\displaystyle E:=\mathrm{S}\mathrm{p}\mathrm{a}\mathrm{n}\lbrace
 e_{a},\ a\in S\rbrace ,$  with the norm inherited from  $\displaystyle
 {\mathcal{L}}(B)\ ;$  we shall extend our definition of interpolation
 to this context.\ \par 

\begin{Dfnt}
 We say that the sequence of bounded vectors  $\displaystyle
 \lbrace e_{a}\rbrace _{a\in S}$  in  $B$  is interpolating for
  ${\mathcal{A}}$  if\par 
\quad \quad \quad  $\exists A>0,\ \forall \lambda \in \ell ^{\infty }(S),\ \exists
 M\in {\mathcal{A}},\ {\left\Vert{M}\right\Vert}_{{\mathcal{L}}(E)}\leq
 A{\left\Vert{\lambda }\right\Vert}_{\ell ^{\infty }}::\forall
 a\in S,\ Me_{a}=\lambda _{a}e_{a}.$ 
\end{Dfnt}
\quad  The first general result is in the special case of Hilbert spaces.\ \par 

\begin{Thrm}
 ~\label{cB9}Let  $H$  be a Hilbert space,  $\displaystyle \lbrace
 e_{a}\rbrace _{a\in S}$  be a sequence of normalized vectors in  $H.$ \par 
\quad  If   $\displaystyle \lbrace e_{a}\rbrace _{a\in S}$  is interpolating
 for  ${\mathcal{A}}$  then  $\displaystyle \lbrace e_{a}\rbrace
 _{a\in S}$  is equivalent to a basic sequence in  $\displaystyle
 \ell ^{2}(S).$ \par 
\quad  If  $\displaystyle \lbrace e_{a}\rbrace _{a\in S}$  is equivalent
 to a basic sequence in  $\displaystyle \ell ^{2}(S),$  set 
 $\displaystyle E:=\mathrm{S}\mathrm{p}\mathrm{a}\mathrm{n}\lbrace
 e_{a},\ a\in S\rbrace $  and  ${\mathcal{D}}$  the algebra of
 operators in  ${\mathcal{L}}(E)$   diagonalizing in  $\displaystyle
 \lbrace e_{a}\rbrace _{a\in S},$  then   $\displaystyle \lbrace
 e_{a}\rbrace _{a\in S}$  is interpolating for  ${\mathcal{D}}.$ 
\end{Thrm}
\quad  This theorem was proved in~\cite{AmarThesis77}, (Proposition
 3, p. 17) en route to a characterisation of interpolating sequences
 in the spectrum of a commutative algebra of operators in  ${\mathcal{L}}(H).$
  I shall reprove it here for the reader's convenience.\ \par 
\quad  Proof.\ \par 
Suppose that  $\displaystyle \lbrace e_{a}\rbrace _{a\in S}$
  is interpolating for  ${\mathcal{A}},$  and take  $\epsilon
 \in {\mathcal{R}}(S)$  a Rademacher sequence. Then  $\displaystyle
 \epsilon \in \ell ^{\infty }(S)$  hence there is an operator
  $M_{\epsilon }\in {\mathcal{A}}$  such that\ \par 
\quad \quad \quad  $M_{\epsilon }e_{a}=\epsilon _{a}e_{a},\ {\left\Vert{M_{\epsilon
 }}\right\Vert}_{{\mathcal{L}}(H)}\leq A.$ \ \par 
Now consider  $\displaystyle h:=\sum_{a\in S}{h_{a}e_{a}}\in
 E\subset H$  we have\ \par 
\quad \quad \quad  $\displaystyle M_{\epsilon }h=\sum_{a\in S}{\epsilon _{a}h_{a}e_{a}},$
  and  $\displaystyle \ {\left\Vert{M_{\epsilon }h}\right\Vert}_{H}\leq
 A{\left\Vert{h}\right\Vert}_{H}$ \ \par 
so\ \par 
\quad \quad \quad  $\displaystyle A^{2}{\left\Vert{h}\right\Vert}_{H}^{2}\geq {\mathbb{E}}({\left\Vert{M_{\epsilon
 }h}\right\Vert}_{H}^{2})=\sum_{a\in S}{\left\vert{h_{a}}\right\vert
 ^{2}{\left\Vert{e_{a}}\right\Vert}_{H}^{2}},$ \ \par 
because the  $\epsilon _{a}$  are independent and of mean  $\displaystyle
 0.$  So we get, the  $\displaystyle e_{a}$  being normalized,\ \par 
\quad \quad \quad  $\displaystyle \ \sum_{a\in S}{\left\vert{h_{a}}\right\vert
 ^{2}}\leq A^{2}{\left\Vert{h}\right\Vert}^{2}.$ \ \par 
Because  $\displaystyle \epsilon _{a}^{2}=1,$  we get  $\displaystyle
 M_{\epsilon }M_{\epsilon }=I_{d}$  on  $E,$  hence by the boundedness
 of  $\displaystyle M_{\epsilon },$ \ \par 
\quad \quad \quad  $\displaystyle \forall h\in E,\ h=M_{\epsilon }(M_{\epsilon
 }h)\Rightarrow {\left\Vert{h}\right\Vert}_{H}\leq A{\left\Vert{M_{\epsilon
 }h}\right\Vert}_{H}$ \ \par 
hence taking again expectation\ \par 
\quad \quad \quad  $\displaystyle \ {\left\Vert{h}\right\Vert}_{H}^{2}\leq A^{2}{\mathbb{E}}({\left\Vert{M_{\epsilon
 }h}\right\Vert}^{2})=A^{2}\sum_{a\in S}{\left\vert{h_{a}}\right\vert
 ^{2}}.$ \ \par 
So we proved\ \par 
\quad \quad \quad  $\displaystyle \ \frac{1}{A^{2}}\sum_{a\in S}{\left\vert{h_{a}}\right\vert
 ^{2}}\leq {\left\Vert{h}\right\Vert}_{H}^{2}\leq A^{2}\sum_{a\in
 S}{\left\vert{h_{a}}\right\vert ^{2}},$ \ \par 
which means that  $\displaystyle \lbrace e_{a}\rbrace _{a\in
 S}$  is equivalent to a basic sequence in  $\displaystyle \ell
 ^{2}(S).$ \ \par 
\ \par 
\quad  Now suppose that  $\displaystyle \lbrace e_{a}\rbrace _{a\in
 S}$  is equivalent to a basic sequence in  $\displaystyle \ell
 ^{2}(S).$  This means (see for instance~\cite{AmarThesis77})
 that there is a bounded operator  $Q$  in  ${\mathcal{L}}(E),$
  with  $\displaystyle Q^{-1}$  also bounded, and an orthonormal
 system  $\lbrace \eta _{a}\rbrace _{a\in S}$  in  $E$  such that\ \par 
\quad \quad \quad  $\displaystyle \forall a\in S,\ Q\eta _{a}=e_{a}.$ \ \par 
Let  $\displaystyle \lambda \in \ell ^{\infty }(S)$  then the
 diagonal operator  $\displaystyle T_{\lambda }\eta _{a}:=\lambda
 _{a}\eta _{a}$  is bounded on  $E$  with  $\displaystyle \ {\left\Vert{T_{\lambda
 }}\right\Vert}\leq {\left\Vert{\lambda }\right\Vert}_{\infty
 }.$  Now set\ \par 
\quad \quad \quad  $\displaystyle R_{\lambda }:=QT_{\lambda }Q^{-1},$ \ \par 
then we get\ \par 
\quad \quad \quad  $\displaystyle \forall a\in S,\ R_{\lambda }e_{a}=QT_{\lambda
 }\eta _{a}=Q\lambda _{a}\eta _{a}=\lambda _{a}e_{a}$ \ \par 
hence  $R\in {\mathcal{D}}$  and\ \par 
\quad \quad \quad \quad  $\ {\left\Vert{R_{\lambda }}\right\Vert}_{{\mathcal{L}}(E)}\leq
 {\left\Vert{Q}\right\Vert}_{{\mathcal{L}}(E)}{\left\Vert{Q^{-1}}\right\Vert}_{{\mathcal{L}}(E)}{\left\Vert{T_{\lambda
 }}\right\Vert}_{{\mathcal{L}}(E)}\leq C{\left\Vert{\lambda }\right\Vert}_{\infty
 }$ \ \par 
hence  $\displaystyle \lbrace e_{a}\rbrace _{a\in S}$  is interpolating
 for  ${\mathcal{D}}.$   $\displaystyle \hfill\blacksquare $ \ \par 
Now we shall need a definition.\ \par 

\begin{Dfnt}
 We shall say that the algebra  ${\mathcal{A}}$  separates the
 points  $\displaystyle \lbrace e_{a}\rbrace _{a\in S}$  if\par 
\quad \quad \quad  $\exists C>0,\ \forall a,b\neq a\in S,\ \exists M_{ab}\in {\mathcal{A}}::M_{ab}e_{a}=e_{a},\
 M_{ab}e_{b}=0$  and  $\ {\left\Vert{M_{ab}}\right\Vert}_{{\mathcal{A}}}\leq
 C.$ 
\end{Dfnt}
\quad  Then we have the following remark.\ \par 

\begin{Rmrq}
 Suppose that  ${\mathcal{A}}$  separates  $\displaystyle \lbrace
 e_{a}\rbrace _{a\in S},$  this implies easily that, for any
 finite set  $S$  and any  $\lambda \in \ell ^{\infty }(S)$ 
 there is a  $M\in {\mathcal{A}}$  such that  $\displaystyle
 \forall a\in S,\ M(a)e_{a}=\lambda _{a}e_{a}.$  Hence if there
 is  $\displaystyle C>0$  such that there is a  $M'\in {\mathcal{A}}$
  with  $\displaystyle M'_{\mid E}=M_{\mid E}$  and  $\ {\left\Vert{M'}\right\Vert}_{{\mathcal{A}}}\leq
 C{\left\Vert{M}\right\Vert}_{{\mathcal{L}}(E)}$  then, as a
 corollary of theorem~\ref{cB9}, we get that if  $\displaystyle
 \lbrace e_{a}\rbrace _{a\in S}$  is equivalent to a basic sequence
 in  $\displaystyle \ell ^{2}(S),$  then   $\displaystyle \lbrace
 e_{a}\rbrace _{a\in S}$  is interpolating for  ${\mathcal{A}}.$
  We say that  ${\mathcal{A}}$  is a Pick algebra if this property
 is true for  ${\mathcal{A}}.$  This is very well studied in
 the nice book by Agler and McCarthy~\cite{AglMCar02}.
\end{Rmrq}
\quad  We shall generalise this result to  $p$  interpolating sequences.\ \par 
Let  $B$  be a Banach space and  $\displaystyle \lbrace e_{a}\rbrace
 _{a\in S}$  be a sequence of normalized vectors in  $B.$ \ \par 
Recall that the Banach  $B$  is of type  $p$  if\ \par 
\quad \quad \quad  $\exists T_{p}>0::\forall N\in {\mathbb{N}},\ \epsilon \in {\mathcal{R}}(\lbrace
 1,...,N\rbrace ),\ \forall f_{1},...,f_{N}\in B,\ {\mathbb{E}}({\left\Vert{\sum_{j=1}^{N}{\epsilon
 _{j}f_{j}}}\right\Vert}^{2})^{1/2}\leq T_{p}(\sum_{j=1}^{N}{{\left\Vert{f_{j}}\right\Vert}_{B}^{p}})^{1/p}.$
 \ \par 

\begin{Thrm}
 ~\label{cB8}If  $\lbrace e_{a}\rbrace _{a\in S}$  is interpolating
 for  ${\mathcal{A}}$  :\par 
\quad  if  $B$  is of type  $p>1$  then  $\displaystyle \lbrace e_{a}\rbrace
 _{a\in S}$  is  $p$  Carleson.\par 
\quad  if  $\displaystyle B'$  is of type  $\displaystyle p'>1,\ $
 then there is a dual sequence  $\lbrace \rho _{a}\rbrace _{a\in
 S}\subset B'$  to  $\displaystyle \lbrace e_{a}\rbrace _{a\in
 S}$  and  $\displaystyle \lbrace \rho _{a}\rbrace _{a\in S}$
  is  $p'$  Carleson, hence  $\displaystyle \lbrace e_{a}\rbrace
 _{a\in S}$  is  $p$  interpolating for  $B$  with a bounded
 linear extension operator ;
\end{Thrm}
\quad  Proof.\ \par 
Because  $\lbrace e_{a}\rbrace _{a\in S}$  is interpolating
 for  ${\mathcal{A}}$  we have\ \par 
\quad \quad \quad  $\forall a\in S,\ \exists M_{a}\in {\mathcal{A}}::M_{a}e_{b}=\delta
 _{ab}e_{b},\ {\left\Vert{M_{a}}\right\Vert}_{{\mathcal{A}}}\leq A.$ \ \par 
Now fix  $\displaystyle a\in S$  and take  $\displaystyle h\in
 B'$  such that  $\displaystyle \ {\left\langle{h,e_{a}}\right\rangle}=1.$
  This  $h$  exists by Hahn Banach with norm  $1$  and\ \par 
\quad \quad \quad  $\displaystyle \ {\left\langle{M_{a}^{*}h,\ e_{b}}\right\rangle}={\left\langle{h,\
 M_{a}e_{b}}\right\rangle}=\delta _{ab}{\left\langle{h,e_{b}}\right\rangle}=\delta
 _{ab}.$ \ \par 
So, setting  $\displaystyle \rho _{a}:=M_{a}^{*}h,$  we get
  $\displaystyle \rho _{a}\in B',\ {\left\langle{\rho _{a},e_{b}}\right\rangle}=\delta
 _{ab}$  and  $\displaystyle \ {\left\Vert{\rho _{a}}\right\Vert}_{B'}\leq
 A{\left\Vert{h}\right\Vert}_{B'}\leq A.$  Doing the same for
 any  $\displaystyle a\in S$  we get that  $\displaystyle \lbrace
 \rho _{a}\rbrace _{a\in S}\subset B'$  exists hence  $\displaystyle
 \lbrace e_{a}\rbrace _{a\in S}$   is dual bounded.\ \par 
\ \par 
Now as above, take  $\epsilon \in {\mathcal{R}}(S).$  Then 
 $\displaystyle \epsilon \in \ell ^{\infty }(S)$  hence there
 is an operator  $M_{\epsilon }\in {\mathcal{A}}$  such that\ \par 
\quad \quad \quad  $M_{\epsilon }e_{a}=\epsilon _{a}e_{a},\ {\left\Vert{M_{\epsilon
 }}\right\Vert}_{{\mathcal{L}}(B)}\leq A.$ \ \par 
By duality,  $\displaystyle M_{\epsilon }^{*}:B'\rightarrow
 B'$  is such that  $\displaystyle M_{\epsilon }^{*}\rho _{a}=\epsilon
 _{a}\rho _{a},$  and  $\ {\left\Vert{M_{\epsilon }^{*}}\right\Vert}_{{\mathcal{L}}(B')}\leq
 A\ ;$  so let\ \par 
\quad \quad \quad  $\displaystyle \forall \mu \in \ell ^{p',\ }h:=\sum_{a\in S}{\mu
 _{a}\rho _{a}}.$ \ \par 
We have\ \par 
\quad \quad \quad  $\displaystyle \ {\left\Vert{M_{\epsilon }^{*}h}\right\Vert}_{B'}={\left\Vert{\sum_{a\in
 S}{\epsilon _{a}\mu _{a}\rho _{a}}}\right\Vert}_{B'}.$ \ \par 
Using  $\displaystyle M_{\epsilon }M_{\epsilon }=I_{d},$  we get\ \par 
\quad \quad \quad  $\displaystyle \ {\left\Vert{h}\right\Vert}_{B'}\leq A{\left\Vert{M_{\epsilon
 }^{*}h}\right\Vert}_{B'}=A{\left\Vert{\sum_{a\in S}{\epsilon
 _{a}\mu _{a}\rho _{a}}}\right\Vert}_{B'},$ \ \par 
hence, taking expectation,\ \par 
\quad \quad \quad  $\ {\left\Vert{h}\right\Vert}_{B'}\leq A{\mathbb{E}}({\left\Vert{\sum_{a\in
 S}{\epsilon _{a}\mu _{a}\rho _{a}}}\right\Vert}_{B'})$ \ \par 
so  $\displaystyle B'$  of type  $p'$  means  $\displaystyle
 {\mathbb{E}}({\left\Vert{\sum_{a\in S}{\epsilon _{a}\mu _{a}\rho
 _{a}}}\right\Vert}_{B'}^{2})^{1/2}\leq T_{p'}(\sum_{a\in S}{\left\vert{\mu
 _{a}}\right\vert ^{p'}{\left\Vert{\rho _{a}}\right\Vert}_{B'}^{p'}})^{1/p'},$
  hence\ \par 
\quad \quad \quad  $\displaystyle {\mathbb{E}}({\left\Vert{\sum_{a\in S}{\epsilon
 _{a}\mu _{a}\rho _{a}}}\right\Vert}_{B'})\leq {\mathbb{E}}({\left\Vert{\sum_{a\in
 S}{\epsilon _{a}\mu _{a}\rho _{a}}}\right\Vert}_{B'}^{2})^{1/2}\leq
 T_{p'}(\sum_{a\in S}{\left\vert{\mu _{a}}\right\vert ^{p'}{\left\Vert{\rho
 _{a}}\right\Vert}_{B'}^{p'}})^{1/p'},$ \ \par 
so\ \par 
\quad \quad \quad  $\displaystyle \ {\left\Vert{h}\right\Vert}_{B'}\leq AT_{p'}(\sum_{a\in
 S}{\left\vert{\mu _{a}}\right\vert ^{p'}{\left\Vert{\rho _{a}}\right\Vert}_{B'}^{p'}})^{1/p'}=AT_{p'}{\left\Vert{\mu
 }\right\Vert}_{\ell ^{p'}},$ \ \par 
which prove that  $\lbrace \rho _{a}\rbrace _{a\in S}$  is 
 $p'$  Carleson hence applying lemma~\ref{cB6} we get that  $\displaystyle
 \lbrace e_{a}\rbrace _{a\in S}$  is  $p$  interpolating for
  $B$  with a bounded linear extension operator.\ \par 
\ \par 
\quad  To get the second part, set  $\displaystyle \varphi :=\sum_{a\in
 S}{\lambda _{a}e_{a}}$  and use again\ \par 
\quad \quad \quad  $\displaystyle \varphi =M_{\epsilon }(M_{\epsilon }\varphi )\Rightarrow
 {\left\Vert{\varphi }\right\Vert}_{B}\leq A{\left\Vert{M_{\epsilon
 }\varphi }\right\Vert}_{B}=A{\left\Vert{\sum_{a\in S}{\epsilon
 _{a}\lambda _{a}e_{a}}}\right\Vert}_{B'},$ \ \par 
hence, taking expectation,\ \par 
\quad \quad \quad  $\ {\left\Vert{\varphi }\right\Vert}_{B}\leq A{\mathbb{E}}({\left\Vert{\sum_{a\in
 S}{\epsilon _{a}\lambda _{a}e_{a}}}\right\Vert}_{B})$ \ \par 
so if  $\displaystyle B$  is of type  $p$  then again\ \par 
\quad \quad \quad  $\displaystyle \ {\left\Vert{\varphi }\right\Vert}_{B}\leq AT_{p}(\sum_{a\in
 S}{\left\vert{\lambda _{a}}\right\vert ^{p}{\left\Vert{e_{a}}\right\Vert}_{B}^{p}})^{1/p}=AT_{p}{\left\Vert{\lambda
 }\right\Vert}_{\ell ^{p}},$ \ \par 
which prove that  $\lbrace e_{a}\rbrace _{a\in S}$  is  $p$
  Carleson.  $\displaystyle \hfill\blacksquare $ \ \par 

\subsection{Application to Hardy Sobolev spaces.}
\quad  Let  $\displaystyle H_{s}^{p}$  be the Hardy Sobolev space and
  $\displaystyle {\mathcal{M}}_{s}^{p}$  its multipliers algebra
 ; let also  $\displaystyle S\subset {\mathbb{B}}$  be a finite
 sequence of points in  $\displaystyle {\mathbb{B}}.$ \ \par 
Set, for  $\displaystyle a\in {\mathbb{B}},\ e_{a}:=\frac{k_{a}}{{\left\Vert{k_{a}}\right\Vert}_{H_{s}^{p}}}$
  the normalized reproducing kernel in  $\displaystyle H_{s}^{p}$
  for functions in  $\displaystyle H_{s}^{p'}.$ \ \par 
Then we have that \ \par 
\quad \quad \quad  $\displaystyle \forall m\in {\mathcal{M}}_{s}^{p},\ \forall
 a\in {\mathbb{B}},\ m^{*}k_{a}={\overline{m(a)}}k_{a}\Rightarrow
 m^{*}e_{a}={\overline{m(a)}}e_{a},$ \ \par 
because\ \par 
\quad \quad \quad \quad  $\displaystyle \forall h\in H_{s}^{p},\ {\left\langle{h,m^{*}k_{a}}\right\rangle}={\left\langle{mh,k_{a}}\right\rangle}=m(a)h(a)=m(a){\left\langle{h,k_{a}}\right\rangle}.$
 \ \par 
\quad  So we have that the adjoint of elements in  $\displaystyle {\mathcal{M}}_{s}^{p}$
  make an algebra diagonalizing in  $\displaystyle \lbrace e_{a}\rbrace
 _{a\in S}$  so we can apply the previous results with the diagonalizing
 algebra  $\displaystyle {\mathcal{A}}:=\lbrace m^{*},\ m\in
 {\mathcal{M}}_{s}^{p}\rbrace $  operating on  $\displaystyle
 H_{s}^{p'}.$ \ \par 
\quad  The first thing to know is that  $\displaystyle H_{s}^{p}$ 
 has the same type and cotype than  $\displaystyle L^{p}.$  We
 prove it directly in theorem~\ref{S3}.\ \par 
\quad  So we have  $H_{s}^{p},\ \forall s\in {\mathbb{R}}_{+},$  is
 of type  $\displaystyle \min \ (2,p)$  and of cotype  $\displaystyle
 \max \ (2,p),$  hence we can apply theorem~\ref{cB8} to get
 directly, for all real values of  $s\in \lbrack 0,n/p\rbrack ,$ \ \par 

\begin{Thrm}
 ~\label{cB10}If  $\lbrace e_{a}\rbrace _{a\in S}$  is interpolating
 for  $\displaystyle {\mathcal{M}}_{s}^{p}$  then  $\displaystyle
 \lbrace e_{a}\rbrace _{a\in S}$  is dual bounded and \par 
\quad  because  $\displaystyle H_{s}^{p}$  is of type  $\min \ (2,p)$
  then  $\displaystyle \lbrace e_{a}\rbrace _{a\in S}$  is  $\min
 \ (2,p)$  Carleson ;\par 
\quad  because  $\displaystyle H_{s}^{p'}$  of type  $\displaystyle
 \min \ (2,p')\ \ $ then  $\lbrace \rho _{a}\rbrace _{a\in S}$
  is  $\min \ (2,p')$  Carleson, hence  $\displaystyle \lbrace
 e_{a}\rbrace _{a\in S}$  is  $p$  interpolating for  $\displaystyle
 H_{s}^{p}$  with a bounded linear extension operator provided
 that  $\displaystyle p\geq 2.$ 
\end{Thrm}
\quad  In fact we shall prove later on a better result by use of harmonic
 analysis for the last case : we shall get rid of the condition
  $\displaystyle p\geq 2.$  Nevertheless we have, in the special
 case  $\displaystyle p=2,$  as an application of theorem~\ref{cB9},
 for all real values of  $s\in \lbrack 0,n/p\rbrack ,$  :\ \par 

\begin{Thrm}
 Let  $\displaystyle \lbrace e_{a}\rbrace _{a\in S}$  be a sequence
 of normalized vectors in  $\displaystyle H_{s}^{2}\ ;$  if 
 $\displaystyle \lbrace e_{a}\rbrace _{a\in S}$  is interpolating
 for  $\displaystyle {\mathcal{M}}_{s}^{2}$  then  $\displaystyle
 \lbrace e_{a}\rbrace _{a\in S}$  is equivalent to a basic sequence
 in  $\displaystyle \ell ^{2}(S).$  If  $\displaystyle {\mathcal{M}}_{s}^{2}$
  is a Pick algebra, i.e. if   $\displaystyle s\geq \frac{n-1}{2},$
  then  $\displaystyle \lbrace e_{a}\rbrace _{a\in S}$  equivalent
 to a basic sequence in  $\displaystyle \ell ^{2}(S)$  implies
 that   $\displaystyle \lbrace e_{a}\rbrace _{a\in S}$  is interpolating
 for  $\displaystyle {\mathcal{M}}_{s}^{2}.$ 
\end{Thrm}

\section{Harmonic analysis.~\label{5HA1}}
\quad  Let  $S$  be an interpolating sequence for the multipliers algebra
  ${\mathcal{M}}_{s}^{p}$  of  $\displaystyle H_{s}^{p}({\mathbb{B}})$
  and recall that the interpolating constant for  $S$  is the
 smallest number  $C=C(S)$  such that\ \par 
\quad \quad \quad \quad \quad 	 $\displaystyle \forall \lambda \in l^{\infty }(S),\ \exists
 m\in {\mathcal{M}}_{s}^{p}::\forall a\in S,\ m(a)=\lambda _{a}$
  and  $\ {\left\Vert{m}\right\Vert}_{{\mathcal{M}}_{s}^{p}}\leq
 C{\left\Vert{\lambda }\right\Vert}_{l^{\infty }}.$ \ \par 
\quad \quad  	We have easily  $\displaystyle {\mathcal{M}}_{s}^{p}\subset
 H^{\infty }({\mathbb{B}})$  with  $\displaystyle \forall m\in
 {\mathcal{M}}_{s}^{p},\ {\left\Vert{m}\right\Vert}_{H^{\infty
 }({\mathbb{B}})}\leq {\left\Vert{m}\right\Vert}_{{\mathcal{M}}_{s}^{p}}.$
 \ \par 
We shall develop here a very useful feature introduced by S.
 Drury~\cite{Drury70}. Consider a finite sequence in  $\displaystyle
 {\mathbb{B}}$  with interpolating constant  $\displaystyle C(S).$ \ \par 
Set  $\displaystyle N=\# S\in {\mathbb{N}},\ S:=\lbrace a_{1},...,a_{N}\rbrace
 \subset {\mathbb{B}}$  and  $\theta :=\exp \frac{2i\pi }{N}.\
 $  $S$  interpolating in  $\displaystyle {\mathcal{M}}_{s}^{p}$
  implies that\ \par 
\quad \quad \quad \quad \quad 	 $\displaystyle \forall j=1,...,N,\ \exists \beta (j,z)\in {\mathcal{M}}_{s}^{p}::\forall
 k=1,...,N,\ \beta (j,a_{k})=\theta ^{jk}$ \ \par 
and  $\displaystyle \forall j=1,...,N,\ {\left\Vert{\beta (j,\cdot
 )}\right\Vert}_{{\mathcal{M}}_{s}^{p}}\leq C(S).$ \ \par 
\quad \quad  	Let\ \par 
\quad \quad \quad \quad \quad 	 $\gamma (l,z):=\frac{1}{N}\sum_{j=1}^{N}{\theta ^{-jl}\beta
 (j,z)}\in {\mathcal{M}}_{s}^{p},\ {\left\Vert{\gamma (l,\cdot
 )}\right\Vert}_{{\mathcal{M}}_{s}^{p}}\leq C(S).$ \ \par 
this is the Fourier transform, on the group of  $\displaystyle
 n^{th}$  roots of unity, of the function  $\beta (\cdot ,z),$  i.e.\ \par 
\quad \quad \quad \quad \quad 	 $\displaystyle \gamma (l,z)=\hat \beta (l,z),$ \ \par 
the parameter  $\displaystyle z\in {\mathbb{B}}$  being fixed.\ \par 
\quad \quad  	We have\ \par 
\quad \quad \quad \quad  	\begin{equation}  \gamma (l,a_{k})=\frac{1}{N}\sum_{j=1}^{N}{\theta
 ^{-jl}\beta (j,a_{k})}=\delta _{lk}.\label{SMH9}\end{equation}\ \par 
Hence the  $\gamma (l,\cdot )$  make a dual bounded sequence
 for  $\displaystyle S,$  with a norm in  $\displaystyle {\mathcal{M}}_{s}^{p}$
  bounded by  $\displaystyle C(S).$ \ \par 
We have by Plancherel on this group\ \par 
\quad \quad \quad \quad  	\begin{equation}  \ \sum_{l=1}^{N}{\left\vert{\gamma (l,z)}\right\vert
 ^{2}}=\frac{1}{N}\sum_{j=1}^{N}{\left\vert{\beta (j,z)}\right\vert
 ^{2}}.\label{SM2}\end{equation}\ \par 
Multiplying on both side by  $\displaystyle \ \left\vert{h}\right\vert
 ^{2}$  with  $\displaystyle h\in H_{s}^{p}({\mathbb{B}}),$  we get\ \par 
\quad \quad \quad \quad \quad 	 $\displaystyle \ \sum_{l=1}^{N}{\left\vert{\gamma (l,z)h(z)}\right\vert
 ^{2}}=\frac{1}{N}\sum_{j=1}^{N}{\left\vert{\beta (j,z)h(z)}\right\vert
 ^{2}}$ \ \par 
and applying  $R^{j}$  on both sides, recalling  $\displaystyle
 R^{j}$  operates only on the holomorphic part,\ \par 
\quad \quad \quad \quad \quad 	 $\displaystyle \ \sum_{l=1}^{N}{\bar \gamma (l,\cdot )\bar
 hR^{j}(\gamma (l,\cdot )h)}=\frac{1}{N}\sum_{l=1}^{N}{\bar \beta
 (l,\cdot )\bar hR^{j}(\beta (l,\cdot )h)}$ \ \par 
and again  $\bar R^{j}$  on both sides\ \par 
\quad \quad \quad \quad  	\begin{equation}  \ \sum_{l=1}^{N}{\left\vert{R^{j}(\gamma
 (l,\cdot )h)}\right\vert ^{2}}=\frac{1}{N}\sum_{l=1}^{N}{\left\vert{R^{j}(\beta
 (l,\cdot )h)}\right\vert ^{2}.\label{SMH10}}\end{equation}\ \par 

\begin{Lmm}
 ~\label{iSH17}Let  $\displaystyle Q_{l}(k,z):=\underbrace{\beta
 *\cdot \cdot \cdot *\beta (k,z)}_{l\ times}$  then  $\ {\left\Vert{Q_{l}(k,\cdot
 )}\right\Vert}_{{\mathcal{M}}_{s}^{p}}\leq C(S)^{l}$  and hence
  $\displaystyle \ {\left\Vert{Q_{l}(k,\cdot )}\right\Vert}_{H^{\infty
 }({\mathbb{B}})}\leq {\left\Vert{Q_{l}(k,\cdot )}\right\Vert}_{{\mathcal{M}}_{s}^{p}}\leq
 C(S)^{l}.$ 
\end{Lmm}
\quad  Proof.\ \par 
Let  $Q_{2}(k,z):=\beta *\beta (k,z)=\frac{1}{N}\sum_{j=1}^{N}{\beta
 (j,z)\beta (k-j,z)}$  then, because  $\displaystyle {\mathcal{M}}_{s}^{p}$
  is an Banach algebra, we have\ \par 
\quad \quad \quad \quad  $\displaystyle \ {\left\Vert{\beta (j,\cdot )\beta (j-k,\cdot
 )}\right\Vert}_{{\mathcal{M}}_{s}^{p}}\leq {\left\Vert{\beta
 (j,\cdot )}\right\Vert}_{{\mathcal{M}}_{s}^{p}}{\left\Vert{\beta
 (j-k,\cdot )}\right\Vert}_{{\mathcal{M}}_{s}^{p}}\leq C(S)^{2.}$ \ \par 
Hence by induction we get the lemma.  $\displaystyle \hfill\blacksquare
 $ \ \par 

\begin{Lmm}
 ~\label{J12}We have\par 
\quad \quad \quad \quad \quad 	 $\displaystyle \ \sum_{k=1}^{N}{\left\vert{R^{j}(\gamma (k,\cdot
 )^{l}h)}\right\vert ^{2}}=\frac{1}{N}\sum_{k=1}^{N}{\left\vert{R^{j}(\underbrace{\beta
 *\beta *\cdot \cdot \cdot *\beta (k,\cdot )}_{l\ times}h)}\right\vert ^{2}}.$ 
\end{Lmm}
\quad \quad  	Proof.\ \par 
We have\ \par 
\quad \quad \quad \quad \quad 	 $\gamma (k,\cdot )^{l}=\widehat{\underbrace{\beta *\beta *\cdot
 \cdot \cdot *\beta (k,\cdot )}_{l\ times}}$ \ \par 
hence by Plancherel\ \par 
\quad \quad \quad \quad \quad 	 $\displaystyle \ \sum_{k=1}^{N}{\left\vert{\gamma (k,z)^{l}}\right\vert
 ^{2}}=\frac{1}{N}\sum_{k=1}^{N}{\left\vert{\beta *\cdot \cdot
 \cdot *\beta }\right\vert ^{2}(k,z)}$ \ \par 
and by lemma~\ref{iSH17}, because the  $\displaystyle {\mathcal{M}}_{s}^{p}$
  norm is bigger than the  $\displaystyle H^{\infty }({\mathbb{B}})$
  one,\ \par 
\quad \quad \quad  \begin{equation}  \ \forall z\in {\mathbb{B}},\ \sum_{k=1}^{N}{\left\vert{\gamma
 (k,z)^{l}}\right\vert ^{2}}\leq C(S)^{2l}.\label{iSH16}\end{equation}\ \par 
Multiplying by  $\displaystyle \ \left\vert{h}\right\vert ^{2}$
  on both sides, we get\ \par 
\quad \quad \quad \quad \quad 	 $\displaystyle \ \sum_{k=1}^{N}{\left\vert{\gamma (k,\cdot
 )^{l}h}\right\vert ^{2}}=\frac{1}{N}\sum_{k=1}^{N}{\left\vert{\beta
 *\cdot \cdot \cdot *\beta (k,\cdot )h}\right\vert ^{2}}$ \ \par 
and taking  $\displaystyle R^{j}$  derivatives, which operate
 only on the holomorphic part\ \par 
\quad \quad \quad \quad \quad 	 $\displaystyle \ \sum_{k=1}^{N}{R^{j}(\gamma (k,\cdot )^{l}h){\overline{\gamma
 (k,\cdot )^{l}h}}}=\frac{1}{N}\sum_{k=1}^{N}{R^{j}(\beta *\cdot
 \cdot \cdot *\beta (k,\cdot )h){\overline{\beta *\cdot \cdot
 \cdot *\beta (k,\cdot )h}}}\ ;$ \ \par 
Now we take  $\displaystyle \bar R^{j}$  derivatives on both
 sides to get the lemma.  $\displaystyle \hfill\blacksquare $ \ \par 
\quad  Let  $\displaystyle S:=\lbrace a_{1},...,a_{N}\rbrace \subset
 {\mathbb{B}}$  be a finite sequence in  $\displaystyle {\mathbb{B}}$
  then we have built the functions  $\displaystyle \lbrace \gamma
 (l,z)\rbrace _{l=1,...,N}\subset {\mathcal{M}}_{s}^{p}$  such that\ \par 
\quad \quad \quad  $\displaystyle \forall k,l=1,...,N,\ \gamma (l,a_{k})=\delta
 _{lk}$  and  $\displaystyle \ {\left\Vert{\gamma (l,\cdot )}\right\Vert}_{{\mathcal{M}}_{s}^{p}}\leq
 C(S)$ \ \par 
where  $\displaystyle C(S)$  is the interpolating constant of
 the sequence  $\displaystyle S.$  Now on we shall also use the
 notation   $\displaystyle \forall a\in S,\ \gamma _{a}(z):=\gamma
 (l,z)$  if  $\displaystyle a=a_{l}$  and we call  $\displaystyle
 \lbrace \gamma _{a}\rbrace _{a\in S}$  the canonical dual sequence
 for  $S$  in  $\displaystyle {\mathcal{M}}_{s}^{p}.$ \ \par 
\quad  The following proposition will be very useful for the sequel.\ \par 
\ \par 

\begin{Prps}
 ~\label{HA0}Let  $\displaystyle \lbrace \gamma _{a}\rbrace _{a\in
 S}$  be the canonical dual sequence for  $S$  in  $\displaystyle
 {\mathcal{M}}_{s}^{p}$  then we have\par 
\quad \quad \quad  $\displaystyle \forall l\geq 1,\ \forall z\in {\mathbb{B}},\
 \sum_{a\in S}{\left\vert{\gamma _{a}(z)}\right\vert ^{2l}}\leq
 C(S)^{2l},$ \par 
where  $\displaystyle C(S)$  is the interpolating constant for
  $\displaystyle S.$ 
\end{Prps}
\quad  Proof.\ \par 
This is just inequality~(\ref{iSH16}) with the new notations
  $\displaystyle \hfill\blacksquare $ \ \par 

\section{Interpolating sequences in the multipliers algebra.~\label{6IM3}}
\quad  We shall generalise theorem~\ref{cB10} valid for  $\displaystyle
 p\geq 2$  to all values of  $\displaystyle p\geq 1,$  but here
  $s$  must be an integer.\ \par 

\begin{Thrm}
 ~\label{lCS1}Let  $S$  be an interpolating sequence for  $\displaystyle
 {\mathcal{M}}_{s}^{p}$  and  $\gamma _{a}$  its canonical dual
 sequence, then, with  $\displaystyle e_{a}$  the normalised
 reproducing kernel for the point  $\displaystyle a\in {\mathbb{B}}$
  in  $\displaystyle H_{s}^{p},$ \par 
\quad \quad \quad \quad \quad 	 $\displaystyle \forall \lambda \in l^{p}(S),\ f:=\sum_{a\in
 S}{\lambda _{a}\gamma _{a}^{l}e_{a}}\in H_{s}^{p}({\mathbb{B}}),\
 {\left\Vert{f}\right\Vert}_{s,p}\lesssim {\left\Vert{\lambda
 }\right\Vert}_{p}.$ \par 
This means that  $S$  is interpolating for  $\displaystyle H_{s}^{p}$
  with the bounded extension property.
\end{Thrm}
\quad \quad  	Proof.\ \par 
As usual  $S$  is finite hence the series is well defined and we have\ \par 
\quad \quad \quad  $\displaystyle \forall b\in S,\ f(b)=\lambda _{b}e_{b}(b)=\lambda
 _{b}{\left\Vert{k_{b}}\right\Vert}_{H_{s}^{p'}}$ \ \par 
because by lemma~\ref{dBS1} :\ \par 
 $\displaystyle e_{b}(z):=\frac{k_{b}(z)}{{\left\Vert{k_{b}}\right\Vert}_{H_{s}^{p}}}\Rightarrow
 e_{b}(b):=\frac{k_{b}(b)}{{\left\Vert{k_{b}}\right\Vert}_{H_{s}^{p}}}=\frac{(1-\left\vert{b}\right\vert
 ^{2})^{2s-n}}{(1-\left\vert{b}\right\vert ^{2})^{s-n/p'}}={\left\Vert{k_{b}}\right\Vert}_{H_{s}^{p'}}.$
 \ \par 
This means that  $f$  interpolates the right values. So it remains
 to show that  $\displaystyle f\in H_{s}^{p}({\mathbb{B}}),\
 {\left\Vert{f}\right\Vert}_{H_{s}^{p}}\leq C{\left\Vert{\lambda
 }\right\Vert}_{p}.$ \ \par 
\quad  So we have to show that\ \par 
\quad \quad \quad  $\displaystyle \forall j\leq s,\ {\left\Vert{R^{j}f}\right\Vert}_{H^{p}}\leq
 C{\left\Vert{\lambda }\right\Vert}_{\ell ^{p}}.$ \ \par 
Fix  $\displaystyle j\leq s$  then\ \par 
\quad \quad \quad  $\displaystyle R^{j}(f)=\sum_{a\in S}{\lambda _{a}R^{j}(\gamma
 _{a}^{l}e_{a})}.$ \ \par 
By the exclusion proposition~\ref{iP2} with  $\displaystyle
 l>s,$  hence  $\displaystyle m:=\min \ (j,l)=j,$  we get\ \par 
\quad \quad \quad  $\displaystyle R^{j}(\gamma _{a}^{l}e_{a})=\sum_{q=0}^{j}{A_{q}\gamma
 ^{l-q}R^{j}(\gamma ^{q}e_{a})},$ \ \par 
because we have at most  $s$  terms, it is enough to control sums like :\ \par 
\quad \quad \quad  $\displaystyle T_{1}:=\sum_{a\in S}{\lambda _{a}\gamma _{a}^{l-q}R^{j}(\gamma
 _{a}^{q}e_{a})}.$ \ \par 
By H\"older we get\ \par 
\quad \quad \quad  $\displaystyle \ \left\vert{T_{1}}\right\vert ^{p}\leq (\sum_{a\in
 S}{\left\vert{\lambda _{a}}\right\vert ^{p}\left\vert{R^{j}(\gamma
 _{a}^{q}e_{a})}\right\vert ^{p}})(\sum_{a\in S}{\left\vert{\gamma
 _{a}}\right\vert ^{(l-q)p'}})^{p/p'}.$ \ \par 
Now we have by proposition~\ref{HA0}, provided that  $\displaystyle
 (l-q)p'\geq 2,$ \ \par 
\quad \quad \quad  $\displaystyle \forall z\in {\mathbb{B}},\ \sum_{a\in S}{\left\vert{\gamma
 _{a}}\right\vert ^{(l-q)p'}}\leq C(S)^{(l-q)p'}$ \ \par 
hence\ \par 
\quad \quad \quad  $\displaystyle \forall z\in {\mathbb{B}},\ \left\vert{T_{1}}\right\vert
 ^{p}\leq C(S)^{(l-q)p'}(\sum_{a\in S}{\left\vert{\lambda _{a}}\right\vert
 ^{p}\left\vert{R^{j}(\gamma _{a}^{q}e_{a})}\right\vert ^{p}}).$ \ \par 
So integrating\ \par 
\quad \quad \quad  $\displaystyle \forall r<1,\ \int_{\partial {\mathbb{B}}}{\left\vert{T_{1}(r\zeta
 )}\right\vert ^{p}d\sigma (\zeta )}\leq C(S)^{(l-q)p'}\int_{\partial
 {\mathbb{B}}}{\sum_{a\in S}{\left\vert{\lambda _{a}}\right\vert
 ^{p}\left\vert{R^{j}(\gamma _{a}^{q}e_{a})(r\zeta )}\right\vert
 ^{p}}d\sigma (\zeta )},$ \ \par 
hence\ \par 
\quad \quad \quad  \begin{equation}  \forall r<1,\ \int_{\partial {\mathbb{B}}}{\left\vert{T_{1}(r\zeta
 )}\right\vert ^{p}d\sigma (\zeta )}\leq C(S)^{(l-q)p'}\sum_{a\in
 S}{\left\vert{\lambda _{a}}\right\vert ^{p}\int_{\partial {\mathbb{B}}}{\left\vert{R^{j}(\gamma
 _{a}^{q}e_{a})(r\zeta )}\right\vert ^{p}d\sigma (\zeta )}}.\label{lCS0}\end{equation}\
 \par 
But we have\ \par 
\quad \quad \quad  $\displaystyle \gamma _{a}\in {\mathcal{M}}_{s}^{p}\Rightarrow
 \gamma _{a}^{q}\in {\mathcal{M}}_{s}^{p}$  with  $\displaystyle
 \ {\left\Vert{\gamma _{a}^{q}}\right\Vert}_{{\mathcal{M}}_{s}^{p}}\leq
 {\left\Vert{\gamma _{a}}\right\Vert}_{{\mathcal{M}}_{s}^{p}}^{q},$ \ \par 
because  $\displaystyle {\mathcal{M}}_{s}^{p}$  is a Banach algebra, so\ \par 
\quad \quad \quad  $\displaystyle \forall j\leq s,\ {\left\Vert{R^{j}(\gamma _{a}^{q}e_{a})}\right\Vert}_{H^{p}}\leq
 {\left\Vert{\gamma _{a}^{q}}\right\Vert}_{{\mathcal{M}}_{s}^{p}}{\left\Vert{e_{a}}\right\Vert}_{H_{s}^{p}}\leq
 {\left\Vert{\gamma _{a}}\right\Vert}_{{\mathcal{M}}_{s}^{p}}^{q},$ \ \par 
because  $\displaystyle e_{a}$  is normalised in  $\displaystyle
 H_{s}^{p}.$ \ \par 
So replacing in~(\ref{lCS0}) we get\ \par 
\quad \quad \quad  $\displaystyle \forall r<1,\ \int_{\partial {\mathbb{B}}}{\left\vert{T_{1}(r\zeta
 )}\right\vert ^{p}d\sigma (\zeta )}\leq C(S)^{(l-q)p'}\sum_{a\in
 S}{\left\vert{\lambda _{a}}\right\vert ^{p}}{\left\Vert{\gamma
 _{a}}\right\Vert}_{{\mathcal{M}}_{s}^{p}}^{q}.$ \ \par 
But  $S$  being interpolating we get\ \par 
\quad \quad \quad  $\displaystyle \ {\left\Vert{\gamma _{a}}\right\Vert}_{{\mathcal{M}}_{s}^{p}}^{q}\leq
 C(S)^{q}$ \ \par 
so finally\ \par 
\quad \quad \quad  $\displaystyle \forall r<1,\ \int_{\partial {\mathbb{B}}}{\left\vert{T_{1}(r\zeta
 )}\right\vert ^{p}d\sigma (\zeta )}\leq C(S)^{q}C(S)^{(l-q)p'}{\left\Vert{\lambda
 }\right\Vert}_{\ell ^{p}(S)}^{p}.$ \ \par 
Adding these  $s$  set of sums we get, because  $p'\geq 1,$ \ \par 
\quad \quad \quad  $\displaystyle \ {\left\Vert{R^{j}f}\right\Vert}_{H^{p}}\leq
 s(\max \ A_{q})C(S)^{lp'/p}{\left\Vert{\lambda }\right\Vert}_{\ell
 ^{p}(S)}$ \ \par 
and we are done.  $\displaystyle \hfill\blacksquare $ \ \par 
\ \par 
\quad  Now we shall improve theorem~\ref{cB10}, for all real values
 of  $s\in \lbrack 0,n/p\rbrack ,$ \ \par 

\begin{Thrm}
 ~\label{6IM2}Let  $S$  be interpolating for  $\displaystyle
 {\mathcal{M}}_{s}^{p}$  and suppose that  $\displaystyle p\leq
 2,$  then  $S$  is Carleson in  $\displaystyle H_{r}^{p},\ \forall r\leq s.$ 
\end{Thrm}
\quad  Proof.\ \par 
We know, by theorem~\ref{cB10}, that if  $S$  is interpolating
 for  $\displaystyle {\mathcal{M}}_{s}^{p}$  and if  $\displaystyle
 p\leq 2,$  then  $S$  is Carleson  $\displaystyle H_{s}^{p}.$
  hence we apply theorem~\ref{3CS1} to get the result.  $\displaystyle
 \hfill\blacksquare $ \ \par 
\ \par 
\quad  Arcozzi, Rochberg and Sawyer in~\cite{ArcoRochSaw06} proved,
 in particular, that if  $S$  is interpolating in  $\displaystyle
 B^{p}=B_{0}^{p},$  where  $\displaystyle B_{\sigma }^{p}$  is
 a Besov space of the ball  $\displaystyle {\mathbb{B}},$  then
 we have that  $S$  is Carleson for  $\displaystyle B^{p}.$ 
 In the case  $\displaystyle p=2,\ B^{2}=H_{n/2}^{2},$  we have
 a better result.\ \par 

\begin{Crll}
 Let  $S$  be an interpolating sequence for  $\displaystyle H_{s}^{2}$
  with  $\displaystyle n-2s\leq 1,$  then  $S$  is Carleson for
  $\displaystyle H_{r}^{2},\ \forall r\leq s.$ 
\end{Crll}
\quad \quad  	Proof.\ \par 
We know that  $\displaystyle H_{s}^{2}=B_{\sigma }^{2}$  where
  $\displaystyle B_{\sigma }^{2}$  is the Besov space of the
 ball  $\displaystyle {\mathbb{B}}$  and where  $\displaystyle
 \sigma =\frac{n}{2}-s.$  We know by~\cite{AglMCar02} that for
  $\displaystyle \sigma \leq 1/2,\ B_{\sigma }^{2}$  has Pick
 kernels hence  $S$  interpolating for  $\displaystyle H_{s}^{2}=B_{\sigma
 }^{2}$  implies  $S$  interpolating for  $\displaystyle {\mathcal{M}}_{s}^{2}$
  so we can apply theorem~\ref{6IM2} to get the result.  $\displaystyle
 \hfill\blacksquare $ \ \par 

\subsection{Union of separated interpolating sequences.}
\quad  In the case  $\displaystyle s=0,\ {\mathcal{M}}_{0}^{p}=H^{\infty
 }({\mathbb{B}}),$  the union  $S$  of two interpolating sequences
 in  $\displaystyle H^{\infty }({\mathbb{B}})$  is still interpolating
 in  $\displaystyle H^{\infty }({\mathbb{B}})$  if  $S$  is separated
 by a theorem of Varopoulos~\cite{Varo71}. We shall generalise
 this fact in the next results.\ \par 

\begin{Dfnt}
 A sequence  $S$  is separated in  $\displaystyle {\mathcal{M}}_{s}^{p}$
  if there is a  $c_{S}>0$  such that\par 
\quad \quad \quad  $\displaystyle \forall a,b\neq a\in S,\ \exists m_{a,b}\in {\mathcal{M}}_{s}^{p}::m_{a,b}(a)=1,\
 m_{a,b}(b)=0,\ {\left\Vert{m_{a,b}}\right\Vert}_{{\mathcal{M}}_{s}^{p}}\leq
 c_{S}.$ \par 

\end{Dfnt}

\begin{Dfnt}
 A sequence  $S$  is strongly separated in  $\displaystyle {\mathcal{M}}_{s}^{p}$
  if there is a  $c_{S}>0$  such that\par 
\quad \quad \quad  $\displaystyle \forall a,b\neq a\in S,\ \exists m_{a,b}\in {\mathcal{M}}_{s}^{p}::m_{a,b}(a)=1,\
 m_{a,b}(b)=0,$ \par 
and
\end{Dfnt}
\quad \quad \quad \quad  $\displaystyle \ \forall h\in H_{s}^{p},\ \forall a\in S,\ \exists
 H\in H_{s}^{p},\ {\left\Vert{H}\right\Vert}_{H_{s}^{p}}\leq
 c_{S}{\left\Vert{h}\right\Vert}_{H_{s}^{p}}::\forall b\in S,\
 b\neq a,\ \forall j\leq s,\ \left\vert{R^{j}(m_{a,b}h)}\right\vert
 \leq \left\vert{R^{j}(H)}\right\vert .$ \ \par 
\quad  Clearly the strong separation in  $\displaystyle {\mathcal{M}}_{s}^{p}$
  implies the separation in  $\displaystyle {\mathcal{M}}_{s}^{p}.$ \ \par 

\begin{Dfnt}
 The sequences  $S_{1},\ S_{2}$  are completely separated in
  $\displaystyle {\mathcal{M}}_{s}^{p}$  if there is a  $c_{A}>0$
  such that\par 
\quad \quad \quad  $\displaystyle \forall a\in S_{1},\ \forall b\in S_{2},\ \exists
 m_{a,b}\in {\mathcal{M}}_{s}^{p}::m_{a,b}(a)=1,\ m_{a,b}(b)=0$ \par 
and\par 
\quad \quad \quad  $\displaystyle \forall h\in H_{s}^{p},\ \exists H\in H_{s}^{p},\
 {\left\Vert{H}\right\Vert}_{H_{s}^{p}}\leq c_{A}{\left\Vert{h}\right\Vert}_{H_{s}^{p}}::\forall
 a\in S_{1},b\in S_{2},\ \forall j\leq s,\ \left\vert{R^{j}(m_{a,b}h)}\right\vert
 \leq \left\vert{R^{j}(H)}\right\vert .$ 
\end{Dfnt}
\quad  This time the vector  $H$  does not depend on  $a$  nor on 
 $\displaystyle b.$ \ \par 

\begin{Thrm}
 Let  $\displaystyle S_{1}$  and  $\displaystyle S_{2}$  be two
 interpolating sequences in  ${\mathcal{M}}_{s}^{p},\ s\in {\mathbb{N}}\cap
 \lbrack 0,n/p\rbrack ,$  then  $S:=S_{1}\cup S_{2}$  is an interpolating
 sequence in  $\displaystyle {\mathcal{M}}_{s}^{p}$  if and only
 if  $\displaystyle S_{1}$  and  $\displaystyle S_{2}$  are completely
 separated.
\end{Thrm}
\quad \quad  	Proof.\ \par 
Suppose first that  $S:=S_{1}\cup S_{2}$  is an interpolating
 sequence in  $\displaystyle {\mathcal{M}}_{s}^{p}$  and take
  $\displaystyle \forall a\in S_{1},\ \lambda _{a}=1,\ \forall
 b\in S_{2},\ \lambda _{b}=0.$  Then  $\displaystyle \lambda
 \in \ell ^{\infty }(S)$  hence there is  function  $\displaystyle
 m\in {\mathcal{M}}_{s}^{p}$  such that\ \par 
\quad \quad \quad \quad  $\displaystyle \forall a\in S,\ m(a)=\lambda _{a},$  i.e. $\displaystyle
 \forall a\in S_{1},\ m(a)=1,\ \forall b\in S_{2},\ m(b)=0.$ \ \par 
Now we choose  $\displaystyle \forall a\in S_{1},\ \forall b\in
 S_{2},\ m_{a,b}:=m$  which works and if we set  $\displaystyle
 \forall h\in H_{s}^{p},\ H:=mh$  then we are done with  $\displaystyle
 c_{A}:={\left\Vert{m}\right\Vert}_{{\mathcal{M}}_{s}^{p}},$
  proving that the complete separation is necessary to have 
 $\displaystyle S:=S_{1}\cup S_{2}$  interpolating.\ \par 
\ \par 
\quad  Now we suppose we have the complete separation. As usual we
 suppose  $\displaystyle S_{1},S_{2}$  finite and we set  $\displaystyle
 \lbrace \gamma _{a}\rbrace _{a\in S_{1}}$  the canonical dual
 sequence for  $\displaystyle S_{1}$  in  $\displaystyle {\mathcal{M}}_{s}^{p}$
  and   $\displaystyle \lbrace \Gamma _{b}\rbrace _{b\in S_{2}}$
  the canonical dual sequence for  $\displaystyle S_{2}$  in
  $\displaystyle {\mathcal{M}}_{s}^{p}$  and we want estimates
 not depending on the number of points in  $\displaystyle S_{1}$
  and in  $\displaystyle S_{2}.$ \ \par 
\quad  Take  $\displaystyle b\in S_{2},$  then by hypothesis we have\ \par 
\quad \quad \quad \quad  $\displaystyle \forall a\in S_{1},\ $  $\exists m_{a,b}(z)\in
 {\mathcal{M}}_{s}^{p}\ ::\ m_{a,b}(a)=1,\ m_{a,b}(b)=0.$ \ \par 
\quad \quad  	We set,  $\displaystyle m_{b}:=\sum_{a\in S_{1}}{\gamma _{a}^{l}m_{a,b}}.$
  Then we have  $\displaystyle \forall a\in S_{1},\ m_{b}(a)=1$
  and  $\displaystyle m_{b}(b)=0.$ \ \par 
Because  $\displaystyle S_{1}$  and  $\displaystyle S_{2}$ 
 are finite, the functions  $\displaystyle m_{b}$  are in  $\displaystyle
 {\mathcal{M}}_{s}^{p}$  and they verify\ \par 
\quad \quad \quad \quad  $\displaystyle \forall a\in S_{1},\ \forall b\in S_{2},\ m_{b}(z)={\left\lbrace{
\begin{matrix}
{1}&{if\
 z=a}\cr 
{0}&{if\ z=b}\cr 
\end{matrix}
}\right.}.$  \ \par 
Now we shall glue them by setting\ \par 
\quad \quad \quad  $\displaystyle m:=\sum_{b\in S_{2}}{\Gamma _{b}^{l}(1-m_{b}).}$ \ \par 
We have\ \par 
\quad \quad \quad  $\displaystyle m(a)=0$  if  $a\in S_{1}$  and  $\displaystyle
 m(b)=1$  if  $\displaystyle b\in S_{2},$ \ \par 
hence if  $\displaystyle m\in {\mathcal{M}}_{s}^{p}$  with a
 norm depending only on the constants of interpolation of  $\displaystyle
 S_{1}$  and  $\displaystyle S_{2}$  and of the complete separation,
 then we shall be done because then :\ \par 
\quad \quad \quad  $\displaystyle \forall \lambda ^{1}\in \ell ^{\infty }(S_{1}),\
 \forall \lambda ^{2}\in \ell ^{\infty }(S_{2}),\ \exists m_{j}\in
 {\mathcal{M}}_{s}^{p},\ \forall c\in S_{j},\ m_{j}(c)=\lambda
 _{c}^{j},\ j=1,2\ ;$ \ \par 
now we set, with  $\displaystyle m_{j},\ j=1,2$  as above,\ \par 
\quad \quad \quad  $\displaystyle M:=(1-m)m_{1}+mm_{2}\in {\mathcal{M}}_{s}^{p}$ \ \par 
because  $\displaystyle {\mathcal{M}}_{s}^{p}$  is an algebra,
 and we get\ \par 
\quad \quad \quad  $\displaystyle \forall a\in S_{1},\ M(a)=(1-m(a))m_{1}(a)+m(a)m_{2}(a)=m_{1}(a)=\lambda
 ^{1}_{a}$ \ \par 
and\ \par 
\quad \quad \quad  $\displaystyle \forall b\in S_{2},\ M(b)=(1-m(b))m_{1}(b)+m(b)m_{2}(b)=m_{2}(b)=\lambda
 ^{2}_{b},$ \ \par 
hence  $\displaystyle M$  interpolates the sequence  $\displaystyle
 (\lambda ^{1},\lambda ^{2})$  on  $\displaystyle S_{1}\cup S_{2}.$ \ \par 
\ \par 
\quad  In order to have  $\displaystyle m\in {\mathcal{M}}_{s}^{p},$
  we have to show that \ \par 
\quad \quad \quad  $\displaystyle \forall h\in H_{s}^{p},\ \forall j\leq s,\ R^{j}(mh)\in
 H^{p}$  with control of the norms.\ \par 
We start the same way we did with the linear extension :\ \par 
\quad \quad \quad  $\displaystyle R^{j}(mh)=R^{j}(\sum_{b\in S_{2}}{\Gamma _{b}^{l}(1-m_{b})h})=R^{j}(\sum_{b\in
 S_{2}}{\Gamma _{b}^{l}(1-\sum_{a\in S_{1}}{\gamma _{a}^{l}(z)m_{ab}(z)})h})$
 \ \par 
so we have two terms\ \par 
\quad \quad \quad  $\displaystyle T_{1}=\sum_{b\in S_{2}}{R^{j}(\Gamma _{b}^{l}h)}$ \ \par 
and\ \par 
\quad \quad \quad  $\displaystyle T=\sum_{b\in S_{2}}{R^{j}(\Gamma _{b}^{l}(\sum_{a\in
 S_{1}}{\gamma _{a}^{l}(z)m_{ab}(z)})h)}.$ \ \par 
For  $\displaystyle T_{1}$  we are exactly in the situation
 of the linear extension with  $\displaystyle \lambda _{b}=1,\
 \forall b\in S_{2}$  so we get\ \par 
\quad \quad \quad  $\displaystyle T_{1}\in {\mathcal{M}}_{s}^{p},\ {\left\Vert{T_{1}}\right\Vert}_{{\mathcal{M}}_{s}^{p}}\leq
 C(S_{2})^{l}\max \ _{j=0,...,s}\sum_{q=0}^{j}{\left\vert{A_{q}}\right\vert
 }.$ \ \par 
\ \par 
\quad  Now for  $\displaystyle T$  this is more delicate. First we
 set  $\displaystyle h_{ab}:=m_{ab}h\in H_{s}^{p}$  so we have\ \par 
\quad \quad \quad  $\displaystyle T=\sum_{b\in S_{2}}{R^{j}(\Gamma _{b}^{l}(\sum_{a\in
 S_{1}}{\gamma _{a}^{l}h_{ab}}))}=\ \sum_{a\in S_{1},\ b\in S_{2}}{R^{j}(\Gamma
 _{b}^{l}\gamma _{a}^{l}h_{ab})}.$ \ \par 
\quad  We have to exit the converging factors  $\displaystyle (\gamma
 _{a}\Gamma _{b})^{l-q}$  by the exclusion proposition~\ref{iP2} :\ \par 
\quad \quad \quad  $\displaystyle R^{j}(\Gamma _{b}^{l}\gamma _{a}^{l}h_{ab})=\sum_{q=0}^{j}{A_{q}(\gamma
 _{a}\Gamma _{b})^{l-q}R^{j}((\gamma _{a}\Gamma _{b})^{q}h_{ab})}.$ \ \par 
Because  $\displaystyle s$  is fixed and  $\displaystyle j\leq
 s,$  we have only less than  $s$  terms in the sum and the constants
  $\displaystyle A_{q}$  are bounded, hence, up to a finite sum,
 it is enough to control terms of the forms\ \par 
\quad \quad \quad  $\displaystyle T_{2}:=\sum_{a\in S_{1},b\in S_{2}}{\left\vert{\gamma
 _{a}}\right\vert ^{l-q}\left\vert{\Gamma _{b}}\right\vert ^{l-q}\left\vert{R^{j}((\gamma
 _{a}\Gamma _{b})^{q}h_{ab})}\right\vert }.$ \ \par 
By the Leibnitz formula we get\ \par 
\quad \quad \quad  $\displaystyle R^{j}((\gamma _{a}\Gamma _{b})^{q}h_{ab})=\sum_{k=0}^{j}{C_{j}^{k}R^{k}((\gamma
 _{a}\Gamma _{b})^{q})R^{j-k}(h_{ab})}.$ \ \par 
But the complete separation assumption gives the domination :\ \par 
\quad \quad \quad  $\displaystyle \ \left\vert{R^{j-k}(h_{ab})}\right\vert \leq
 \left\vert{R^{j-k}(H)}\right\vert $  with  $\displaystyle H\in
 H_{s}^{p},\ {\left\Vert{H}\right\Vert}_{H_{s}^{p}}\leq C_{S1}{\left\Vert{h}\right\Vert}_{H_{s}^{p}},$
  and  $H$  independent of  $a,\ b.$ \ \par 
So again up to finite number of terms and bounded constants,
 we are lead to control terms of the form\ \par 
\quad \quad \quad  $\displaystyle T_{3}:=\sum_{a\in S_{1},b\in S_{2}}{\left\vert{\gamma
 _{a}}\right\vert ^{l-q}\left\vert{\Gamma _{b}}\right\vert ^{l-q}\left\vert{R^{k}((\gamma
 _{a}\Gamma _{b})^{q})}\right\vert \left\vert{R^{j-k}(H)}\right\vert }.$ \ \par 
Let  $\displaystyle H':=R^{j-k}(H),$  still independent of 
 $a$  and  $b,$  we have  $\displaystyle H'\in H_{s-j+k}^{p}$  so\ \par 
\quad \quad \quad  $\displaystyle T_{3}:=\sum_{a\in S_{1},b\in S_{2}}{\left\vert{\gamma
 _{a}}\right\vert ^{l-q}\left\vert{\Gamma _{b}}\right\vert ^{l-q}\left\vert{R^{k}((\gamma
 _{a}\Gamma _{b})^{q})}\right\vert \left\vert{H'}\right\vert }.$ \ \par 
Now the inclusion lemma~\ref{iP4} gives\ \par 
\quad \quad \quad  $\displaystyle R^{k}((\gamma _{a}\Gamma _{b})^{q})H'=\sum_{m=0}^{k}{A_{km}R^{m}((\gamma
 _{a}\Gamma _{b})^{q}R^{k-m}(H'))}$ \ \par 
so again it is enough to deal with terms of the form\ \par 
\quad \quad \quad  $\displaystyle T_{4}:=\sum_{a\in S_{1},b\in S_{2}}{\left\vert{\gamma
 _{a}}\right\vert ^{l-q}\left\vert{\Gamma _{b}}\right\vert ^{l-q}\left\vert{R^{m}((\gamma
 _{a}\Gamma _{b})^{q}R^{k-m}(H'))}\right\vert }.$ \ \par 
But  $\displaystyle H':=R^{j-k}(H)$  hence   $\displaystyle
 H'':=R^{k-m}(H')=R^{j-m}(H)$  with  $\displaystyle H\in H_{s}^{p},$
  so  $\displaystyle H''\in H_{s-j+m}^{p},$  with  $\displaystyle
 \ {\left\Vert{H''}\right\Vert}_{H_{s-j+m}^{p}}\leq C_{4}{\left\Vert{H}\right\Vert}_{H_{s}^{p}},$
  still independent of  $a$  and  $\displaystyle b.$  So we have\ \par 
\quad \quad \quad  $\displaystyle T_{4}=\sum_{a\in S_{1},b\in S_{2}}{\left\vert{\gamma
 _{a}}\right\vert ^{l-q}\left\vert{\Gamma _{b}}\right\vert ^{l-q}\left\vert{R^{m}(\gamma
 _{a}^{q}(\Gamma _{b}^{q}H''))}\right\vert }.$ \ \par 
By the Leibnitz formula again we get\ \par 
\quad \quad \quad  $\displaystyle R^{m}(\gamma _{a}^{q}(\Gamma _{b}^{q}H''))=\sum_{k=0}^{m}{C_{k}^{m}R^{k}(\gamma
 _{a}^{q})R^{m-k}(\Gamma _{b}^{q}H'')}$ \ \par 
hence by the finiteness of the number of terms,  it is enough
 to control terms of the form\ \par 
\quad \quad \quad  $\displaystyle T_{5}:=\sum_{a\in S_{1},b\in S_{2}}{\left\vert{\gamma
 _{a}}\right\vert ^{l-q}\left\vert{\Gamma _{b}}\right\vert ^{l-q}\left\vert{R^{k}(\gamma
 _{a}^{q})}\right\vert \left\vert{R^{m-k}(\Gamma _{b}^{q}H'')}\right\vert
 }.$ \ \par 
But the sequence  $\displaystyle S_{2}$  is interpolating for
  $\displaystyle {\mathcal{M}}_{s}^{p}$  hence, still by theorem~\ref{iSH19}
 we have that  $\displaystyle S_{2}$  is interpolating for  $\displaystyle
 {\mathcal{M}}_{r}^{p},\ \forall r\leq s$  so we can apply the
 domination lemma~\ref{iP0} from the appendix to  $\displaystyle
 R^{m-k}(\Gamma _{b}^{q}H'')$  :\ \par 
\quad \quad \quad  $\displaystyle \ \left\vert{R^{m-k}(\Gamma _{b}^{q}H'')}\right\vert
 \leq \frac{1}{N_{2}}\sum_{\mu =1}^{N_{2}}{R^{m-k}(H_{\mu })}$ \ \par 
with  $\displaystyle H_{\mu }\in H_{s-j+m}^{p},\ {\left\Vert{H_{\mu
 }}\right\Vert}_{H_{s-j+m}^{p}}\leq C(S_{2})^{q}{\left\Vert{H''}\right\Vert}_{H_{s-j+m}^{p}}$
  and  $\displaystyle H_{\mu }$  independent of  $a$  and  $\displaystyle
 b.$ \ \par 
Because of the  $\displaystyle \ \frac{1}{N}$  it is enough
 to control uniformly in  $\mu $  terms of the form\ \par 
\quad \quad \quad  $\displaystyle T_{6}:=\sum_{a\in S_{1},b\in S_{2}}{\left\vert{\gamma
 _{a}}\right\vert ^{l-q}\left\vert{\Gamma _{b}}\right\vert ^{l-q}\left\vert{R^{k}(\gamma
 _{a}^{q})}\right\vert \left\vert{R^{m-k}(H_{\mu })}\right\vert }.$ \ \par 
We use the inclusion lemma~\ref{iP4} to get\ \par 
\quad \quad \quad  $\displaystyle R^{k}(\gamma _{a}^{q})R^{m-k}(H_{\mu })=\sum_{r=0}^{k}{A_{kr}R^{r}(\gamma
 _{a}^{q}R^{k-r}(R^{m-k}(H_{\mu })))}=\sum_{r=0}^{k}{A_{kr}R^{r}(\gamma
 _{a}^{q}R^{m-r}(H_{\mu }))}.$ \ \par 
So it remains to control terms of the form\ \par 
\quad \quad \quad  $\displaystyle T_{7}:=\sum_{a\in S_{1},b\in S_{2}}{\left\vert{\gamma
 _{a}}\right\vert ^{l-q}\left\vert{\Gamma _{b}}\right\vert ^{l-q}\left\vert{R^{r}(\gamma
 _{a}^{q}R^{m-r}(H_{\mu }))}\right\vert }.$ \ \par 
Set  $\displaystyle V_{\mu }:=R^{m-r}(H_{\mu })\ ;$  because
  $\displaystyle H_{\mu }\in H_{s-j+m}^{p},$  we have that  $\displaystyle
 V_{\mu }\in H_{s-j+r}^{p}$  with control of its norm,\ \par 
\quad \quad \quad  $\displaystyle \ {\left\Vert{V_{\mu }}\right\Vert}_{H_{s-j+m}^{p}}\leq
 C_{7}{\left\Vert{H_{\mu }}\right\Vert}_{H_{s-j+m}^{p}}.$ \ \par 
So we have\ \par 
\quad \quad \quad  $\displaystyle T_{7}=\sum_{a\in S_{1},b\in S_{2}}{\left\vert{\gamma
 _{a}}\right\vert ^{l-q}\left\vert{\Gamma _{b}}\right\vert ^{l-q}\left\vert{R^{r}(\gamma
 _{a}^{q}V_{\mu })}\right\vert }.$ \ \par 
Now we shall use that  $\displaystyle S_{1}$  is interpolating
 for  $\displaystyle {\mathcal{M}}_{s}^{p}$  hence, still by
 theorem~\ref{iSH19} we have that  $\displaystyle S_{1}$  is
 interpolating for  $\displaystyle {\mathcal{M}}_{r}^{p},\ \forall
 r\leq s$  so we can apply the domination lemma~\ref{iP0} :\ \par 
\quad \quad \quad  $\displaystyle \ \left\vert{R^{r}(\gamma _{a}^{q}V_{\mu })}\right\vert
 \leq \frac{1}{N_{1}}\sum_{\nu =1}^{N_{1}}{\left\vert{R^{r}(H_{\nu
 \mu })}\right\vert }$ \ \par 
with  $\displaystyle 1\leq \nu \leq N,\ {\left\Vert{H_{\nu \mu
 }}\right\Vert}_{H_{s-j+r}^{p}}\leq C(S_{1})^{q}{\left\Vert{V_{\mu
 }}\right\Vert}_{H_{s-j+r}^{p}}$  and  $\displaystyle H_{\nu
 \mu }$  not depending on  $\displaystyle a\in S_{1}$  nor on
  $\displaystyle b\in S_{2}.$ \ \par 
So, because of the  $\displaystyle \ \frac{1}{N_{1}}$  we need
 to control uniformly in  $\nu ,$  terms of the form\ \par 
\quad \quad \quad  $\displaystyle T_{8}:=\sum_{a\in S_{1},b\in S_{2}}{\left\vert{\gamma
 _{a}}\right\vert ^{l-q}\left\vert{\Gamma _{b}}\right\vert ^{l-q}\left\vert{R^{r}(H_{\nu
 \mu })}\right\vert }.$ \ \par 
But now we use proposition~\ref{HA0} which tells us for  $\displaystyle
 l-q\geq 2$  :\ \par 
\quad \quad \quad  $\displaystyle \ \sum_{a\in S_{1}}{\left\vert{\gamma _{a}}\right\vert
 ^{l-q}}\leq C(S_{1})^{l-q}$ \ \par 
and the same for  $\displaystyle S_{2}$ \ \par 
\quad \quad \quad  $\displaystyle \ \sum_{b\in S_{2}}{\left\vert{\Gamma _{b}}\right\vert
 ^{l-q}}\leq C(S_{2})^{l-q}.$ \ \par 
Hence porting in  $\displaystyle T_{8}$ \ \par 
\quad \quad \quad  $\displaystyle T_{8}\leq (C(S_{1})C(S_{2}))^{l-q}\left\vert{R^{r}(H_{\nu
 \mu })}\right\vert .$ \ \par 
Now taking the  $\displaystyle H^{p}$  norm we get\ \par 
\quad \quad \quad  $\displaystyle \ {\left\Vert{T_{8}}\right\Vert}_{H^{p}}\leq
 (C(S_{1})C(S_{2}))^{l-q}{\left\Vert{R^{r}(H_{\nu \mu })}\right\Vert}_{H^{p}}$
 \ \par 
but recall that\ \par 
\quad \quad \quad  $\displaystyle \ {\left\Vert{H_{\nu \mu }}\right\Vert}_{H_{s-j+r}^{p}}\leq
 C(S_{1})^{q}{\left\Vert{V_{\mu }}\right\Vert}_{H_{s-j+r}^{p}}\Rightarrow
 {\left\Vert{R^{r}(H_{\nu \mu })}\right\Vert}_{H^{p}}\leq C(S_{1})^{q}{\left\Vert{V_{\mu
 }}\right\Vert}_{H_{s-j+r}^{p}}$ \ \par 
and\ \par 
\quad \quad \quad  $\displaystyle \ {\left\Vert{V_{\mu }}\right\Vert}_{H_{s-j+m}^{p}}\leq
 C_{7}{\left\Vert{H_{\mu }}\right\Vert}_{H_{s-j+m}^{p}}$ \ \par 
and\ \par 
\quad \quad \quad  $\displaystyle \ {\left\Vert{H_{\mu }}\right\Vert}_{H_{s-j+m}^{p}}\leq
 C(S_{2}){\left\Vert{H''}\right\Vert}_{H_{s-j+m}^{p}}$ \ \par 
and\ \par 
\quad \quad \quad  $\displaystyle \ {\left\Vert{H''}\right\Vert}_{H_{s-j+m}^{p}}\leq
 C_{4}{\left\Vert{H}\right\Vert}_{H_{s}^{p}}$ \ \par 
and\ \par 
\quad \quad \quad  $\displaystyle \ {\left\Vert{H}\right\Vert}_{H_{s}^{p}}\leq
 C_{A}{\left\Vert{h}\right\Vert}_{H_{s}^{p}},$ \ \par 
so concatenating we get\ \par 
\quad \quad \quad  $\displaystyle \ {\left\Vert{T_{8}}\right\Vert}_{H^{p}}\leq
 C_{A}C_{4}C_{7}(C(S_{1})C(S_{2}))^{l}{\left\Vert{h}\right\Vert}_{H_{s}^{p}}$
 \ \par 
and the proof is complete.  $\displaystyle \hfill\blacksquare $ \ \par 

\section{Dual boundedness and interpolating sequences in  $\displaystyle
 H_{s}^{p}.$ ~\label{iSH13}}
\quad  The Sobolev embedding theorem gives, in  ${\mathbb{R}}^{n},$ \ \par 
\quad \quad \quad \quad \quad 	 $\displaystyle f\in W_{s}^{p}({\mathbb{R}}^{n})\Rightarrow
 f\in L^{q}({\mathbb{R}}^{n}),\ \frac{1}{q}=\frac{1}{p}-\frac{s}{n}.$ \ \par 
Here we are on the manifold  $\displaystyle \partial {\mathbb{B}}$
  which is of dimension  $\displaystyle 2n-1,$  and with complex
 tangential derivatives of order  $\displaystyle 2s$  and normal
 conjugate  derivative of order  $\displaystyle s.$ \ \par 
\quad  Thanks to Folland and Stein~\cite{FolStein74}, theorem 2, which
 we iterate and which we apply with  $\displaystyle \alpha =0$
  or by use of Romanovskii~\cite{Romanovskii05}, theorem 7, we
 have a Sobolev anisotropic embedding in Heisenberg group, which
 is also a representation of the boundary of the ball  $\displaystyle
 {\mathbb{B}}\ :$ \ \par 
\quad \quad \quad \quad  	\begin{equation}  f\in H_{s}^{p}({\mathbb{B}})\Rightarrow f\in
 H^{q}({\mathbb{B}}),\ \frac{1}{q}=\frac{1}{p}-\frac{s}{n}.\label{SH5}\end{equation}\
 \par 

\begin{Thrm}
 ~\label{SH3}Let  $\displaystyle S\subset {\mathbb{B}}$  be a
 dual  bounded sequence for  $\displaystyle H_{s}^{p}({\mathbb{B}}),$
  then  $S$  is dual bounded for  $\displaystyle H^{q}({\mathbb{B}})$
  with  $\displaystyle \ \frac{1}{q}=\frac{1}{p}-\frac{s}{n}.$ 
\end{Thrm}
\quad \quad  	Proof.\ \par 
Saying  $S$  dual bounded in  $\displaystyle H_{s}^{p}({\mathbb{B}})$
  means, with  $\displaystyle k_{s,a}$  the reproducing kernel
 for  $\displaystyle H_{s}^{2}({\mathbb{B}}),$ \ \par 
\quad \quad \quad \quad \quad 	 $\displaystyle \exists C>0,\ \forall a\in S,\ \exists \rho
 _{a}\in H_{s}^{p}({\mathbb{B}})::\rho _{a}(b)=\delta _{a,b}{\left\Vert{k_{s,a}}\right\Vert}_{s,p},\
 {\left\Vert{\rho _{a}}\right\Vert}_{s,p}\leq C.$ \ \par 
But by use of anisotropic Sobolev embeddings~(\ref{SH5}) we
 have, with  $\displaystyle \ \frac{1}{q}=\frac{1}{p}-\frac{s}{n},$ \ \par 
\quad \quad \quad \quad \quad 	 $\displaystyle \exists C>0,\ f\in H_{s}^{p}({\mathbb{B}})\Rightarrow
 f\in H^{q}({\mathbb{B}}),\ {\left\Vert{f}\right\Vert}_{q}\leq
 C{\left\Vert{f}\right\Vert}_{s,p}.$ \ \par 
\quad \quad  	On the other hand we have\ \par 
\quad \quad \quad \quad \quad 	 $\displaystyle \ {\left\Vert{k_{s,a}}\right\Vert}_{s,p}=(1-\left\vert{a}\right\vert
 ^{2})^{s-\frac{n}{p'}},$ \ \par 
hence with  $\displaystyle \ \frac{1}{q}=\frac{1}{p}-\frac{s}{n},$
  we get\ \par 
\quad \quad \quad \quad \quad 	 $\displaystyle \ {\left\Vert{k_{0,a}}\right\Vert}_{q'}=(1-\left\vert{a}\right\vert
 ^{2})^{-\frac{n}{q}}=(1-\left\vert{a}\right\vert ^{2})^{-n(\frac{1}{p}-\frac{s}{n})}=(1-\left\vert{a}\right\vert
 ^{2})^{s-\frac{n}{p}}={\left\Vert{k_{s,a}}\right\Vert}_{s,p'}.$ \ \par 
So we have a dual sequence for  $S$  in  $\displaystyle H^{q}({\mathbb{B}}),$
  namely  $\displaystyle \lbrace \rho _{a}\rbrace _{a\in S}$
  itself, doing\ \par 
\quad \quad \quad \quad \quad 	 $\displaystyle \exists C>0,\ \forall a\in S,\ \exists \rho
 _{a}\in H^{q}({\mathbb{B}})::\rho _{a}(b)\simeq \delta _{a,b}{\left\Vert{k_{a}}\right\Vert}_{q},\
 {\left\Vert{\rho _{a}}\right\Vert}_{q}\leq C,$ \ \par 
which means that  $S$  is dual bounded in  $\displaystyle H^{q}({\mathbb{B}}).$
   $\displaystyle \hfill\blacksquare $ \ \par 
\ \par 
\quad \quad   $S$  interpolating for  $\displaystyle H_{s}^{p}({\mathbb{B}})$  means\ \par 
\quad \quad \quad \quad 	 $\displaystyle \forall \lambda \in l^{p}(S),\ \exists f\in
 H_{s}^{p}({\mathbb{B}})::\forall a\in S,\ f(a)=\lambda _{a}{\left\Vert{k_{a}}\right\Vert}_{s,p'}$
 \ \par 
so we have  $\displaystyle f\in H^{q}({\mathbb{B}}),\ {\left\Vert{f}\right\Vert}_{q}\leq
 C{\left\Vert{f}\right\Vert}_{s,p}$  such that\ \par 
\quad \quad \quad \quad \quad 	 $\displaystyle \forall a\in S,\ f(a)=\lambda _{a}{\left\Vert{k_{a}}\right\Vert}_{q}$
 \ \par 
hence we interpolate  $l^{p}(S)$  sequences in  $\displaystyle
 H^{q}({\mathbb{B}})$  for  $\displaystyle \ \frac{1}{q}=\frac{1}{p}-\frac{s}{n},$
  but {\sl not}  $\displaystyle \ell ^{q}(S)$  sequences so this
 is {\sl not} the  $\displaystyle H^{q}({\mathbb{B}})$  interpolation !\ \par 

\begin{Crll}
 ~\label{7BS5}Let  $\displaystyle S\subset {\mathbb{B}}$  be
 a dual  bounded sequence for  $\displaystyle H_{s}^{p}({\mathbb{B}}),$
  then  $S$  is Carleson for  $\displaystyle H^{p}({\mathbb{B}}).$ 
\end{Crll}
\quad  Proof.\ \par 
This is exactly the result in~\cite{AmarWirtBoule07}, because
  $S$  is dual bounded in  $\displaystyle H^{q}({\mathbb{B}})$
  hence Carleson for all  $\displaystyle H^{r}({\mathbb{B}}).$ \ \par 
 $\displaystyle \hfill\blacksquare $ \ \par 
\ \par 
\quad  The first structural hypothesis (see ~\cite{AmarExtInt06}) is
 true for these spaces :\ \par 

\begin{Lmm}
 ~\label{dBS1}we have\par 
\quad \quad \quad \quad \quad 	 $\displaystyle \ \forall r>1,\ k_{a}(a)={\left\Vert{k_{a}}\right\Vert}_{H_{s}^{2}}^{2}={\left\Vert{k_{a}}\right\Vert}_{H_{s}^{r}}{\left\Vert{k_{a}}\right\Vert}_{H_{s}^{r'}}.$
 
\end{Lmm}
\quad  Proof.\ \par 
We have\ \par 
\quad \quad \quad  $\displaystyle \ {\left\Vert{k_{a}}\right\Vert}_{H_{s}^{p}}=(1-\left\vert{a}\right\vert
 ^{2})^{s-n/p'}$ \ \par 
and\ \par 
\quad \quad \quad  $\displaystyle \ {\left\Vert{k_{a}}\right\Vert}_{H_{s}^{r}}{\left\Vert{k_{a}}\right\Vert}_{H_{s}^{r'}}=(1-\left\vert{a}\right\vert
 ^{2})^{s-n/r'}(1-\left\vert{a}\right\vert ^{2})^{s-n/r}=(1-\left\vert{a}\right\vert
 ^{2})^{2s-n}$ \ \par 
hence\ \par 
\quad \quad \quad  $\displaystyle \ {\left\Vert{k_{a}}\right\Vert}_{H_{s}^{2}}=(1-\left\vert{a}\right\vert
 ^{2})^{s-n/2}\Rightarrow {\left\Vert{k_{a}}\right\Vert}_{H_{s}^{2}}^{2}=(1-\left\vert{a}\right\vert
 ^{2})^{2s-n}$ \ \par 
which proves the lemma.  $\displaystyle \hfill\blacksquare $ \ \par 
\quad  Recall that we shall say that  $S$  is a  $\displaystyle H_{s}^{p}$
  weighted interpolating sequence for the weight  $\displaystyle
 w=\lbrace w_{a}\rbrace _{a\in S}$  if \ \par 
\quad \quad \quad  $\displaystyle \forall \lambda \in \ell ^{p}(S),\ \exists f\in
 H_{s}^{p}::\forall a\in S,\ f(a)=\lambda _{a}w_{a}{\left\Vert{k_{a}}\right\Vert}_{H_{s}^{p'}}.$
 \ \par 
\quad  By use of lemma~\ref{dBS1} we get\ \par 

\begin{Thrm}
 Let  $\displaystyle p>1$  and suppose that  $S$  is dual bounded
 in  $\displaystyle H_{s}^{p},$  then  $S$  is a  $\displaystyle
 H_{s}^{1}$  weighted interpolating sequence for the weight 
 $\displaystyle \lbrace (1-\left\vert{a}\right\vert ^{2})^{s}\rbrace
 _{a\in S}$  with the bounded linear extension property.
\end{Thrm}
\quad  Proof.\ \par 
Consider the dual sequence  $\displaystyle \lbrace \rho _{a}\rbrace
 _{a\in S}$  in  $\displaystyle H_{s}^{p},$  given by the dual
 boundedness, it verifies\ \par 
\quad \quad \quad  $\displaystyle \exists C>0,\ \forall a\in S,\ {\left\Vert{\rho
 _{a}}\right\Vert}_{H_{s}^{p}}\leq C,\ \forall b\in S,\ \rho
 _{a}(b)=\delta _{ab}{\left\Vert{k_{a}}\right\Vert}_{H_{s}^{p'}}.$ \ \par 
Let, for  $\displaystyle \lambda \in \ell ^{1}(S),$ \ \par 
\quad \quad \quad \quad  $\displaystyle h:=\sum_{a\in S}{\lambda _{a}\rho _{a}\frac{k_{a}}{{\left\Vert{k_{a}}\right\Vert}_{H_{s}^{p'}}}},$
 \ \par 
we have\ \par 
\quad \quad \quad  $\displaystyle h(a)=\lambda _{a}k_{a}(a)=\lambda _{a}(1-\left\vert{a}\right\vert
 ^{2})^{2s-n}$ \ \par 
which is the right value. As its norm we get\ \par 
\quad \quad \quad  $\displaystyle \ {\left\Vert{h}\right\Vert}_{H_{s}^{1}}\leq
 \sum_{a\in S}{\left\vert{\lambda _{a}}\right\vert {\left\Vert{\rho
 _{a}\frac{k_{a}}{{\left\Vert{k_{a}}\right\Vert}_{H_{s}^{p'}}}}\right\Vert}_{H_{s}^{1}}}\leq
 C{\left\Vert{\lambda }\right\Vert}_{\ell ^{1}},$ \ \par 
because, using proposition~\ref{S2} we get\ \par 
\quad \quad \quad  $\displaystyle \ {\left\Vert{\rho _{a}\frac{k_{a}}{{\left\Vert{k_{a}}\right\Vert}_{H_{s}^{p'}}}}\right\Vert}_{H_{s}^{1}}\leq
 C_{s}{\left\Vert{\rho _{a}}\right\Vert}_{H_{s}^{p}}{\left\Vert{\frac{k_{a}}{{\left\Vert{k_{a}}\right\Vert}_{H_{s}^{p'}}}}\right\Vert}_{H_{s}^{p'}}\leq
 C_{s}{\left\Vert{\rho _{a}}\right\Vert}_{H_{s}^{p}}\leq C_{s}C.$
   $\displaystyle \hfill\blacksquare $ \ \par 
\ \par 
\quad  For the second structural hypothesis we have\ \par 

\begin{Lmm}
 ~\label{dBS4}Let  $\displaystyle p,r\in \lbrack 1,\infty \rbrack
 $  and  $q$  such that  $\displaystyle \ \frac{1}{r}=\frac{1}{p}+\frac{1}{q}$
  then we have\par 
\quad \quad \quad  $\displaystyle \ {\left\Vert{k_{a}}\right\Vert}_{H_{s}^{r'}}\simeq
 (1-\left\vert{a}\right\vert ^{2})^{-s}{\left\Vert{k_{a}}\right\Vert}_{H_{s}^{p'}}{\left\Vert{k_{a}}\right\Vert}_{H_{s}^{q'}}.$
 
\end{Lmm}
\quad  Proof.\ \par 
We have\ \par 
\quad \quad \quad  $\displaystyle \ {\left\Vert{k_{a}}\right\Vert}_{H_{s}^{r'}}\simeq
 (1-\left\vert{a}\right\vert ^{2})^{s-n/r}$ \ \par 
and\ \par 
\quad \quad \quad  $\displaystyle \ {\left\Vert{k_{a}}\right\Vert}_{H_{s}^{p'}}{\left\Vert{k_{a}}\right\Vert}_{H_{s}^{q'}}\simeq
 (1-\left\vert{a}\right\vert ^{2})^{s-n/p}(1-\left\vert{a}\right\vert
 ^{2})^{s-n/q}=(1-\left\vert{a}\right\vert ^{2})^{2s-n/r}$ \ \par 
hence\ \par 
\quad \quad \quad  $\displaystyle \ {\left\Vert{k_{a}}\right\Vert}_{H_{s}^{r'}}\simeq
 (1-\left\vert{a}\right\vert ^{2})^{-s}{\left\Vert{k_{a}}\right\Vert}_{H_{s}^{p'}}{\left\Vert{k_{a}}\right\Vert}_{H_{s}^{q'}}$
 \ \par 
which proves the lemma.  $\displaystyle \hfill\blacksquare $ \ \par 
\ \par 
\ \par 
\quad  Now we are in position to get an analogous result to theorem
 6.1 in~\cite{AmarExtInt06} by an analogous proof.\ \par 

\begin{Thrm}
 Let  $S$  be a sequence of points in  $\displaystyle {\mathbb{B}}$
  such that, with  $\displaystyle \ \frac{1}{r}=\frac{1}{p}+\frac{1}{q},$
  and  $\displaystyle p\leq 2,$  \par 
\quad \quad \quad 	 $\displaystyle \bullet $   $S$  is dual bounded in  $\displaystyle
 H_{s}^{p}.$ \par 
\quad \quad \quad  $\displaystyle \bullet $    $S$  is Carleson in  $\displaystyle
 H_{s}^{q}({\mathbb{B}}).$ \par 
Then  $S$  is a  $\displaystyle H_{s}^{r}$  weighted interpolating
 sequence for the weight  $\displaystyle \lbrace (1-\left\vert{a}\right\vert
 ^{2})^{s}\rbrace _{a\in S}$  with the bounded linear extension property.
\end{Thrm}
\quad  Proof.\ \par 
Consider the dual sequence  $\displaystyle \lbrace \rho _{a}\rbrace
 _{a\in S}$  in  $\displaystyle H_{s}^{p},$  given by the hypothesis,
 it verifies\ \par 
\quad \quad \quad  $\displaystyle \exists C>0,\ \forall a\in S,\ {\left\Vert{\rho
 _{a}}\right\Vert}_{H_{s}^{p}}\leq C,\ \forall b\in S,\ \rho
 _{a}(b)=\delta _{ab}{\left\Vert{k_{a}}\right\Vert}_{H_{s}^{p'}}.$ \ \par 
Now we set \ \par 
\quad \quad \quad \quad \quad 	 $\displaystyle h(z):=\sum_{a\in S}{\lambda _{a}\rho _{a}(1-\left\vert{a}\right\vert
 ^{2})^{s}\frac{k_{a}(z)}{{\left\Vert{k_{a}}\right\Vert}_{H_{s}^{p'}}{\left\Vert{k_{a}}\right\Vert}_{H_{s}^{r}}}}.$
 \ \par 
We get\ \par 
\quad \quad \quad \quad \quad 	 $\displaystyle \forall a\in S,\ h(a)=\lambda _{a}(1-\left\vert{a}\right\vert
 ^{2})^{s}\rho _{a}(a)\frac{k_{a}(a)}{{\left\Vert{k_{a}}\right\Vert}_{H_{s}^{p'}}{\left\Vert{k_{a}}\right\Vert}_{H_{s}^{r}}}$
 \ \par 
and using the first structural hypothesis, lemma~\ref{dBS1}, we get\ \par 
\quad \quad \quad  $\displaystyle h(a)=\lambda _{a}(1-\left\vert{a}\right\vert
 ^{2})^{s}{\left\Vert{k_{a}}\right\Vert}_{H_{s}^{r'}},$ \ \par 
hence  $h$  interpolates the correct values.\ \par 
\quad  Clearly  $h$  is linear in  $\displaystyle \lambda ,$  and it
 remains to estimate the norm of  $\displaystyle h.$ \ \par 
\ \par 
Proof of the estimates.\ \par 
\quad  In order to do this, we proceed as in~\cite{AmarExtInt06} :\ \par 
let  $\lbrace \epsilon _{a}\rbrace _{a\in S}\in {\mathcal{R}}(S)$
  be a Rademacher sequence of random variables, we set, with
  $\forall a\in S,\ \lambda _{a}=\mu _{a}\nu _{a},\ \mu \in l^{p},\
 \nu \in \ell ^{q}$  explicitly :\ \par 
\quad \quad \quad  $\displaystyle \mu _{a}:=\frac{\lambda _{a}}{\left\vert{\lambda
 _{a}}\right\vert ^{\alpha }},\ \nu _{a}:=\left\vert{\lambda
 _{a}}\right\vert ^{\alpha }$  with  $\displaystyle \alpha =\frac{r}{q}\
 ;$ \ \par 
then we get  $\displaystyle \lambda _{a}=\mu _{a}\nu _{a},\
 \mu \in l^{p},\ \nu \in \ell ^{q}$  and  $\displaystyle \ {\left\Vert{\mu
 }\right\Vert}_{\ell ^{p}}^{p}={\left\Vert{\nu }\right\Vert}_{\ell
 ^{q}}^{q}={\left\Vert{\lambda }\right\Vert}_{\ell ^{r}}^{r}\Rightarrow
 {\left\Vert{\lambda }\right\Vert}_{\ell ^{r}}={\left\Vert{\nu
 }\right\Vert}_{\ell ^{q}}{\left\Vert{\mu }\right\Vert}_{\ell ^{q}}.$ \ \par 
\quad  Now set\ \par 
\quad \quad \quad \quad \quad 	 $\displaystyle f(\epsilon ,z):=\sum_{a\in S}{\mu _{a}\epsilon
 _{a}\rho _{a}(z)}\ ;\ g(\epsilon ,z):=\sum_{a\in S}{\nu _{a}(1-\left\vert{a}\right\vert
 ^{2})^{s}\epsilon _{a}\frac{k_{a}(z)}{{\left\Vert{k_{a}}\right\Vert}_{H_{s}^{p'}}{\left\Vert{k_{a}}\right\Vert}_{H_{s}^{r}}}}.$
 \ \par 
We have	 	 ${\mathbb{E}}(fg)=h$  hence\ \par 
\quad \quad \quad  $\displaystyle \ {\left\Vert{h}\right\Vert}_{H_{s}^{r}}={\left\Vert{{\mathbb{E}}(fg)}\right\Vert}_{H_{s}^{r}}.$
 \ \par 
By lemma~\ref{S1} we get  $\displaystyle \ {\left\Vert{{\mathbb{E}}(fg)}\right\Vert}_{H_{s}^{r}}^{r}\leq
 {\mathbb{E}}({\left\Vert{fg}\right\Vert}_{H_{s}^{r}}^{r})$ 
 and by proposition~\ref{S2} we get  $\displaystyle \ {\left\Vert{fg}\right\Vert}_{H_{s}^{r}}\leq
 C_{s}{\left\Vert{f}\right\Vert}_{H_{s}^{p}}{\left\Vert{g}\right\Vert}_{H_{s}^{q}},$
  so\ \par 
\quad \quad \quad  \begin{equation}  \ {\left\Vert{{\mathbb{E}}(fg)}\right\Vert}_{H_{s}^{r}}^{r}\leq
 {\mathbb{E}}({\left\Vert{fg}\right\Vert}_{H_{s}^{r}}^{r})\leq
 C_{s}^{r}{\mathbb{E}}({\left\Vert{f}\right\Vert}_{H_{s}^{p}}^{r}{\left\Vert{g}\right\Vert}_{H_{s}^{q}}^{r}).\label{dBS2}\end{equation}\
 \par 
\quad  Set  $\displaystyle \gamma _{a}:=\frac{(1-\left\vert{a}\right\vert
 ^{2})^{s}{\left\Vert{k_{a}}\right\Vert}_{H_{s}^{q}}}{{\left\Vert{k_{a}}\right\Vert}_{H_{s}^{p'}}{\left\Vert{k_{a}}\right\Vert}_{H_{s}^{r}}}$
  we have  $\displaystyle g(\epsilon ,z):=\sum_{a\in S}{\nu _{a}\epsilon
 _{a}\gamma _{a}\frac{k_{a}(z)}{{\left\Vert{k_{a}}\right\Vert}_{H_{s}^{q}}}}.$
 \ \par 
Because  $S$  is  $q$  Carleson by assumption we get\ \par 
\quad \quad \quad  $\displaystyle \exists C>0::\forall \epsilon ,\ {\left\Vert{g}\right\Vert}_{H_{s}^{q}}^{q}\leq
 c_{q}^{q}\sum_{a\in S}{\left\vert{\nu _{a}}\right\vert ^{q}\gamma
 _{a}^{q}.}$ \ \par 
Let us compute  $\gamma _{a}\ :$ \ \par 
\quad \quad \quad  $\displaystyle \gamma _{a}:=\frac{(1-\left\vert{a}\right\vert
 ^{2})^{s}{\left\Vert{k_{a}}\right\Vert}_{H_{s}^{q}}}{{\left\Vert{k_{a}}\right\Vert}_{H_{s}^{p'}}{\left\Vert{k_{a}}\right\Vert}_{H_{s}^{r}}}=\frac{(1-\left\vert{a}\right\vert
 ^{2})^{s}(1-\left\vert{a}\right\vert ^{2})^{s-n/q'}}{(1-\left\vert{a}\right\vert
 ^{2})^{s-n/p}(1-\left\vert{a}\right\vert ^{2})^{s-n/r'}}=1$ \ \par 
because  $\displaystyle \ \frac{1}{r}=\frac{1}{p}+\frac{1}{q}.$
  We see here that the weight  $\displaystyle (1-\left\vert{a}\right\vert
 ^{2})^{s}$  compensates the second structural hypothesis, which
 is given by lemma~\ref{dBS4}.\ \par 
So we get\ \par 
\quad \quad \quad  $\displaystyle \ {\left\Vert{g}\right\Vert}_{H_{s}^{q}}\leq
 c_{q}{\left\Vert{\nu }\right\Vert}_{\ell ^{q}(S)}.$ \ \par 
Putting this in~(\ref{dBS2}) we get\ \par 
\quad  \begin{equation}  \ {\left\Vert{h}\right\Vert}_{H_{s}^{r}}^{r}={\left\Vert{{\mathbb{E}}(fg)}\right\Vert}_{H_{s}^{r}}^{r}\leq
 C_{s}^{r}{\mathbb{E}}({\left\Vert{f}\right\Vert}_{H_{s}^{p}}^{r}{\left\Vert{g}\right\Vert}_{H_{s}^{q}}^{r})\leq
 C_{s}^{r}c_{q}^{r}{\left\Vert{\nu }\right\Vert}_{\ell ^{q}(S)}^{r}{\mathbb{E}}({\left\Vert{f}\right\Vert}_{H_{s}^{p}}^{r}).\label{dBS3}\end{equation}\
 \par 
Now we use that  $\displaystyle p\leq 2$  to get, because  $\displaystyle
 r<p\leq 2,$ \ \par 
\quad \quad \quad  ${\mathbb{E}}({\left\Vert{f}\right\Vert}_{H_{s}^{p}}^{r})\leq
 ({\mathbb{E}}({\left\Vert{f}\right\Vert}_{H_{s}^{p}}^{2}))^{r/2}$ \ \par 
and  $\displaystyle H_{s}^{p}$  is of type  $p$  so, with  $\displaystyle
 f=\sum_{a\in S}{\mu _{a}\epsilon _{a}\rho _{a}(z)},$  we get\ \par 
\quad \quad \quad  $\displaystyle ({\mathbb{E}}({\left\Vert{f}\right\Vert}_{H_{s}^{p}}^{2}))^{1/2}\leq
 T_{p}(\sum_{a\in S}{\left\vert{\mu _{a}}\right\vert ^{p}{\left\Vert{\rho
 _{a}}\right\Vert}_{H_{s}^{p}}^{p}})^{1/p}$ \ \par 
hence, because  $\displaystyle \forall a\in S,\ {\left\Vert{\rho
 _{a}}\right\Vert}_{H_{s}^{p}}\leq C,$  we get\ \par 
\quad \quad \quad  $\displaystyle ({\mathbb{E}}({\left\Vert{f}\right\Vert}_{H_{s}^{p}}^{2}))^{1/2}\leq
 T_{p}C{\left\Vert{\mu }\right\Vert}_{\ell ^{q}}.$ \ \par 
Putting this in~(\ref{dBS3}) we get\ \par 
\quad \quad \quad  $\displaystyle \ {\left\Vert{{\mathbb{E}}(fg)}\right\Vert}_{H_{s}^{r}}\leq
 C_{s}c_{q}T_{p}C{\left\Vert{\nu }\right\Vert}_{\ell ^{q}(w_{s})}{\left\Vert{\mu
 }\right\Vert}_{\ell ^{q}}.$ \ \par 
Hence finally\ \par 
\quad \quad \quad  $\displaystyle \ {\left\Vert{h}\right\Vert}_{H_{s}^{r}}={\left\Vert{{\mathbb{E}}(fg)}\right\Vert}_{H_{s}^{r}}\leq
 C_{s}c_{q}T_{p}C{\left\Vert{\nu }\right\Vert}_{\ell ^{q}}{\left\Vert{\mu
 }\right\Vert}_{\ell ^{q}}$ \ \par 
which proves the theorem because  $\displaystyle \ {\left\Vert{\lambda
 }\right\Vert}_{\ell ^{r}}={\left\Vert{\nu }\right\Vert}_{\ell
 ^{q}}{\left\Vert{\mu }\right\Vert}_{\ell ^{q}}.$   $\displaystyle
 \hfill\blacksquare $ \ \par 

\section{Appendix.}

\subsection{Technical lemmas.}
\quad  With the notations of section~\ref{5HA1} and ~\ref{6IM3}, let
  $\displaystyle S=\lbrace a_{1},...,\ a_{N}\rbrace ,$  fix 
 $\displaystyle a\in S$  and set  $\gamma =\gamma _{a}$  to ease
 the notations. Also if  $\displaystyle f\in H_{s}^{p},$  we
 set  $\displaystyle f^{(j)}:=R^{j}f.$ \ \par 

\begin{Lmm}
 ~\label{J10}We have, with  $m:=\min \ (l,j),\ \forall j,l\in
 {\mathbb{N}},$ \par 
\quad \quad \quad \quad  	\begin{equation}  R^{j}(\gamma ^{l}h)=\gamma ^{l}F_{0,j}(z)+l\gamma
 ^{l-1}F_{1,j}(z)+\cdot \cdot \cdot +l(l-1)\cdot \cdot \cdot
 (l-m+1)\gamma ^{l-m}F_{m,j}(z)\label{J8}\end{equation}\par 
where the functions  $\displaystyle F_{k,j}(z)$  do not depend
 on  $\displaystyle l$  for  $\displaystyle k\leq m.$ 
\end{Lmm}
\quad \quad  	Proof.\ \par 
By induction on  $\displaystyle j.$  For  $\displaystyle j=1$  we have :\ \par 
\quad \quad \quad  $\displaystyle R(\gamma ^{l}h)=\gamma ^{l}h^{(1)}+l\gamma ^{l-1}\gamma
 ^{(1)}h,$ \ \par 
hence\ \par 
\quad \quad \quad \quad  $\displaystyle \forall l\geq 1,\ F_{0,1}=h^{(1)},\ F_{1,1}=\gamma
 ^{(1)}h,$ \ \par 
so ~(\ref{J8}) is true.\ \par 
\quad  Suppose that~(\ref{J8}) is true for  $j$  and let us see for
  $\displaystyle j+1.$ \ \par 
Suppose that  $\displaystyle l>j,$  we have\ \par 
\quad \quad \quad  $\displaystyle R^{j+1}(\gamma ^{l}h)=R(R^{j}(\gamma ^{l}h))=\gamma
 ^{l}R(F_{0,j})+l\gamma ^{l-1}(\gamma ^{(1)}F_{0,j}(z)+R(F_{1,j}))+\cdot
 \cdot \cdot +$ \ \par 
\quad \quad \quad \quad \quad \quad \quad  $\displaystyle +l(l-1)\cdot \cdot \cdot (l-k+1)\gamma ^{l-k}(\gamma
 ^{(1)}F_{k-1,j}(z)+R(F_{k,j}))+...+$ \ \par 
\quad \quad \quad \quad \quad \quad \quad \quad \quad  $\displaystyle +l(l-1)\cdot \cdot \cdot (l-j+1)\gamma ^{l-j}(\gamma
 ^{(1)}F_{j-1,j}(z)+R(F_{j,j}))+l(l-1)\cdot \cdot \cdot (l-j)\gamma
 ^{l-j-1}(\gamma ^{(1)}F_{j,j}).$ \ \par 
Hence we set\ \par 
\quad \quad \quad \quad \quad 	 $\displaystyle F_{0,j+1}:=R(F_{0,j})=h^{(j+1)},$ \ \par 
and\ \par 
\quad \quad \quad  $\displaystyle \ \forall k,\ 1\leq k\leq m,\ F_{k,j+1}:=\gamma
 ^{(1)}F_{k-1,j}(z)+R(F_{k,j})$ \ \par 
and the last one\ \par 
\quad \quad \quad  $\displaystyle F_{j+1,j+1}:=\gamma ^{(1)}F_{j,j}(z).$ \ \par 
If  $\displaystyle l=j\ :$  the formula~(\ref{J8}) read\ \par 
\quad \quad \quad  $\displaystyle R^{j}(\gamma ^{l}h)=\gamma ^{l}F_{0,j}(z)+l\gamma
 ^{l-1}F_{1,j}(z)+\cdot \cdot \cdot +l!F_{l,j}(z)$ \ \par 
hence we get\ \par 
\quad \quad \quad  $\displaystyle R^{j+1}(\gamma ^{l}h)=R(\gamma ^{l}F_{0,j}(z)+l\gamma
 ^{l-1}F_{1,j}(z)+\cdot \cdot \cdot +l!\gamma F_{l-1,j})+\ l!R(F_{l,j}).$
 \ \par 
So again\ \par 
\quad \quad \quad \quad \quad 	 $\displaystyle F_{0,j+1}=R(F_{0,j}),$ \ \par 
and\ \par 
\quad \quad \quad  $\displaystyle \ \forall k,\ 1\leq k\leq l-1,\ F_{k,j+1}=\gamma
 ^{(1)}F_{k-1,j}(z)+R(F_{k,j})$ \ \par 
but\ \par 
\quad \quad \quad  $\displaystyle F_{l,j+1}=\gamma ^{(1)}F_{l-1,j}$ \ \par 
and\ \par 
\quad \quad \quad \quad  $\displaystyle F_{l+1,j+1}=R(F_{l,j})$ \ \par 
which is formula~(\ref{J8}) with  $\displaystyle m=l=\min \ (j+1,l)$ .\ \par 
If  $\displaystyle l<j$  : by use of formula~(\ref{J8}) with
  $\displaystyle m=l=\min \ (j,l)$  we get\ \par 
\quad \quad \quad  $\displaystyle R^{j+1}(\gamma ^{l}h)=R(\gamma ^{l}F_{0,j}+l\gamma
 ^{l-1}F_{1,j}+\cdot \cdot \cdot +l!F_{l,j})$ \ \par 
hence again\ \par 
\quad \quad \quad \quad \quad 	 $\displaystyle F_{0,j+1}=R(F_{0,j}),$ \ \par 
and\ \par 
\quad \quad \quad  $\displaystyle \ \forall k,\ 1\leq k\leq l-1,\ F_{k,j+1}=\gamma
 ^{(1)}F_{k-1,j}(z)+R(F_{k,j})$ \ \par 
and\ \par 
\quad \quad \quad  $\displaystyle F_{l,j+1}:=\gamma ^{(1)}F_{l-1,j}+R(F_{l,j}).$ \ \par 
Clearly the  $\displaystyle F_{k,j+1}$  do not depend on  $\displaystyle
 l,$  for  $\displaystyle k\leq m,$  because the  $\displaystyle
 F_{k,j}$  do not.  $\displaystyle \hfill\blacksquare $ \ \par 

\begin{Lmm}
 ~\label{J11}We have, with  $\displaystyle \alpha _{m}$  constants
 independent of  $\gamma $  and of  $h$  :\par 
\quad \quad \quad \quad  	\begin{equation}  \forall k\leq j,\ F_{k,j}=\alpha _{k}R^{j}(\gamma
 ^{k}h)+\alpha _{k-1}\gamma R^{j}(\gamma ^{k-1}h)+\cdot \cdot
 \cdot +\alpha _{1}\gamma ^{k-1}R^{j}(\gamma h)+\alpha _{0}\gamma
 ^{k}h^{(j)}.\label{J9}\end{equation}
\end{Lmm}
\quad \quad  	Proof.\ \par 
To get  $\displaystyle F_{1,j}$  we take  $\displaystyle l=1$
  in~(\ref{J8}) so we get\ \par 
\quad \quad \quad \quad \quad 	 $\displaystyle R^{j}(\gamma h)=\gamma h^{(j)}+F_{1,j}\Rightarrow
 F_{1,j}=R^{j}(\gamma h)-\gamma h^{(j)}.$ \ \par 
So it is true for  $\displaystyle k=1$  and any  $\displaystyle
 l\geq 1$  because  $\displaystyle F_{1,j}$  is independent of
  $\displaystyle l.$  \ \par 
\quad  Suppose it is true up to  $k$  ; let us see for  $\displaystyle k+1.$ \ \par 
We choose  $\displaystyle l=k+1,\ j\geq l$  in~(\ref{J8}), we get\ \par 
\quad \quad \quad \quad \quad 	 $\displaystyle R^{j}(\gamma ^{k+1}h)=\gamma ^{k+1}h^{(j)}+(k+1)\gamma
 ^{k}F_{1,j}(z)+\cdot \cdot \cdot +(k+1)!\gamma F_{k,j}+(k+1)!F_{k+1,j},$
 \ \par 
hence\ \par 
\quad \quad \quad \quad \quad 	 $\displaystyle (k+1)!F_{k+1,j}=R^{j}(\gamma ^{k+1}h)-\gamma
 ^{k+1}h^{(j)}-(k+1)\gamma ^{k}F_{1,j}(z)-\cdot \cdot \cdot -(k+1)!\gamma
 F_{k,j}$ \ \par 
and assuming the decomposition~(\ref{J9}) for all the  $\displaystyle
 F_{m,j},\ m\leq k,$  we get that the formula is true for  $\displaystyle
 k+1.$   $\displaystyle \hfill\blacksquare $ \ \par 

\begin{Prps}
 ~\label{iP2} (Exclusion) We have, with  $\displaystyle m:=\min \ (l,j),$ \par 
\quad \quad \quad \quad \quad 	 $\forall j,l\in {\mathbb{N}},\ R^{j}(\gamma ^{l}h)=\sum_{q=0}^{m}{A_{q}\gamma
 ^{l-q}R^{j}(\gamma ^{q}h)},$ \par 
where the  $\displaystyle A_{q}$  are constants independent
 of  $\gamma $  and  $\displaystyle h.$ 
\end{Prps}
\quad  Proof.\ \par 
This is trivial if  $\displaystyle l\leq j,$  just take  $\displaystyle
 A_{q}=0$  for  $\displaystyle q<l$  and  $\displaystyle A_{l}=1.$
  Now take  $\displaystyle l>j+1.$ \ \par 
From lemma~\ref{J10} we get\ \par 
\quad \quad \quad  $\displaystyle \forall j,l\in {\mathbb{N}},\ R^{j}(\gamma ^{l}h)=\gamma
 ^{l}F_{0,j}(z)+l\gamma ^{l-1}F_{1,j}(z)+\cdot \cdot \cdot +l(l-1)\cdot
 \cdot \cdot (l-m+1)\gamma ^{l-m}F_{m,j}(z)$ \ \par 
and with lemma~\ref{J11} we replace the functions  $\displaystyle
 F_{k,j}$  to get what we want\ \par 
\quad \quad \quad  $\displaystyle \forall j,l\in {\mathbb{N}},\ R^{j}(\gamma ^{l}h)=\sum_{q=0}^{m}{A_{q}\gamma
 ^{l-q}R^{j}(\gamma ^{q}h)}.$   $\displaystyle \hfill\blacksquare $ \ \par 

\begin{Lmm}
 ~\label{iP0}(Domination) Let  $S=\lbrace a_{m}\rbrace _{m=1,...,N}$
  be an interpolating sequence in  $\displaystyle {\mathbb{B}}$
  for  $\displaystyle H_{s}^{p},$  of interpolating constant
  $\displaystyle C(S),$  and  $\displaystyle \lbrace \gamma _{a}\rbrace
 _{a\in S}$  its canonical dual sequence. Then\par 
\quad \quad \quad  $\displaystyle \forall l\in {\mathbb{N}},\ \forall j\leq s,\
 \forall h\in H_{s}^{p},\ \forall q\leq N,\ \exists H_{q}\in
 H_{s}^{p}::\forall a\in S,\ \ \left\vert{R^{j}(\gamma _{a}^{l}h)}\right\vert
 \leq \frac{1}{N}\sum_{q=1}^{N}{\left\vert{R^{j}(H_{q})}\right\vert }.$ \par 
So  $H_{q}$  depends on  $\displaystyle l,j$  and  $h,$  but
 {\bf not} on  $a$  and we have  $\displaystyle \ 1\leq q\leq
 N,\ {\left\Vert{H_{q}}\right\Vert}_{H_{s}^{p}}\leq C(S)^{l}{\left\Vert{h}\right\Vert}_{H_{s}^{p}}.$
 
\end{Lmm}
\quad  Proof.\ \par 
We have by definition of  $\gamma _{a}$ \ \par 
\quad \quad \quad  $\displaystyle \gamma _{a_{m}}(z):=\frac{1}{N}\sum_{q=1}^{N}{\theta
 ^{-qm}\beta (q,z)}\in {\mathcal{M}}_{s}^{p},\ {\left\Vert{\beta
 (q,\cdot )}\right\Vert}_{{\mathcal{M}}_{s}^{p}}\leq C(S).$ \ \par 
By lemma~\ref{iSH17} with  $\displaystyle Q_{l}(k,z):=\underbrace{\beta
 *\cdot \cdot \cdot *\beta (k,z)}_{l\ times}$  and  $\displaystyle
 \ {\left\Vert{Q_{l}(k,\cdot )}\right\Vert}_{{\mathcal{M}}_{s}^{p}}\leq
 C(S)^{l},$ \ \par 
\quad \quad \quad  $\displaystyle \gamma _{a_{m}}(z)^{l}=\widehat{Q_{l}(m,z)}=\frac{1}{N}\sum_{q=1}^{N}{\theta
 ^{-qm}Q_{l}(q,z)}$ \ \par 
so\ \par 
\quad \quad \quad  $\displaystyle \gamma _{a_{m}}(z)^{l}h=\frac{1}{N}\sum_{q=1}^{N}{\theta
 ^{-qm}Q_{l}(q,z)h}$ \ \par 
and\ \par 
\quad \quad \quad  $\displaystyle R^{j}(\gamma _{a_{m}}^{l}h)=R^{j}(\frac{1}{N}\sum_{q=1}^{N}{\theta
 ^{-qm}Q_{l}(q,z)h})=\frac{1}{N}\sum_{q=1}^{N}{\theta ^{-qm}R^{j}(Q_{l}(q,z)h)}.$
 \ \par 
So\ \par 
\quad \quad \quad  $\displaystyle \ \left\vert{R^{j}(\gamma _{a}^{l}h)}\right\vert
 \leq \frac{1}{N}\sum_{q=1}^{N}{\left\vert{R^{j}(Q_{l}(q,z)h)}\right\vert
 },$ \ \par 
hence setting\ \par 
\quad \quad \quad  $\displaystyle \forall z\in {\mathbb{B}},\ H_{q}(z):=Q_{l}(q,z)h(z)$ \ \par 
we have that  $H_{q}$  is independent of  $\displaystyle a\in S$  and\ \par 
\quad \quad \quad  $\displaystyle \ {\left\Vert{H_{q}}\right\Vert}_{H_{s}^{p}}\leq
 {\left\Vert{Q(q,\cdot )}\right\Vert}_{{\mathcal{M}}_{s}^{p}}{\left\Vert{h}\right\Vert}_{H_{s}^{p}}\leq
 C(S)^{l}{\left\Vert{h}\right\Vert}_{H_{s}^{p}}.$ \ \par 
This ends the proof of the lemma.  $\displaystyle \hfill\blacksquare $ \ \par 
\ \par 

\begin{Lmm}
 ~\label{iP4} (Inclusion) Let  $\displaystyle \gamma \in {\mathcal{M}}_{s}^{p}$
  and  $\displaystyle h\in H_{s}^{p},$  then there are constants
  $\displaystyle A_{q}$  such that\par 
\quad \quad \quad  $\displaystyle \forall j,l,\ R^{j}(\gamma ^{l})h=\sum_{q=0}^{j}{A_{j,q}R^{q}(\gamma
 ^{l}R^{j-q}(h))}.$ 
\end{Lmm}
\quad  Proof.\ \par 
By induction on  $j.$  For  $\displaystyle j=1$  we have  $\displaystyle
 R(\gamma ^{l}h)=R(\gamma ^{l})h+\gamma ^{l}R(h)$  hence\ \par 
\quad \quad \quad  $\displaystyle R(\gamma ^{l})h=R(\gamma ^{l}h)-\gamma ^{l}R(h),$ \ \par 
so it is true. Suppose it is true for any  $\displaystyle q<j$
  then we have\ \par 
\quad \quad \quad  $\displaystyle R^{j}(\gamma ^{l}h)=\sum_{q=0}^{j}{C_{j}^{q}R^{q}(\gamma
 ^{l})R^{j-q}(h)}=R^{j}(\gamma ^{l})h+\sum_{q=0}^{j-1}{C_{j}^{q}R^{q}(\gamma
 ^{l})R^{j-q}(h)}$ \ \par 
hence\ \par 
\quad \quad \quad  \begin{equation}  R^{j}(\gamma ^{l})h=R^{j}(\gamma ^{l}h)-\sum_{q=0}^{j-1}{C_{j}^{q}R^{q}(\gamma
 ^{l})R^{j-q}(h)}.\label{iP3}\end{equation}\ \par 
Now because  $\displaystyle q<j$  we have, with  $\displaystyle
 k:=R^{j-q}(h),$ \ \par 
\quad \quad \quad  $\displaystyle R^{q}(\gamma ^{l})k=\sum_{m=0}^{q}{A_{q,m}R^{m}(\gamma
 ^{l}R^{q-m}(k))}$ \ \par 
hence\ \par 
\quad \quad \quad  $\displaystyle R^{q}(\gamma ^{l})k=\sum_{m=0}^{q}{A_{q,m}R^{m}(\gamma
 ^{l}R^{q-m}(R^{j-q}(h)))}=\sum_{m=0}^{q}{A_{q,m}R^{m}(\gamma
 ^{l}R^{j-m}(h))}.$ \ \par 
Replacing in~(\ref{iP3}) we get the lemma.  $\displaystyle \hfill\blacksquare
 $ \ \par 
\ \par 

\bibliographystyle{C:/texlive/2012/texmf-dist/bibtex/bst/base/plain}

\begin{thebibliography}{10}

\bibitem{AglMCar02}
J.~Agler and J.~MacCarthy.
\newblock {\em Pick interpolation and {H}ilbert functions spaces.}, volume~44
  of {\em Graduate Studies in Mathematics}.
\newblock 2002.

\bibitem{DenAm78}
D.~Amar and E.~Amar.
\newblock Sur les suites d'interpolation en plusieurs variables.
\newblock {\em Pacific J. Math.}, 15:15--20, 1978.

\bibitem{AmarThesis77}
E.~Amar.
\newblock Suites d'interpolation dans le spectre d'une alg\`ebre
  d'op\'erateurs.
\newblock 1977.
\newblock Th\`ese.

\bibitem{amBerg78}
E.~Amar.
\newblock Suites d'interpolation pour les classes de {B}ergman de la boule et
  du polydisque de $\mathbb{C}^n$.
\newblock {\em Canadian J. Math.}, 30:711--737, 1978.

\bibitem{AmarExtInt06}
E.~Amar.
\newblock On linear extension for interpolating sequences.
\newblock {\em Studia Mathematica}, 186(3):251--265, 2008.

\bibitem{AmarCarlType07}
E.~Amar.
\newblock A {C}arleson type condition for interpolating sequences in the unit
  ball of $\mathbb{C}^n$.
\newblock {\em Publ. Mat.}, 53:481--488, 2009.

\bibitem{AmarWirtBoule07}
Eric Amar.
\newblock Interpolating sequences, {C}arleson mesures and {W}irtinger
  inequality.
\newblock {\em Annales Polonici Mathematici}, 94(1):79--87, 2008.

\bibitem{ArcoRochSaw06}
N.~Arcozzi, R.~Rochberg, and E.~Sawyer.
\newblock Carleson measures and interpolating sequences for {B}esov spaces on
  complex balls.
\newblock {\em Mem. Amer. Math. Soc.}, 2006.

\bibitem{ArcoRochSaw08}
N.~Arcozzi, R.~Rochberg, and E.~Sawyer.
\newblock Carleson measures for the {D}rury-{A}rveson {H}ardy space and other
  {B}esov-{S}obolev spaces on complex balls.
\newblock {\em Adv. Math.}, 218(4):1107--1180, 2008.

\bibitem{Bernard71}
A.~Bernard.
\newblock Alg\`ebre quotient d'alg\`ebre uniforme.
\newblock {\em C.R.A.S. de Paris}, 272:1101--1104, 1971.

\bibitem{PBeurling62}
P.~Beurling and L.~Carleson.
\newblock Research on interpolation problems.
\newblock Preprint, Uppsala, 1962.

\bibitem{CarlInt58}
L.~Carleson.
\newblock An interpolation problem for bounded analytic functions.
\newblock {\em Amer. J. Math.}, 1958.

\bibitem{CascOrte95}
C.~Cascante and Ortega J.
\newblock Carleson measures on spaces of {H}ardy-{S}obolev type.
\newblock {\em Canad. J. Math.}, 47:1177--1200, 1995.

\bibitem{Cobos86}
F.~Cobos.
\newblock Clarkson's inequalities for {S}obolev spaces.
\newblock {\em Math. Japon.}, 31(1):17--22, 1986.

\bibitem{Drury70}
S.~Drury.
\newblock Sur les ensembles de {S}idon.
\newblock {\em C. R. Acad. Sci. Paris}, 271:A162--A163, 1970.

\bibitem{FolStein74}
G.~B. Folland and E.~M. Stein.
\newblock Estimates for the $\bar \partial _{b}$ complex and analysis on the
  {H}eisenberg group.
\newblock {\em Comm. Pure Appl. Math.}, 27:429--522, 1974.

\bibitem{HormPSH67}
L.~Hormander.
\newblock A {$L^p$} estimates for (pluri-) subharmonic functions.
\newblock {\em Math. Scand.}, 20:65--78, 1967.

\bibitem{Ligocka87}
E.~Ligocka.
\newblock Estimates in {S}obolev norms for harmonic and holomorphic functions
  and interpolation between {S}obolev and {H}\"older spaces of harmonic
  functions.
\newblock {\em Studia Math.}, 86:255--271, 1987.

\bibitem{OrtFab97}
J.~Ortega and J.~Fabrega.
\newblock Holomorphic {T}riebel {L}izorkin spaces.
\newblock {\em Journal of functional analysis}, 151:177--212, 1997.

\bibitem{Romanovskii05}
N.~N. Romanovski.
\newblock Integral representations and embedding theorems for functions defined
  on the {H}eisenberg groups $\mathbb {H}^n$.
\newblock {\em St. Petersburg Math. J.}, 16(2):349--375, 2005.

\bibitem{RudinBall81}
W.~Rudin.
\newblock Function theory in the unit ball of $ {C}^n$.
\newblock {\em Grundenlehren}, 1981.

\bibitem{Thomas87}
P.~J. Thomas.
\newblock Hardy space interpolation in the unit ball.
\newblock {\em Indagationes Mathematicae}, 90(3):325--351, 1987.

\bibitem{Varo71}
N.~Varopoulos.
\newblock Sur la r\'eunion de deux ensembles d'interpolation d'une alg\`ebre
  uniforme.
\newblock {\em C.R.A.S. Paris}, 272:950--952, 1971.

\bibitem{Varo72}
N.~Varopoulos.
\newblock Sur un probl\`eme d'interpolation.
\newblock {\em C.R.A.S. Paris}, 274:1539--1542, 1972.

\bibitem{VolbWick12}
A.~Volberg and B.~Wick.
\newblock Bergman-type singular integral operators and the characterization of
  {Carleson measures for Besov-Sobolev} spaces in the complex ball.
\newblock {\em Amer. J. Math.}, 134(4):949--992, 2012.

\end{thebibliography}

\end{document}